\documentclass[final,hidelinks,onefignum,onetabnum]{siamart251104}
\usepackage[left=3cm, right=2.5cm, top=2.5cm, bottom=2.5cm]{geometry} 
\usepackage{lipsum}
\usepackage{algorithmic}
\usepackage{fancyhdr}
\usepackage{graphicx}
\usepackage{subfigure}
\usepackage{amssymb}
\usepackage{amsmath}
\usepackage{amsfonts}
\usepackage{theorem}
\usepackage{mathrsfs}
\usepackage{mathtools}
\usepackage{bm}
\usepackage{color}
\usepackage{setspace}
\usepackage{exscale}
\usepackage{relsize}
\usepackage{tikz}
\usetikzlibrary{positioning}
\usepackage{rotating}
\usepackage{enumerate}
\usepackage{hyperref}
\usepackage{epstopdf}
\usepackage{float}
\ifpdf
  \DeclareGraphicsExtensions{.eps,.pdf,.png,.jpg}
\else
  \DeclareGraphicsExtensions{.eps}
\fi
\usepackage{cite}

\usepackage{nicefrac}
\usepackage{setspace}
\usepackage{hyperref}

\usepackage{booktabs,multirow} 
\usepackage{array} 
\usepackage{paralist} 
\usepackage{verbatim} 
\usepackage{subfigure} 

\newsiamremark{remark}{Remark}
\newsiamremark{hypothesis}{Hypothesis}
\crefname{hypothesis}{Hypothesis}{Hypotheses}
\newsiamthm{claim}{Claim}
\newsiamremark{fact}{Fact}
\crefname{fact}{Fact}{Facts}

\renewcommand{\theequation}{\arabic{section}.\arabic{equation}}
\renewcommand{\thefigure}{\arabic{section}.\arabic{figure}}
\renewcommand{\thetable}{\arabic{section}.\arabic{table}}

\allowdisplaybreaks[1]

\numberwithin{equation}{section}
\numberwithin{figure}{section}
\numberwithin{table}{section}

\newcommand\eref[1]{(\ref{#1})}

\newcommand*\xbar[1]{%
  \hbox{%
    \vbox{%
      \hrule height 0.5pt 
      \kern0.4ex
      \hbox{%
        \kern-0.05em
        \ensuremath{#1}%
        \kern-0.00em
      }%
    }%
  }%
}

\newcommand{\bmF}{\bm{\mathcal{F}}}
\newcommand{\bmG}{\bm{\mathcal{G}}}

\newcommand{\bmK}{\bm{\mathcal{K}}}
\newcommand{\bmL}{\bm{\mathcal{L}}}
\newcommand{\mF}{\bm{F}}
\newcommand{\mK}{\bm{K}}

\newcommand{\mG}{\bm{G}}
\newcommand{\mL}{\bm{L}}

\newcommand{\mU}{\bm{U}}

\newcommand{\mV}{\bm{V}}

\newcommand{\mW}{\bm{W}}
\newcommand{\mo}{\bm{0}}

\newcommand{\dx}{\Delta x}
\newcommand{\dy}{\Delta y}

\newcommand{\bnabla}{\bm\nabla}
\newcommand{\hf}{{\frac{1}{2}}}

\newcommand{\jph}{{j+\frac{1}{2}}}
\newcommand{\jmh}{{j-\frac{1}{2}}}
\newcommand{\kph}{{k+\frac{1}{2}}}
\newcommand{\kmh}{{k-\frac{1}{2}}}

\def\softd{{\leavevmode\setbox1=\hbox{d}%
          \hbox to 1.05\wd1{d\kern-0.4ex{\char039}\hss}}}

\title{A Locally Divergence-Free Local Characteristic Decomposition Based Path-Conservative Central-Upwind Scheme for Ideal
Magnetohydrodynamics\thanks{Submitted to the editors DATE.
\funding{The work of S. Chu was funded by the DFG--SPP 2183: Eigenschaftsgeregelte Umformprozesse with the Project(s) HE5386/19-2,19-3 Entwicklung eines flexiblen isothermen Reckschmiedeprozesses f\"ur die eigenschaftsgeregelte Herstellung von Turbinenschaufeln aus Hochtemperaturwerkstoffen (424334423) and by the Deutsche Forschungsgemeinschaft (DFG, German Research Foundation)--SPP 2410 Hyperbolic Balance Laws in Fluid Mechanics: Complexity, Scales, Randomness (CoScaRa) within the Project(s) HE5386/27-1 (Zuf\"allige kompressible Euler Gleichungen: Numerik und ihre Analysis, 525853336).  The work of A. Kurganov was supported in
part by NSFC grants 12171226 and W2431004. The work of M. Luk\'{a}\v{c}ov\'{a}-Medvi{\softd}ov\'{a} was supported by the Deutsche Forschungsgemeinschaft (DFG, German Research
Foundation)--SPP 2410  within the Projects
LU1470/10-1  (Random Euler Equations and their Numerical Analysis, 525853336)  and LU1470/9-1 (An Active Flux method for the Euler equations, 525800857), and partially by the Gutenberg Research College and by the DFG - project number 233630050 - TRR 146. She  is grateful to  the  Mainz Institute of Multiscale Modelling  for supporting her research.}}}

\author{Shaoshuai Chu\thanks{Department of Mathematics, RWTH Aachen University, Aachen, 52056, Germany (\email{chu@igpm.rwth-aachen.de}).}
\and Alexander Kurganov\thanks{Department of Mathematics and Shenzhen International Center for Mathematics, Southern University of Science
and Technology, Shenzhen, 518055, China (\email{alexander@sustech.edu.cn}).}
\and M\'{a}ria Luk\'{a}\v{c}ov\'{a}-Medvid'ov\'{a}\thanks{Institute of Mathematics, Johannes Gutenberg University Mainz, Staudingerweg 9,
55128 Mainz, Germany (\email{lukacova@uni-mainz.de}).}
\and Mingye Na\thanks{Department of Mathematics, Southern University of Science and Technology, Shenzhen, 518055, China and Institute
of Mathematics, Johannes Gutenberg University Mainz, Staudingerweg 9, 55128 Mainz, Germany (\email{namingye@uni-mainz.de}).}}

\usepackage{amsopn}

\begin{document}

\maketitle

\begin{abstract}
We introduce a locally divergence-free local characteristic decomposition based path-conservative central-upwind (LCD-PCCU) scheme for ideal
magnetohydrodynamics (MHD) equations. The proposed method is a low-dissipation extension of the recently proposed locally divergence-free
PCCU scheme. To reduce the numerical dissipation, we incorporate the LCD into the PCCU framework. The resulting LCD-PCCU method enhances the
resolution of numerical solutions as demonstrated through a series of benchmark tests.
\end{abstract}

\begin{keywords}
Ideal magnetohydrodynamics, divergence-free constraints, local characteristic decomposition, \\
 path-conservative central-upwind scheme.
\end{keywords}

\begin{MSCcodes}
65M08, 76W05, 76M12, 35L65
\end{MSCcodes}

\section{Introduction}
This paper focuses on the development of a novel and low-dissipation numerical method for ideal magnetohydrodynamics (MHD) equations, which
play a central role in modeling a wide range of physical phenomena in astrophysics, plasma physics, space physics, and engineering. These
models describe the dynamics of electrically conducting fluids interacting with magnetic fields and consist of hyperbolic systems of partial
differential equations (PDEs) that couple fluid flow with electromagnetic effects. A crucial property of these models is a constraint on the
magnetic field, which has to remain divergence-free if it is divergence-free initially. Numerically, however, this condition is nontrivial
to maintain, and improperly handling the divergence-free constraint at the discrete level can lead to numerical instabilities or the
development of nonphysical structures in the solution; see, e.g., \cite{BS1999,BB1980,LS2005,Toth2000}.

Over the past decades, numerous approaches have been developed to address the divergence-free constraint. Among them are the projection
method (see, e.g., \cite{BB1980}), the constrained transport (CT) method (see, e.g.,
\cite{EH1988,BS1999,DeVore1991,GS2005,XBD2016,CRT2014,HRT2013,MT2012,Rossmanith2006}), locally divergence-free discontinuous Galerkin
\cite{LS2005,YXL2013} and finite-volume \cite{Chertock24,CKRZ2024} methods (these methods maintain zero divergence within each computational
cell), and globally divergence-free high-order finite-volume and discontinuous Galerkin methods; see, e.g.,
\cite{Balsara2009,BKC2021,DBTF2019,FLX2018,LX2012,LXY2011}. 

Alternatively, instead of enforcing the divergence-free constraint explicitly, one can reduce divergence errors through the inclusion of
additional Godunov-Powell terms. This leads to the so-called eight-wave formulation of the ideal MHD equations; see, e.g.,
\cite{Godunov2025,PRLGD1999,PRMGD1995,Powell1997}. This formulation introduces nonconservative source terms proportional to the divergence
of the magnetic field. Although these terms vanish analytically, they help to control numerical divergence errors by advecting them with the
flow and preventing their accumulation. Moreover, the modified system gains important properties such as Galilean invariance and entropy
symmetrizability, making it well-suited for the development of entropy-stable schemes; see, e.g.,
\cite{CK2016,DGWW2018,LSZ2018,PRLGD1999,PRMGD1995,Powell1997}.

In \cite{Chertock24}, the Godunov-Powell modification of the ideal MHD and shallow water MHD equations was utilized to develop a locally
divergence-free second-order semi-discrete path-conservative central-upwind (PCCU) scheme, which was later extended to the magnetic rotating
shallow water model in \cite{CKRZ2024}. In \cite{Chertock24}, the studied systems were augmented by evolution equations for the spatial
derivatives of the magnetic field components, and the resulting systems were numerically solved by a PCCU scheme, which was originally
developed in \cite{CKM} as a ``black-box'' solver for general nonconservative hyperbolic systems. We stress that the PCCU schemes are, like
any central and central-upwind (CU) schemes, Riemann-problem-solver-free, and at the same time, they are designed to handle the
nonconservative product terms across cell interfaces in a stable manner.

Although the PCCU scheme for the ideal MHD system \cite{Chertock24} is quite accurate, efficient, and robust, its resolution can be further
improved by reducing the amount of numerical dissipation. This can be done with the help of the local characteristic decomposition (LCD)
based PCCU (LCD-PCCU) scheme, which was recently introduced in \cite{CHK25} as an extension of the LCD-based CU (LCD-CU) scheme proposed in
\cite{CCHKL_22} for hyperbolic systems of conservation laws. Compared with the CU and PCCU schemes, the LCD-CU and LCD-PCCU schemes achieve
higher resolution by aligning the numerical flux computation with the characteristic structure of the system; see
\cite{CHK25,CCHKL_22,CK2023}. In this paper, we develop the LCD-PCCU scheme for the ideal MHD equations and test it on a number of numerical
experiments, which confirm that the proposed scheme achieves high resolution, while being robust and positivity-preserving, and effectively
maintaining divergence control. The obtained numerical results also demonstrate that the new scheme outperforms the PCCU scheme from
\cite{Chertock24}.

The rest of the paper is organized as follows. In \S\ref{sec2}, we present the Godunov-Powell modification of the ideal MHD equations and
its augmented form. In \S\ref{sec3}, we apply the two-dimensional (2-D) LCD-PCCU scheme for the studied MHD system. Finally, in
\S\ref{sec4}, we present results of several numerical experiments.

\section{Ideal MHD Equations}\label{sec2}
The ideal MHD equations read as
\begin{equation}
\begin{aligned}
&\rho_t+\bnabla\!\cdot\!(\rho\bm u)=0,\\
&(\rho\bm u)_t+\bnabla\!\cdot\!\Big[\rho\bm u\bm u^\top+\Big(p+\hf|\bm b|^2\Big)I-\bm b\bm b^\top\Big]=\bm0,\\
&{\cal E}_t+\bnabla\!\cdot\!\Big[\Big({\cal E}+p+\hf|\bm b|^2\Big)\bm u-\bm b(\bm u\!\cdot\!\bm b)\Big]=0,\\
&\bm b_t-\bnabla\!\times\!(\bm u\!\times\!\bm b)=\bm0,
\end{aligned}
\label{1.1}
\end{equation}
where $t$ is time, $\rho$ is the density, $p$ is the pressure, $\bm u=(u,v,w)^\top$ is the fluid velocity, $\bm b=(b_1,b_2,b_3)^\top$ is the
magnetic field, and $\cal E$ is the total energy. Additionally, $I$ is the identity matrix and $\gamma$ is the ratio of specific heats. The
system \eref{1.1} is completed through the equation of state (EOS)
\begin{equation}
{\cal E}=\frac{p}{\gamma-1}+\frac{\rho}{2}|\bm u|^2+\hf|\bm b|^2.
\label{1.2}
\end{equation}
It is easy to show that provided that the magnetic field is initially divergence-free, then the magnetic field satisfies
\begin{equation}
\bnabla\!\cdot\bm b=0.
\label{div}
\end{equation}

In this paper, we will develop a new numerical method for the Godunov-Powell modified ideal MHD equations, which read as
\begin{equation}
\begin{aligned}
&\rho_t+\bnabla\!\cdot\!(\rho\bm u)=0,\\
&(\rho\bm u)_t+\bnabla\!\cdot\!\Big[\rho\bm u\bm u^\top+\Big(p+\hf|\bm b|^2\Big)I-\bm b\bm b^\top\Big]=-\bm b(\bnabla\!\cdot\bm b),\\
&{\cal E}_t+\bnabla\!\cdot\!\Big[\Big({\cal E}+p+\hf|\bm b|^2\Big)\bm u-\bm b(\bm u\!\cdot\!\bm b)\Big]=
-(\bm u\!\cdot\!\bm b)(\bnabla\!\cdot\bm b),\\
&\bm b_t-\bnabla\!\times\!(\bm u\!\times\!\bm b)=-\bm u(\nabla\!\cdot\bm b),
\end{aligned}
\label{1.3}
\end{equation}
which is completed through the EOS \eref{1.2}. We stress that the system \eref{1.3}, \eref{1.2}---unlike the original system
\eref{1.1}--\eref{1.2}---has a complete set of eight eigenvalues with eight corresponding eigenvectors: This allows for an LCD and thus for
designing an LCD-PCCU scheme.

As in \cite{Chertock24}, we restrict our attention to the 2-D case, where all the quantities of interest depend on the spatial variables $x$
and $y$ and time $t$ only. In this case, the divergence-free condition \eref{div} reads as $(b_1)_x+(b_2)_y=0$, and we augment the
Godunov-Powell modified ideal MHD system \eref{1.3}, \eref{1.2} by adding the equations for the auxiliary variables $A:=(b_1)_x$ and
$B:=(b_2)_y$:
\begin{equation}
\begin{aligned}
&A_t+(uA-b_2u_y)_x+(vA+b_1v_x)_y=0,\\
&B_t+(uB+b_2u_y)_x+(vB-b_1v_x)_y=0,
\end{aligned}
\label{1.4}
\end{equation}
which are obtained by differentiating  the $b_1$- and $b_2$-equations in \eref{1.1}.

One can write the system \eref{1.3}--\eref{1.4}, \eref{1.2} in the following vector form:
\begin{align}
&\mU_t+\mF(\mU)_x+\mG(\mU)_y=Q^x(\mU)\mU_x+Q^y(\mU)\mU_y,\label{2.5}\\
&\widetilde{\mU}_t+\widetilde{\mF}(\mW)_x+\widetilde{\mG}(\mW)_y=\bm0,\label{2.6}
\end{align}
where
\begin{equation}
\begin{aligned}
\mU&=(\rho,\rho u,\rho v,\rho w,b_1,b_2,b_3,{\cal E})^\top,\quad\widetilde{\mU}=(A,B)^\top,\quad\mW=(\mU^\top,\widetilde{\mU}^\top)^\top,\\
\mF(\mU)&=\Big(\rho u,\rho u^2+p+\hf|\bm b|^2-b_1^2,\rho uv-b_1b_2,\rho uw-b_1b_3,0,ub_2-vb_1,\\
&\hspace*{0.8cm}ub_3-wb_1,\Big({\cal E}+p+\hf|\bm b|^2\Big)u-(\bm u\!\cdot\!\bm b)b_1\Big)^\top,\\
\mG(\mU)&=\Big(\rho v,\rho uv-b_1b_2,\rho v^2+p+\hf|\bm b|^2-b_2^2,\rho vw-b_2b_3,vb_1-ub_2,0,\\
&\hspace*{0.8cm}vb_3-wb_2,\Big({\cal E}+p+\hf|\bm b|^2\Big)v-(\bm u\!\cdot\!\bm b)b_2\Big)^\top,\\
Q^x(\mU)&=\bm q\mathbf e_5^\top,~~Q^y(\mU)=\bm q\mathbf e_6^\top,\quad\bm q:=-(0,b_1,b_2,b_3,u,v,w,\bm u\!\cdot\!\bm b)^\top,\\
\widetilde{\mF}(\mW)&=\Big(uA-b_2u_y,uB+b_2u_y\Big)^\top,\quad\widetilde{\mG}(\mW)=\Big(vA+b_1v_x,vB-b_1v_x\Big)^\top,
\end{aligned}
\label{3.2a}
\end{equation}
and $\mathbf e_5$ and $\mathbf e_6$ are the fifth and sixth unit vectors in $\mathbb R^8$, respectively. 

Note that for smooth solutions, the system \eref{1.3}, \eref{1.2} can be rewritten in an equivalent quasi-linear form:
\begin{equation}
\mU_t+C^x(\mU)\mU_x+C^y(\mU)\mU_y=\bm0,
\label{2.7}
\end{equation}
where the matrices $C^x(\mU)$ and $C^y(\mU)$ are specified in Appendix \ref{appA1}. Furthermore, one can switch to the primitive variables
$\bm V=(\rho,u,v,w,p,b_1,b_2,b_3)^\top$ and rewrite the system \eref{2.7} in a different quasi-linear form:
\begin{equation}
\mV_t+D^x(\mV)\mV_x+D^y(\mV)\mV_y=\bm0,
\label{2.8}
\end{equation}
where the matrices $D^x(\mV)$ and $D^y(\mV)$ are specified in Appendix \ref{appA2}. In the following, we will use the forms \eref{2.7} and
\eref{2.8} to design the LCD-CU numerical fluxes and to perform a piecewise linear reconstruction, respectively.

\section{2-D Flux Globalization Based LCD-PCCU Scheme}\label{sec3}
In this section, we apply the 2-D flux globalization based LCD-PCCU scheme from \cite{CHK25} to the studied augmented ideal MHD system
\eref{2.5}--\eref{3.2a}. To this end, we first rewrite the system \eref{2.5} in the following quasi-conservative form:
\begin{equation}
\bm U_t+\bm K(\bm U)_x+\bm L(\bm U)_y=\mo,\quad\bm K(\bm U)=\bm F(\bm U)-\bm I^x(\bm U),\quad\bm L(\bm U)=\bm G(\bm U)-\bm I^y(\bm U),
\label{3.2}
\end{equation}
where
\begin{equation*}
\bm I^x(\bm U):=\int\limits^x_{\widehat x}\left[Q^x(\bm U)\bm U_\xi(\xi,y,t)\right]{\rm d}\xi,\quad
\bm I^y(\bm U):=\int\limits^y_{\widehat y}\left[Q^y(\bm U)\bm U_\eta(x,\eta,t)\right]{\rm d}\eta,
\end{equation*}
with $\widehat x$ and $\widehat y$ being arbitrary numbers.

We cover the computational domain with uniform cells $C_{j,k}:=[x_\jmh,x_\jph]\times[y_\kmh,y_\kph]$ centered at
$(x_j,y_k)=\big((x_\jmh+x_\jph)/2$, $(y_\kph+y_\kmh)/2\big)$ with $x_\jph-x_\jmh\equiv\dx$ and $y_\kph-y_\kmh\equiv\dy$ for $j=1,\dots,N_x$,
$k=1,\dots,N_y$, and assume that the computed cell averages of $\mW$ over the corresponding cells $C_{j,k}$,
$$
\xbar\mW_{j,k}(t):\approx\frac{1}{\dx\dy}\int\limits_{C_{j,k}}\mW(x,y,t)\,{\rm d}x{\rm d}y,
$$
are available at a certain time level $t\ge0$. Note that $\xbar\mW_{j,k}$ as well as many of the indexed quantities introduced below are
time-dependent, but from here on, we suppress this dependence for the sake of brevity.

According to \cite{CHK25}, the cell averages $\,\xbar\mW_{j,k}=\big(\,\xbar\mU^\top_{j,k},\,\xbar{\widetilde{\mU}}^\top_{j,k}\big)^\top$ are
evolved in time by numerically solving the following system of ODEs:
\begin{equation}
\begin{aligned}
&\frac{{\rm d}\xbar\mU_{j,k}}{{\rm d}t}=-\frac{\bmK^{\rm LCD}_{\jph,k}-\bmK^{\rm LCD}_{\jmh,k}}{\dx}-
\frac{\bmL^{\rm LCD}_{j,\kph}-\bmL^{\rm LCD}_{j,\kmh}}{\dy},\\
&\frac{{\rm d}\xbar{\widetilde{\mU}}_{j,k}}{{\rm d}t}=-\frac{\widetilde{\bmF}_{\jph,k}-\widetilde{\bmF}_{\jmh,k}}{\dx}-
\frac{\widetilde{\bmG}_{j,\kph}-\widetilde{\bmG}_{j,\kmh}}{\dy},
\end{aligned}
\label{3.3}
\end{equation}
where the numerical fluxes $\widetilde{\bmF}_{\jph,k}$ and $\widetilde{\bmG}_{j,\kph}$ are evaluated by the CU scheme from \cite{KNP} 
\begin{equation}
\begin{aligned}
&\widetilde{\bmF}_{\jph,k}=\frac{s^+_{\jph,k}\widetilde{\mF}(\mW^{\rm E}_{j,k})-s^-_{\jph,k}\widetilde{\mF}(\mW^{\rm W}_{j+1,k})}
{s^+_{\jph,k}-s^-_{\jph,k}}+\frac{s^+_{\jph,k}s^-_{\jph,k}}{s^+_{\jph,k}-s^-_{\jph,k}}
\Big(\widetilde{\mU}^{\rm W}_{j+1,k}-\widetilde{\mU}^{\rm E}_{j,k}\Big),\\
&\widetilde{\bmG}_{j,\kph}=\frac{s^+_{j,\kph}\widetilde{\mG}(\mW^{\rm N}_{j,k})-s^-_{j,\kph}\widetilde{\mG}(\mW^{\rm S}_{j,k+1})}
{s^+_{j,\kph}-s^-_{j,\kph}}+\frac{s^+_{j,\kph}s^-_{j,\kph}}{s^+_{j,\kph}-s^-_{j,\kph}}
\Big(\widetilde{\mU}^{\rm S}_{j,k+1}-\widetilde{\mU}^{\rm N}_{j,k}\Big),
\end{aligned}
\label{3.3a}
\end{equation}
the global numerical fluxes  $\bmK^{\rm LCD}_{\jph,k}$ and $\bmL^{\rm LCD}_{j,\kph}$ are given by
\begin{equation}
\begin{aligned}
\bmK^{\rm LCD}_{\jph,k}&=R^x_{\jph,k}P^{\rm LCD}_{\jph,k}\big(R^x_{\jph,k}\big)^{-1}\mK^{\rm E}_{j,k}+
R^x_{\jph,k}M^{\rm LCD}_{\jph,k}\big(R^x_{\jph,k}\big)^{-1}\mK^{\rm W}_{j+1,k}\\
&+R^x_{\jph,k}Q^{\rm LCD}_{\jph,k}\big(R^x_{\jph,k}\big)^{-1}\big(\mU^{\rm W}_{j+1,k}-\mU^{\rm E}_{j,k}\big),\\[0.8ex]
\bmL^{\rm LCD}_{j,\kph}&=R^y_{j,\kph} P^{\rm LCD}_{j,\kph} \big(R^y_{j,\kph}\big)^{-1} \mL^{\rm N}_{j,k}+
R^y_{j,\kph}M^{\rm LCD}_{j,\kph}\big(R^y_{j,\kph}\big)^{-1}\mL^{\rm S}_{j,k+1}\\
&+R^y_{j,\kph}Q^{\rm LCD}_{j,\kph}\big(R^y_{j,\kph}\big)^{-1}\big(\mU^{\rm S}_{j,k+1}-\mU^{\rm N}_{j,k}\big),
\end{aligned}
\label{3.4}
\end{equation}
and the global fluxes $\bm K^{\rm E,W}_{j,k}$ and $\bm{{\cal L}}^{\rm N,S}_{j,k}$ in \eref{3.4} are obtained using \eref{3.2}:
\begin{equation}
\bm K^{\rm E,W}_{j,k}=\bm F\big(\bm U^{\rm E,W}_{j,k}\big)-(\bm I^x)^{\rm E,W}_{j,k},\quad
\bm L^{\rm N,S}_{j,k}=\bm G\big(\bm U^{\rm N,S}_{j,k}\big)-(\bm I^y)^{\rm N,S}_{j,k}.
\label{3.5a}
\end{equation}
In \eref{3.3a}--\eref{3.5a}, the following quantities have been used.

\smallskip
\noindent
$\bullet$ $\mU^{\rm E,W,N,S}_{j,k}$ are the point values of $\mU$ at midpoints of the cell interfaces of $C_{j,k}$. They are obtained using
a piecewise linear reconstruction applied to the primitive variables $\bm V$ using the corresponding LCD. To this end, we first compute
$$
\begin{aligned}
&u_{j,k}=\frac{(\xbar{\rho u})_{j,k}}{\xbar\rho_{j,k}},\quad v_{j,k}=\frac{(\xbar{\rho v})_{j,k}}{\xbar\rho_{j,k}},\quad
w_{j,k}=\frac{(\xbar{\rho w})_{j,k}}{\xbar\rho_{j,k}},\\
&p_{j,k}=(\gamma-1)\Big[\,\xbar{\cal E}_{j,k}-\frac{\,\xbar\rho_{j,k}}{2}(u_{j,k}^2+v_{j,k}^2+w_{j,k}^2)-
\hf\big((\xbar b_1)_{j,k}^2+(\xbar b_2)_{j,k}^2+(\xbar b_3)_{j,k}^2\big)\Big],
\end{aligned}
$$
where the latter expression has been obtained using the EOS \eref{1.2}.

We then switch to the local characteristic variables $\bm\Gamma$ at the midpoints of each of the cell interfaces $(x_\jph,y_k)$ and
$(x_j,y_\kph)$:
$$
\begin{aligned}
\bm\Gamma^x_{\ell,k}&=\big(T^x_{\jph,k}\big)^{-1}\mV_{\ell,k},\quad\ell=j-1,j,j+1,j+2,\\
\bm\Gamma^y_{j,m}&=\big(T^y_{j,\kph}\big)^{-1}\mV_{j,m},\quad m=k-1,k,k+1,k+2,
\end{aligned}
$$
where the matrices $T^x_{\jph,k}$ and $T^y_{j,\kph}$ are obtained using the LCD for the primitive system \eref{2.8} and they are given in
Appendix \ref{appA2}.

Next, we perform generalized minmod reconstructions in the $x$- and $y$-directions to evaluate the slopes
\begin{equation}
\begin{aligned}
(\bm\Gamma^x_x)_{j,k}&={\rm minmod}\left(\theta\,\frac{\bm\Gamma^x_{j+1,k}-\bm\Gamma^x_{j,k}}{\dx},\,
\frac{\bm\Gamma^x_{j+1,k}-\bm\Gamma^x_{j-1,k}}{2\dx},\,\theta\,\frac{\bm\Gamma^x_{j,k}-\bm\Gamma^x_{j-1,k}}{\dx}\right),\\
(\bm\Gamma^x_x)_{j+1,k}&={\rm minmod}\left(\theta\,\frac{\bm\Gamma^x_{j+2,k}-\bm\Gamma^x_{j+1,k}}{\dx},\,
\frac{\bm\Gamma^x_{j+2,k}-\bm\Gamma^x_{j,k}}{2\dx},\,\theta\,\frac{\bm\Gamma^x_{j+1,k}-\bm\Gamma^x_{j,k}}{\dx}\right),
\end{aligned}
\label{3.5}
\end{equation}
and
\begin{equation}
\begin{aligned}
(\bm\Gamma^y_y)_{j,k}&={\rm minmod}\left(\theta\,\frac{\bm\Gamma^y_{j,k+1}-\bm\Gamma^y_{j,k}}{\dy},\,
\frac{\bm\Gamma^y_{j,k+1}-\bm\Gamma^y_{j,k-1}}{2\dy},\,\theta\,\frac{\bm\Gamma^y_{j,k}-\bm\Gamma^y_{j,k-1}}{\dy}\right),\\
(\bm\Gamma^y_y)_{j,k+1}&={\rm minmod}\left(\theta\,\frac{\bm\Gamma^y_{j,k+2}-\bm\Gamma^y_{j,k+1}}{\dy},\,
\frac{\bm\Gamma^y_{j,k+2}-\bm\Gamma^y_{j,k}}{2\dy},\,\theta\,\frac{\bm\Gamma^y_{j,k+1}-\bm\Gamma^y_{j,k}}{\dy}\right),
\end{aligned}
\label{3.6}
\end{equation}
respectively. In \eref{3.5} and \eref{3.6}, the minmod function is defined as 
$$
{\rm minmod}(z_1,z_2,\ldots):=\begin{cases}\min_j\{z_j\}&\mbox{if}~z_j>0~~\forall j,\\\max_j\{z_j\}&\mbox{if}~z_j<0~~\forall j,\\
0&\mbox{otherwise},\end{cases}
$$
and it is applied in a component-wise manner. The parameter $\theta\in[1,2]$ is used to control the non-oscillatory property of the
resulting scheme: larger $\theta$ typically leads to a sharper, but more oscillatory computed solution.

We then evaluate the corresponding one-sided point values:
$$
\begin{aligned}
&(\bm\Gamma^x_{j,k})^{\rm E}=\bm\Gamma^x_{j,k}+\frac{\dx}{2}(\bm\Gamma^x_x)_{j,k},\quad
(\bm\Gamma^x_{j+1,k})^{\rm W}=\bm\Gamma^x_{j+1,k}-\frac{\dx}{2}(\bm\Gamma^x_x)_{j+1,k},\\
&(\bm\Gamma^y_{j,k})^{\rm N}=\bm\Gamma^y_{j,k}+\frac{\dy}{2}(\bm\Gamma^y_y)_{j,k},\quad
(\bm\Gamma^y_{j,k+1})^{\rm S}=\bm\Gamma^y_{j,k+1}-\frac{\dy}{2}(\bm\Gamma^y_y)_{j,k+1},
\end{aligned}
$$
switch back to the primitive variables:
\begin{equation*}
\bm V^{\rm E}_{j,k}=T^x_{\jph,k}\big(\bm\Gamma^x_{j,k}\big)^{\rm E},~~
\bm V^{\rm W}_{j+1,k}=T^x_{\jph,k}\big(\bm\Gamma^x_{j+1,k}\big)^{\rm W},~~
\bm V^{\rm N}_{j,k}=T^y_{j,\kph}\big(\bm\Gamma^y_{j,k}\big)^{\rm N},~~
\bm V^{\rm S}_{j,k+1}=T^y_{j,\kph}\big(\bm\Gamma^y_{j,k+1}\big)^{\rm S},
\end{equation*}
and then transform $\mV_{j,k}^{\rm E,W,N,S}$ into $\mU_{j,k}^{\rm E,W,N,S}$, which are non-oscillatory, but they do not necessarily satisfy
the local divergence-free requirement, which can be written as 
\begin{equation}
(\bnabla\!\cdot\!\bm b)_{j,k}:=
\frac{(b_1)^{\rm E}_{j,k}-(b_1)^{\rm W}_{j,k}}{\dx}+\frac{(b_2)^{\rm N}_{j,k}-(b_2)^{\rm S}_{j,k}}{\dy}\equiv0,~~\forall j,k.
\label{4.1}
\end{equation}

We thus need to correct the point values $(b_1)^{\rm E,W}_{j,k}$ and $(b_2)^{\rm N,S}_{j,k}$. To this end, we proceed similarly to
\cite[\S2.2.1]{Chertock24} by setting the slopes
\begin{equation}
((b_1)_x)_{j,k}=\sigma_{j,k}\xbar A_{j,k}\quad\mbox{and}\quad
((b_2)_y)_{j,k}=\sigma_{j,k}\xbar B_{j,k},
\label{corr1}
\end{equation}
where
\begin{equation}
\sigma_{j,k}=\min\big\{1,\sigma^x_{j,k},\sigma^y_{j,k}\big\},
\label{corr2}
\end{equation}
and the scaling factors $\sigma^x_{j,k}$ and $\sigma^y_{j,k}$ are computed by
\begin{equation}
\sigma^x_{j,k}:=\left\{\begin{aligned}&\min\big\{1,\sigma^{x,1}_{j,k},\sigma^{x,2}_{j,k}\big\}&&\mbox{if }\sigma^{x,1}_{j,k}>0,\,
\sigma^{x,2}_{j,k}>0,\mbox{ and }\xbar A_{j,k}\ne0\\&0&&\mbox{otherwise}\end{aligned}\right.
\label{corr3}
\end{equation}
and
\begin{equation}
\sigma^y_{j,k}:=\left\{\begin{aligned}&\min\big\{1,\sigma^{y,1}_{j,k},\sigma^{y,2}_{j,k}\big\}&&\mbox{if }\sigma^{y,1}_{j,k}>0,\,
\sigma^{y,2}_{j,k}>0,\mbox{ and }\xbar B_{j,k}\ne0\\&0&&\mbox{otherwise},\end{aligned}\right.
\label{corr4}
\end{equation}
where 
\begin{equation}
\begin{aligned}
&\sigma^{x,1}_{j,k}=\frac{2\big((\widehat b_1)^{\rm E}_{j,k}-(\xbar{b_1})_{j,k}\big)}{\dx\,\xbar A_{j,k}},\quad
\sigma^{x,2}_{j,k}=\frac{2\big((\xbar{b_1})_{j,k}-(\widehat b_1)^{\rm W}_{j,k}\big)}{\dx\,\xbar A_{j,k}},\\
&\sigma^{y,1}_{j,k}=\frac{2\big((\widehat b_2)^{\rm N}_{j,k}-(\xbar{b_2})_{j,k}\big)}{\dy\,\xbar B_{j,k}},\quad
\sigma^{y,2}_{j,k}=\frac{2\big((\xbar{b_2})_{j,k}-(\widehat b_2)^{\rm S}_{j,k}\big)}{\dy\,\xbar B_{j,k}},
\end{aligned}
\label{corr5}
\end{equation}
and $(\widehat b_1)^{\rm E,W}_{j,k}$ and $(\widehat b_2)^{\rm N,S}_{j,k}$ denote the point values of $b_1$ and $b_2$, which have been
reconstructed as described above. We then correct the corresponding one-sided point values:
\begin{equation}
\begin{aligned}
&((b_1)_{j,k})^{\rm E}=(\xbar{b_1})_{j,k}+\frac{\dx}{2}((b_1)_x)_{j,k},\quad
((b_1)_{j+1,k})^{\rm W}=(\xbar{b_1})_{j+1,k}-\frac{\dx}{2}((b_1)_x)_{j+1,k},\\
&((b_2)_{j,k})^{\rm N}=(\xbar{b_2})_{j,k}+\frac{\dy}{2}((b_2)_y)_{j,k},\quad
((b_2)_{j,k+1})^{\rm S}=(\xbar{b_2})_{j,k+1}-\frac{\dy}{2}((b_2)_y)_{j,k+1}.
\end{aligned}
\label{corr}
\end{equation}

\smallskip
\noindent
$\bullet$ The one-sided point values $\widetilde{\mU}^{\rm E,W,N,S}_{j,k}$ are obtained by applying the generalized minmod reconstruction
directly to the $A$ and $B$ fields.

\smallskip
\noindent
$\bullet$ The point values of the global variables $\bm I^x$ and $\bm I^y$ in \eref{3.5a} are computed recursively. We first set
$\widehat x=x_\hf$ and $\widehat y=y_\hf$ so that $(\bm I^x)^-_{\hf,k}:=\bm0$ and $(\bm I^y)^-_{j,\hf}:=\bm0$, and then evaluate
$(\bm I^x)^+_{\hf,k}=\bm Q^x_{\bm\Psi,\hf,k}$, $(\bm I^y)^+_{j,\hf}=\bm Q^y_{\bm\Psi,j,\hf}$, and
\begin{equation*}
\begin{aligned}
&(\bm I^x)^-_{\jph,k}=(\bm I^x)^+_{\jmh,k}+\bm Q^x_{j,k},\quad(\bm I^x)^+_{\jph,k}=(\bm I^x)^-_{\jph,k}+\bm Q^x_{\bm\Psi,\jph,k},\\
&(\bm I^y)^-_{j,\kph}=(\bm I^y)^+_{j,\kmh}+\bm Q^y_{j,k},\quad(\bm I^y)^+_{j,\kph}=(\bm I^y)^-_{j,\kph}+\bm Q^y_{\bm\Psi,j,\kph},
\end{aligned}
\end{equation*}
for $j=1,\dots,N_x$, $k=1,\dots,N_y$. Here, $\bm Q^x_{j,k}$, $\bm Q^x_{\bm\Psi,\jph,k}$, $\bm Q^y_{j,k}$, and $\bm Q^y_{\bm\Psi,j,\kph}$ are
the terms reflecting the contribution of the nonconservative terms $Q^x(\mU)\mU_x$ and $Q^y(\mU)\mU_y$. For the details on evaluating these
terms, we refer the readers to \cite[\S2.2.3]{Chertock24}.

\smallskip
\noindent
$\bullet$ $P^{\rm LCD}_{\jph,k}$, $M^{\rm LCD}_{\jph,k}$, $Q^{\rm LCD}_{\jph,k}$, $P^{\rm LCD}_{j,\kph}$, $M^{\rm LCD}_{j,\kph}$, and
$Q^{\rm LCD}_{j,\kph}$ in \eref{3.4} are diagonal matrices
$$
\begin{aligned}
&P_{\jph,k}^{\rm LCD}={\rm diag}\left(\left(P_1^{\rm LCD}\right)_{\jph,k},\ldots,\left(P_8^{\rm LCD}\right)_{\jph,k}\right)\!,\!&&
P_{j,\kph}^{\rm LCD}={\rm diag}\left(\left(P_1^{\rm LCD}\right)_{j,\kph},\ldots,\left(P_8^{\rm LCD}\right)_{j,\kph}\right)\!,\\
&M_{\jph,k}^{\rm LCD}={\rm diag}\left(\left(M_1^{\rm LCD}\right)_{\jph,k},\ldots,\left(M_8^{\rm LCD}\right)_{\jph,k}\right)\!,\!&&
M_{j,\kph}^{\rm LCD}={\rm diag}\left(\left(M_1^{\rm LCD}\right)_{j,\kph},\ldots,\left(M_8^{\rm LCD}\right)_{j,\kph}\right)\!,\\
&Q_{\jph,k}^{\rm LCD}={\rm diag}\left(\left(Q_1^{\rm LCD}\right)_{\jph,k},\ldots,\left(Q_8^{\rm LCD}\right)_{\jph,k}\right)\!,\!&&
Q_{j,\kph}^{\rm LCD}={\rm diag}\left(\left(Q_1^{\rm LCD}\right)_{j,\kph},\ldots,\left(Q_8^{\rm LCD}\right)_{j,\kph}\right)\!,
\end{aligned}
$$
where
$$
\begin{aligned}
&\big((P^{\rm LCD}_i)_{\jph,k},(M^{\rm LCD}_i)_{\jph,k},(Q^{\rm LCD}_i)_{\jph,k}\big)=
\frac{\Big((\lambda^+_i)_{\jph,k},-(\lambda^-_i)_{\jph,k},(\lambda^+_i)_{\jph,k}(\lambda^-_i)_{\jph,k}\Big)}
{(\lambda^+_i)_{\jph,k}-(\lambda^-_i)_{\jph,k}},\\
&\big((P^{\rm LCD}_i)_{j,\kph},(M^{\rm LCD}_i)_{j,\kph},(Q^{\rm LCD}_i)_{j,\kph}\big)=
\frac{\Big((\lambda^+_i)_{j,\kph},-(\lambda^-_i)_{j,\kph},(\lambda^+_i)_{j,\kph}(\lambda^-_i)_{j,\kph}\Big)}
{(\lambda^+_i)_{j,\kph}-(\lambda^-_i)_{j,\kph}},
\end{aligned}
$$
with
\begin{equation}
\begin{aligned}
&(\lambda^+_i)_{\jph,k}=\max\left\{\lambda_i\big(C^x(\mU^{\rm E}_{j,k})\big),\,\lambda_i\big(C^x(\mU^{\rm W}_{j+1,k})\big),\,
\varepsilon\right\},\\
&(\lambda^-_i)_{\jph,k}=\min\left\{\lambda_i\big(C^x(\mU^{\rm E}_{j,k})\big),\,\lambda_i\big(C^x(\mU^{\rm W}_{j+1,k})\big),\,
-\varepsilon\right\},\\
&(\lambda^+_i)_{j,\kph}=\max\left\{\lambda_i\big(C^y(\mU^{\rm N}_{j,k})\big),\,\lambda_i\big(C^y(\mU^{\rm S}_{j,k+1})\big),\,
\varepsilon\right\},\\
&(\lambda^-_i)_{j,\kph}=\min\left\{\lambda_i\big(C^y(\mU^{\rm N}_{j,k})\big),\,\lambda_i\big(C^y(\mU^{\rm S}_{j,k+1})\big),\,
-\varepsilon\right\},
\end{aligned}
\label{3.7}
\end{equation}
and $\lambda_i\big(C^x(\mU))$ and $\lambda_i\big(C^y(\mU))$ are eigenvalues of $C^x$ and $C^y$, $i=1,\dots,8$, respectively; see Appendix
\ref{appA1} for details. In \eref{3.7}, $\varepsilon$ is a small desingularization constant, which is taken to be $10^{-8}$ in all of the
numerical examples reported in \S\ref{sec4}. 

\smallskip
\noindent
$\bullet$ $R^x_{\jph,k}$ and $R^y_{j,\kph}$ are the matrices of right eigenvectors of $\widehat C^x_{\jph,k}=C^x(\widehat\mU_{\jph,k})$ and
$\widehat C^y_{j,\kph}=C^y(\widehat\mU_{j,\kph})$, respectively. Here, we take
$\widehat\mU_{\jph,k}=(\,\xbar\mU_{j,k}+\,\xbar\mU_{j+1,k})/2$ and $\widehat\mU_{j,\kph}=(\,\xbar\mU_{j,k}+\,\xbar\mU_{j,k+1})/2$. 

\smallskip
\noindent
$\bullet$ $s^\pm_{\jph,k}$ and $s^\pm_{j,\kph}$ are one-sided local speeds of propagation in the $x$- and $y$-direction, respectively. They
are estimated as in \cite[\S2.2.2]{Chertock24}.
\begin{remark}
We stress that the correction of the point values \eref{corr1}--\eref{corr} together with the result proven in
\cite[Theorem 2.2]{Chertock24} enforce the local divergence-free condition \eref{4.1}.
\end{remark}

\section{Numerical Examples}\label{sec4}
In this section, we test the developed LCD-PCCU scheme on a number of numerical examples and compare the obtained results with those
computed by the PCCU scheme from \cite{Chertock24}. We numerically integrate the ODE systems \eref{3.3} by the three-stage third-order
strong stability preserving (SSP) Runge-Kutta method (see, e.g., \cite{Gottlieb11,Gottlieb12}), use the CFL number $0.25$, and set the
minmod parameter $\theta=1.3$ (except for Example 5, where we take $\theta=1$ to reduce oscillations). The specific heat ratio $\gamma$ is
either $2$ (Example 1), $5/3$ (Examples 2--4), or $1.4$ (Example 5).

In Examples 2 and 3, we will demonstrate how the discrete divergence $(\bnabla\cdot\bm b)_{j,k}$, defined in \eref{4.1}, increases in time
if the correction \eref{corr1}--\eref{corr5} of the slopes $((b_1)_x)_{j,k}$ and $((b_2)_y)_{j,k}$ is not implemented (the corresponding
scheme will be referred to as Uncorrected LCD-PCCU scheme).

\subsubsection*{Example 1---Brio-Wu Shock-Tube Problem}
In the first example taken from \cite{Brio88}, we consider the one-dimensional (1-D) Riemann problem, which is a benchmark widely used to
test the ability of schemes to capture compound waves that emerge out of the initial data,
\begin{equation*}
(\rho,u,v,w,b1,b2,b3,p)(x,0)=\left\{\begin{aligned}&(1,0,0,0,0.75,1,0,1)&&\mbox{if}~x<0,\\&(0.125,0,0,0,0.75,-1,0,0.1)&&\mbox{otherwise},
\end{aligned}\right.
\end{equation*}
which depend on $x$ only. We conduct a 2-D computation on the domain $[-1,1]\times[-0.01,0.01]$ subject to the free boundary conditions. 

\begin{figure}[ht!]
\centerline{\includegraphics[trim=0.6cm 0.4cm 1.cm 0.3cm, clip, width=5.cm]{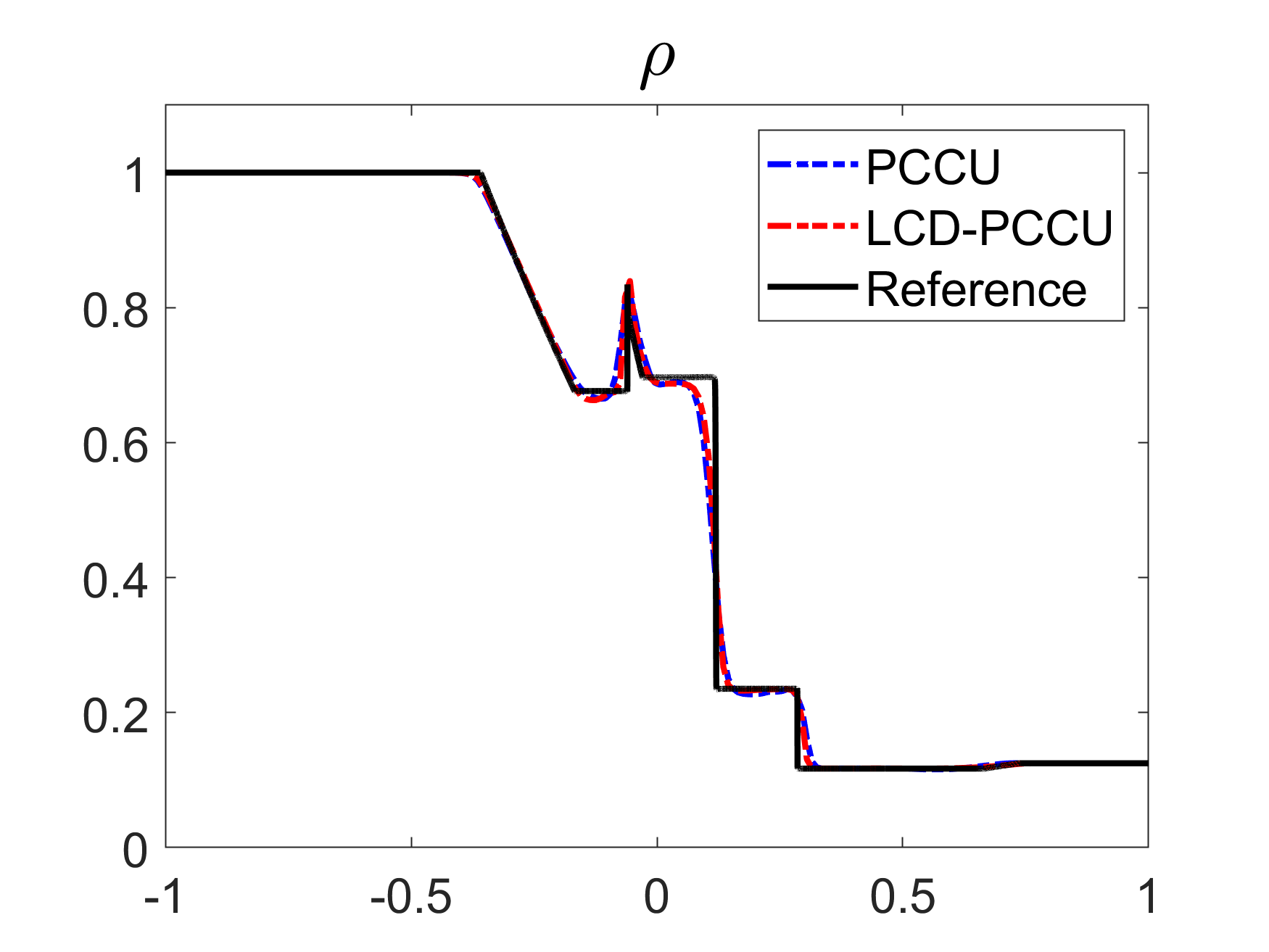}\hspace*{0.5cm}
            \includegraphics[trim=0.6cm 0.4cm 1.cm 0.3cm, clip, width=5.cm]{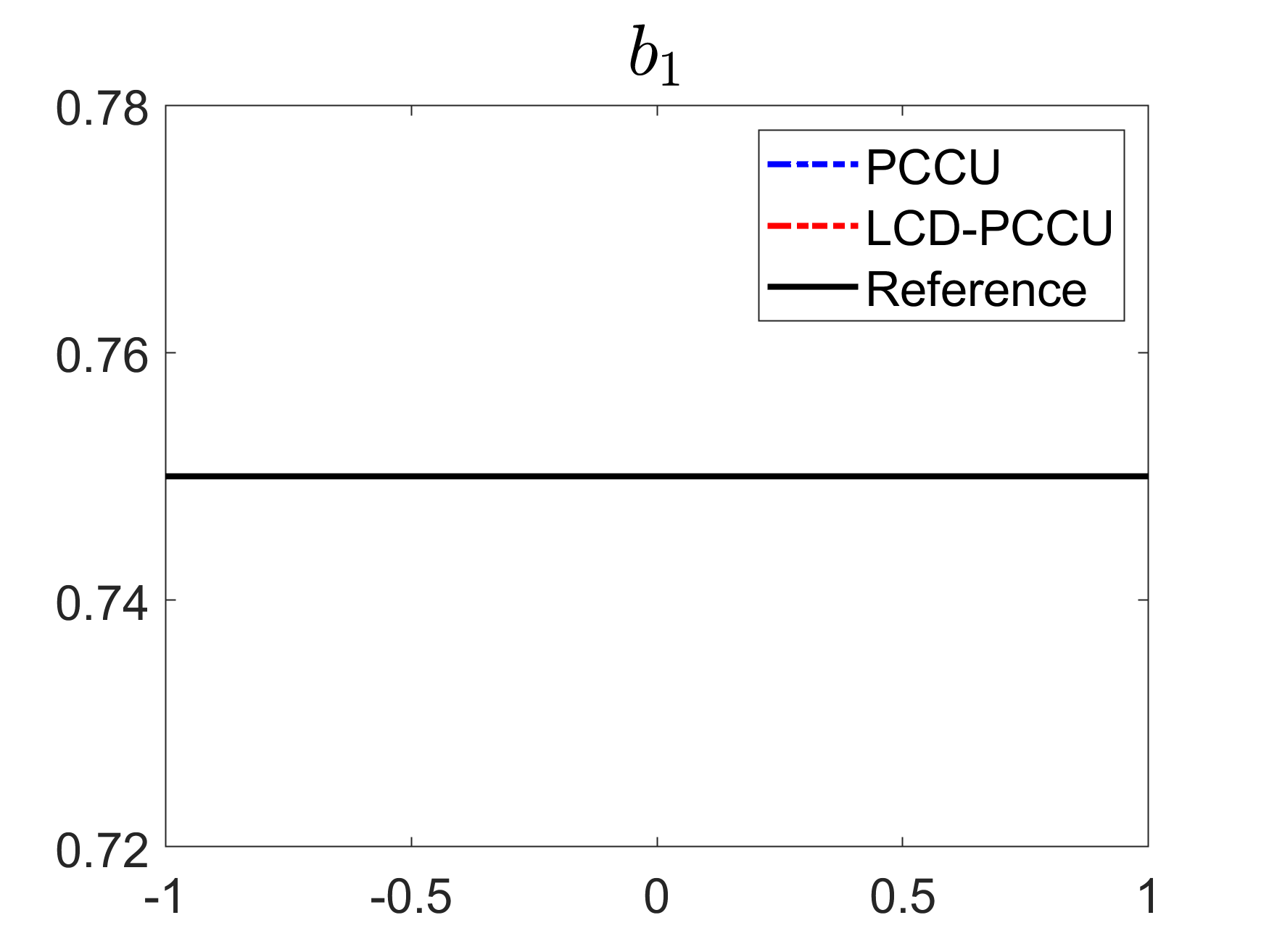}\hspace*{0.5cm}
            \includegraphics[trim=0.6cm 0.4cm 1.cm 0.3cm, clip, width=5.cm]{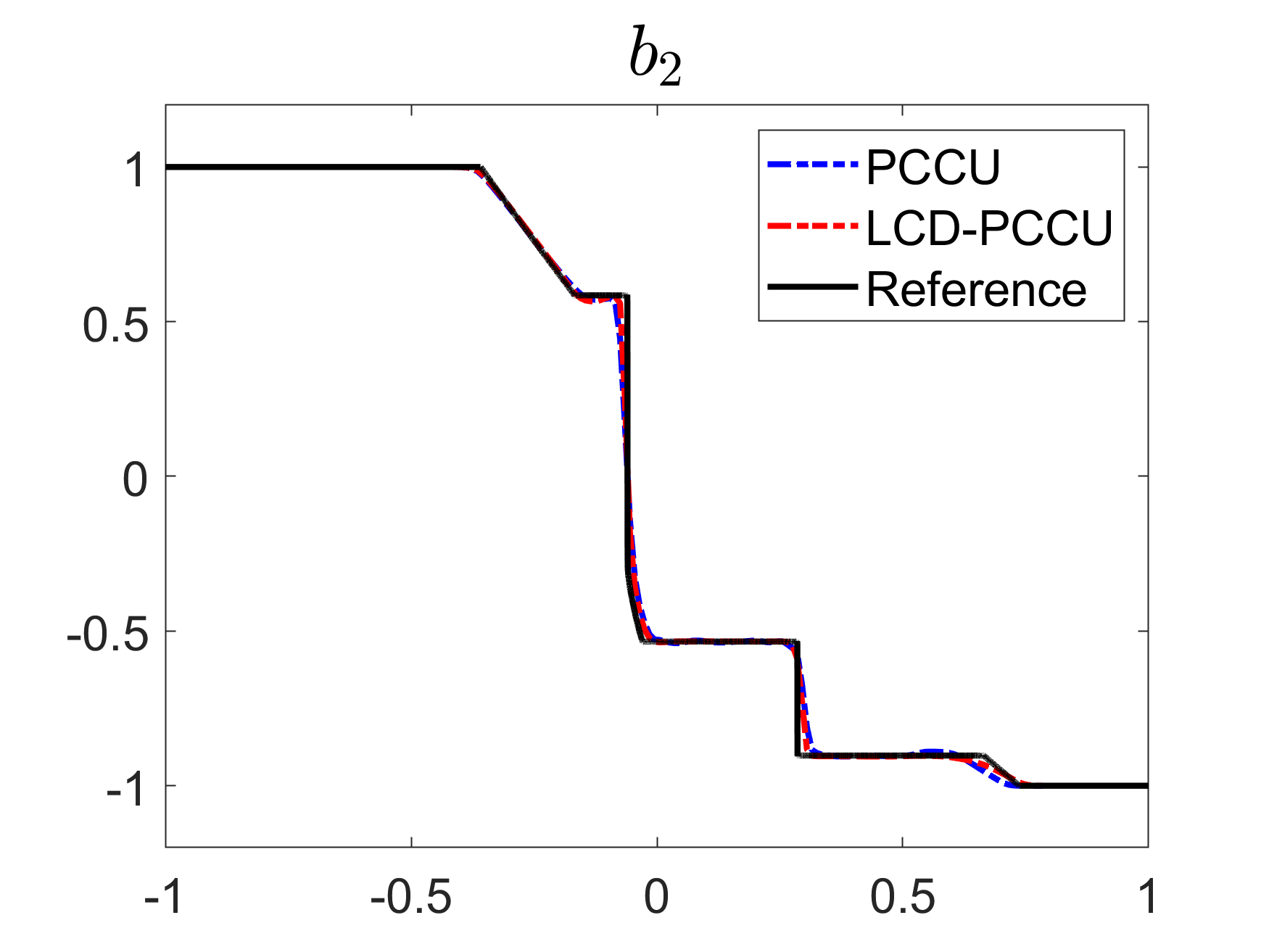}}
\vskip 12pt
\centerline{\includegraphics[trim=0.6cm 0.4cm 1.cm 0.3cm, clip, width=5.cm]{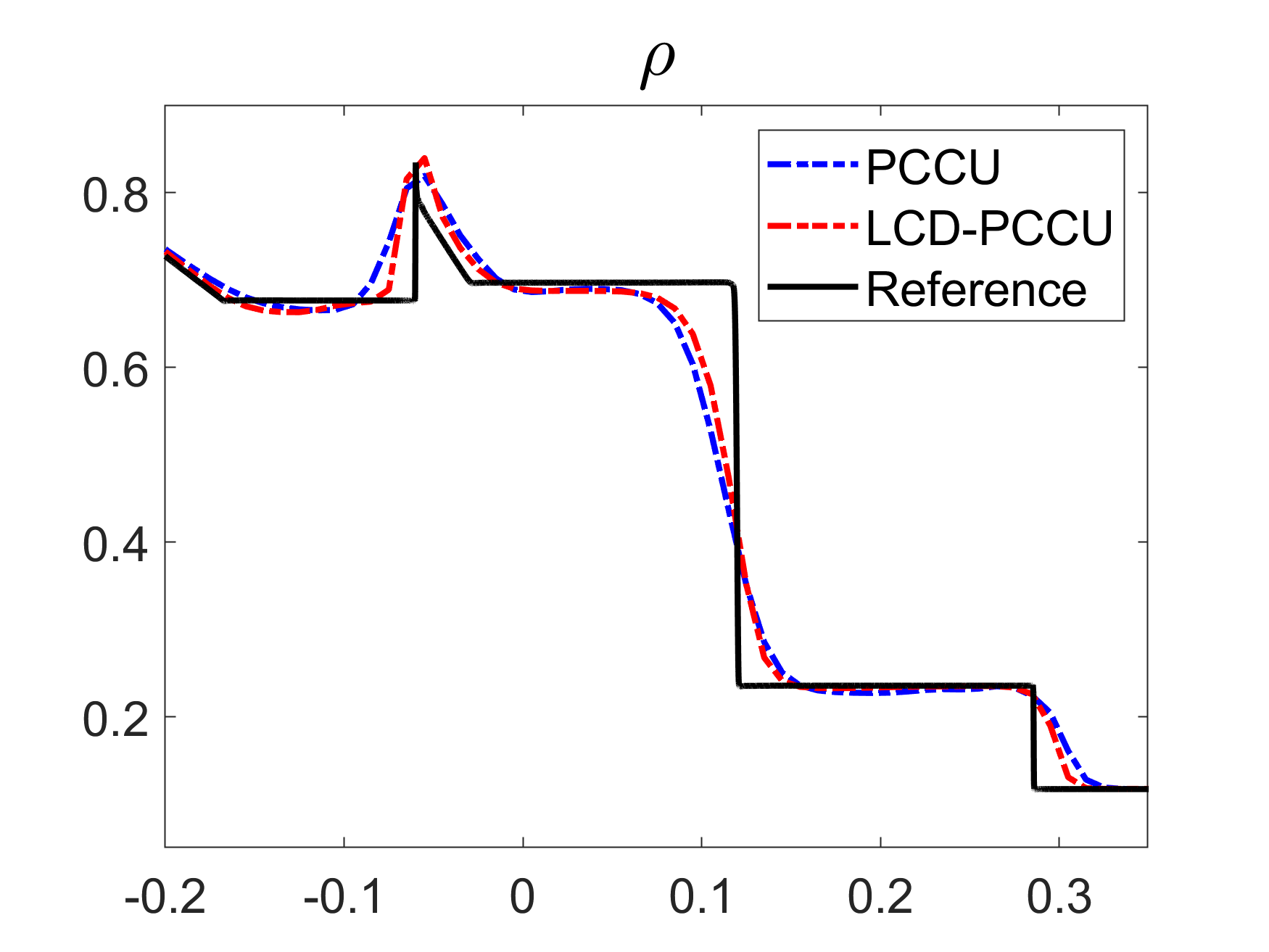}\hspace*{0.5cm}
            \includegraphics[trim=0.6cm 0.4cm 1.cm 0.3cm, clip, width=5.cm]{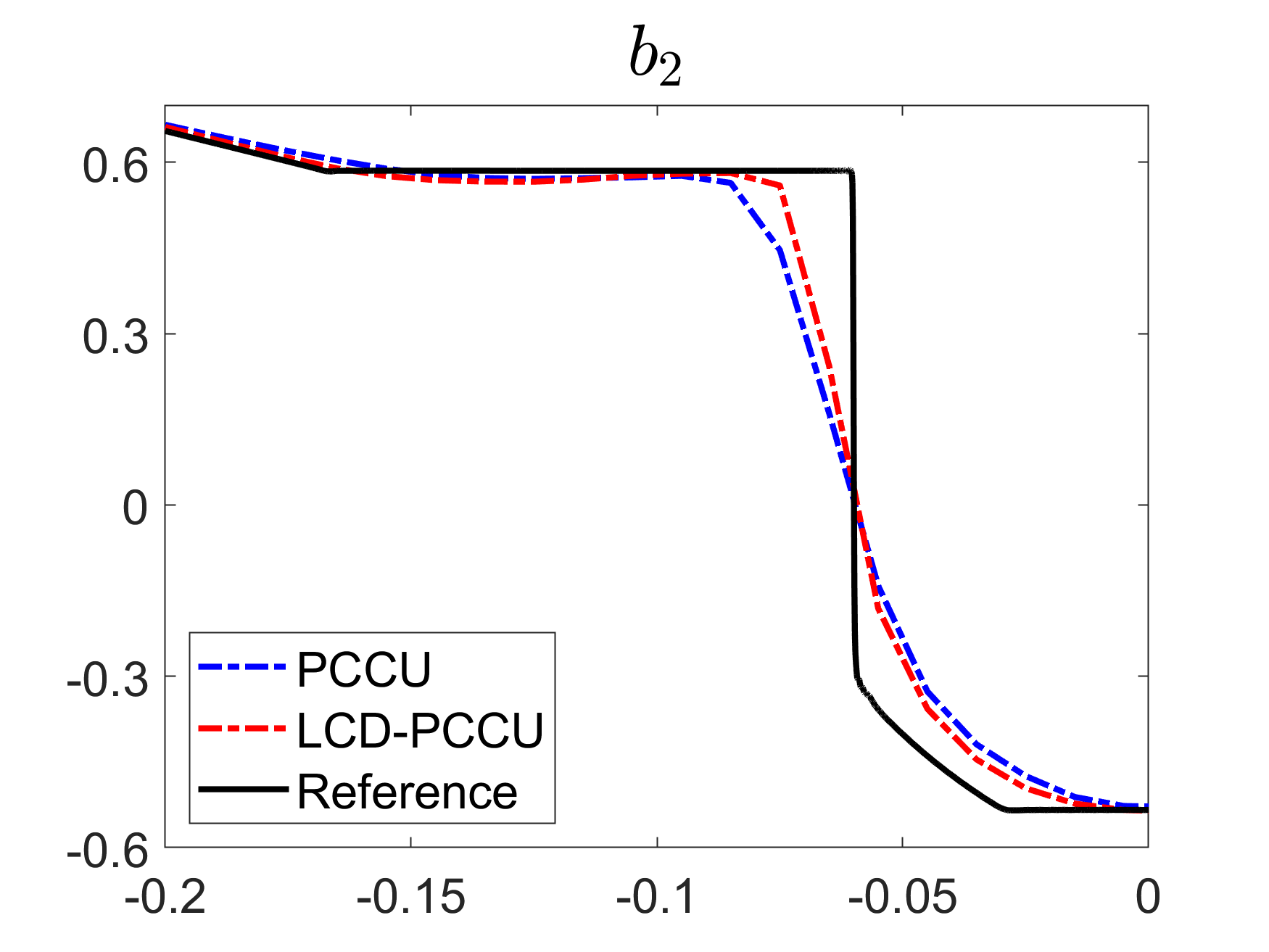}\hspace*{0.5cm}
            \includegraphics[trim=0.6cm 0.4cm 1.cm 0.3cm, clip, width=5.cm]{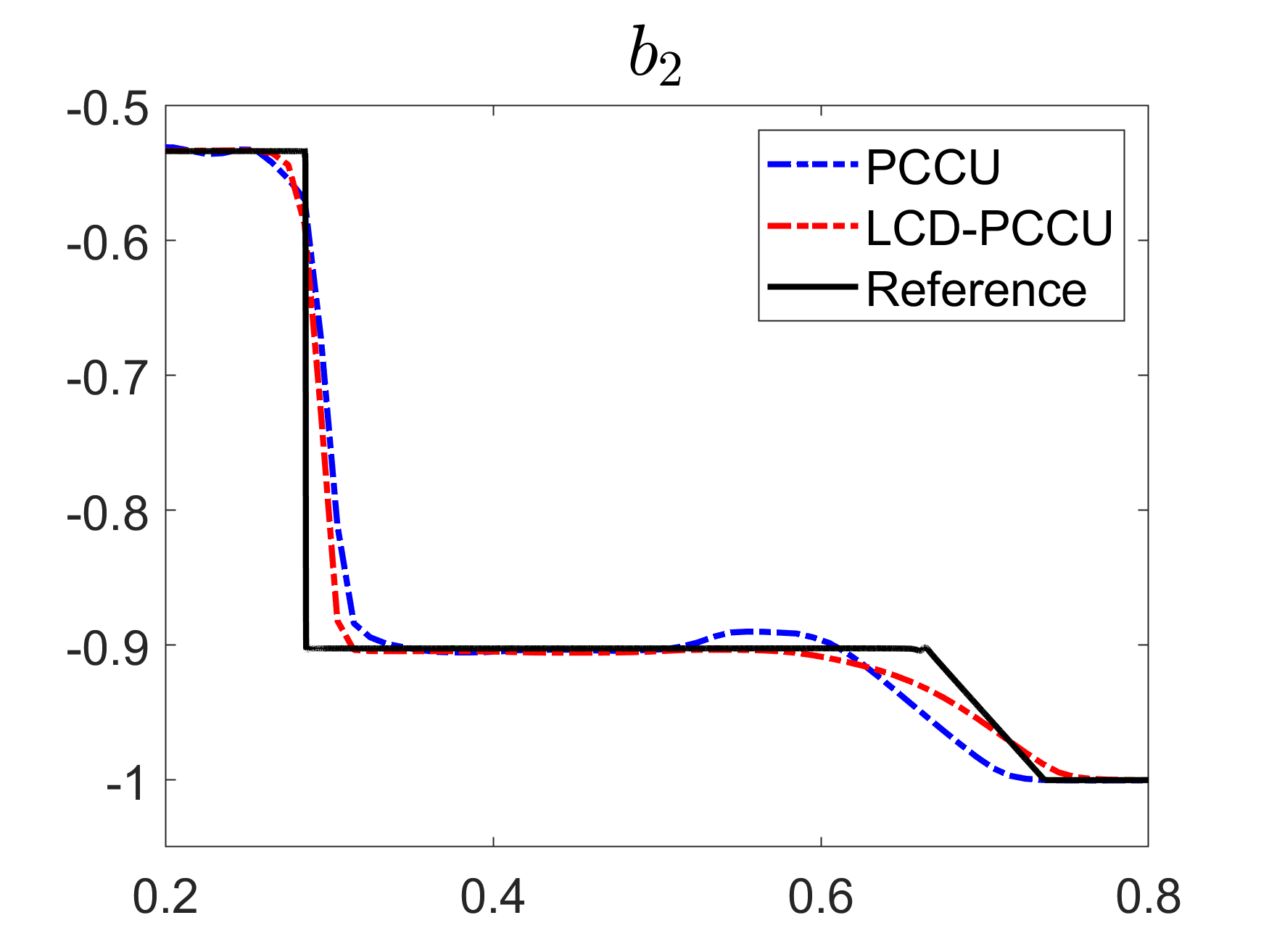}}            
\caption{\sf Example 1: $\rho$, $b_1$, and $b_2$ computed by the LCD-PCCU and PCCU schemes (top row) and zooms for $\rho$ and $b_2$ at
$x\in[-0.2,0.35]$, $[-0.2,0]$, and $[0.2,0.8]$ (bottom row).\label{fig11}}
\end{figure}

We compute the solutions by the LCD-PCCU and PCCU schemes until the final time $t=0.2$ on a uniform mesh consisting of $200\times2$ cells.
The cross-sectional profiles at $y=0$ of $\rho$, $b_1$, and $b_2$ are presented in Figure \ref{fig11} along with the reference solution
computed by the PCCU scheme on a significantly finer mesh of $10000\times2$ cells. One can observe that the solution consists of several
nonsmooth structures, including rarefaction waves, shock waves traveling at various speeds, a contact discontinuity, and a compound shock
wave. Both the LCD-PCCU and PCCU schemes successfully capture all of these complex structures. However, the numerical results obtained by
the LCD-PCCU scheme exhibit somewhat higher resolution compared to those produced by the PCCU scheme. 

\subsubsection*{Example 2---Circularly Polarized Alfv\'en Wave} 
In the second example taken from \cite{Toth2000}, we consider the time evolution of a circularly polarized Alfv\' en wave that travels at a
constant speed at an angle of $\alpha=\pi/6$ with respect to the $x$-axis. In this example, designed to check the experimental order of
accuracy of the studied schemes, the initial conditions are
$$
\begin{aligned} 
&\rho(x,y,0)\equiv1,\quad~~u(x,y,0)=v_\|\cos\alpha+v_\perp\sin\alpha,\quad v(x,y,0)=v_\|\sin\alpha-v_\perp\cos\alpha,\\
&p(x,y,0)\equiv0.1,\quad b_1(x,y,0)=b_\|\cos\alpha+b_\perp\sin\alpha,\quad b_2(x,y,0)=b_\|\sin\alpha-b_\perp\cos\alpha,\\
&\qquad \qquad\qquad w(x,y,0)=b_3(x,y,0)=0.1\cos\big[2\pi(x\cos\alpha+y\sin\alpha)\big], 
\end{aligned}
$$
where
$$
v_{\|}=0,\quad b_{\|}=1,\quad v_\perp=b_\perp=0.1\sin\big[2\pi(x\cos\alpha+y\sin\alpha)\big],
$$
and the periodic boundary conditions are imposed in the computational domain
$\big[0,\frac{1}{\cos\alpha}\big]\times\big[0,\frac{1}{\sin\alpha}\big]$. It is easy to show that the solution of the resulting
initial-boundary value problem is a traveling wave, which returns to its initial position at any integer time $t$.

We compute the solutions by the LCD-PCCU and PCCU schemes until the final time $t=5$ on a sequence of uniform meshes with $20\times20$,
$40\times40$, $80\times80$, $160\times160$, and $320\times320$ cells, and compute the $L^1$-norm of the differences between the numerical
and exact solutions. We report the $L^1$-errors and corresponding experimental rates of convergence for both $u$ and $b_3$ in Table
\ref{tab21}, where one can see that while both the LCD-PCCU and PCCU schemes achieve the expected second order of accuracy, the magnitudes
of the errors are slightly smaller for the proposed LCD-PCCU scheme.
\begin{table}[ht!]
\caption{\sf Example 2: $L^1$-errors and experimental convergence rates for $u$ and $b_3$ computed by the LCD-PCCU and PCCU schemes.\label{tab21}}
\centering
\begin{tabular}{|c|cc|cc|cc|cc|}
\hline
\multirow{3}{3em}{Mesh}&\multicolumn{4}{c|}{LCD-PCCU Scheme}&\multicolumn{4}{c|}{PCCU Scheme}\\
\cline{2-9}
&\multicolumn{2}{c|}{$u$}&\multicolumn{2}{c|}{$b_3$}&\multicolumn{2}{c|}{$u$}&\multicolumn{2}{c|}{$b_3$}\\
\cline{2-9}
&Error&Rate&Error&Rate&Error&Rate&Error&Rate\\
\hline
$20\times20$&2.69e-2&--&4.50e-2&--&3.38e-2&--&6.45e-2&--\\
$40\times40$&7.83e-3&1.78&1.18e-2&1.93&7.98e-3&2.08&1.67e-2&1.95\\
$80\times80$&2.29e-3&1.78&3.29e-3&1.84&2.46e-3&1.70&5.32e-3&1.65\\
$160\times160$&5.75e-4&1.99&8.31e-4&1.98&6.48e-4&1.93&1.41e-3&1.92\\
$320\times320$&1.34e-4&2.10&2.19e-4&1.93&1.54e-4&2.08&3.34e-4&2.07\\
\hline
\end{tabular}
\end{table}

\begin{figure}[ht!]
\centerline{\includegraphics[trim=0.6cm 0.4cm 1.cm 0.3cm, clip, width=6.0cm]{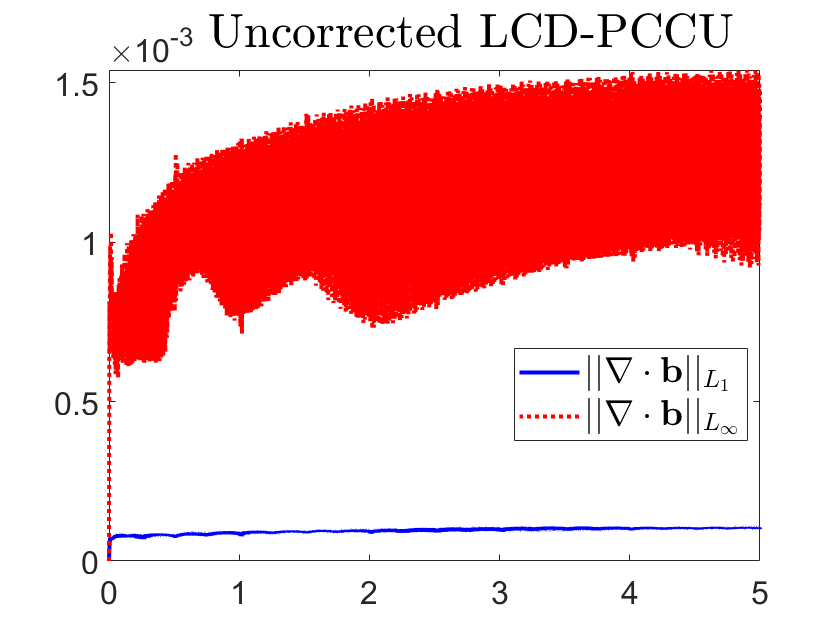}}
\caption{\sf Example 2: Time evolution of the $L^1$- and $L^\infty$-norms of $(\bnabla\cdot\bm b)_{j,k}$ computed by the Uncorrected
LCD-PCCU scheme on a uniform $320\times320$ mesh.\label{fig22}}
\end{figure}

Figure \ref{fig22} presents the time evolution of the $L^1$- and $L^\infty$-norms of $(\bnabla\cdot\bm b)_{j,k}$ computed by the Uncorrected
LCD-PCCU scheme on a uniform $320\times320$ mesh. As one can see, the magnitudes of both norms are comparable or even exceed the size of the
formal truncation error, which is about $10^{-5}$ on this grid. This suggests that the use of the Uncorrected LCD-PCCU scheme may lead to
a substantial numerical inaccuracy, and thus applying the divergence-free correction might be essential for ensuring the physical
consistency and long-term stability of the simulation.

\subsubsection*{Example 3---Orszag-Tang Vortex Problem}
In this example taken from \cite{Orszag_Tang_1979}, we consider the Orszag-Tang vortex problem, which has been widely used as a benchmark
due to the formation and interaction of multiple shocks as the system evolves in time and to the presence of many important features of MHD
turbulence. The initial conditions,
$$
\begin{aligned}
&\rho(x,y,0)\equiv\gamma^2,\quad u(x,y,0)=-\sin y,\quad v(x,y,0)=\sin x,\quad w(x,y,0)\equiv0,\\
&b_1(x,y,0)=-\sin y,\quad b_2(x,y,0)=\sin(2x),\quad b_3(x,y,0)\equiv0,\quad p(x,y,0)\equiv\gamma,
\end{aligned} 
$$
are prescribed in the computational domain $[0,2\pi]\times[0,2\pi]$ subject to the periodic boundary conditions.
 
\begin{figure}[ht!]
\centerline{\includegraphics[trim=1.cm 1.5cm 0.9cm 0.9cm, clip, width=6.0cm]{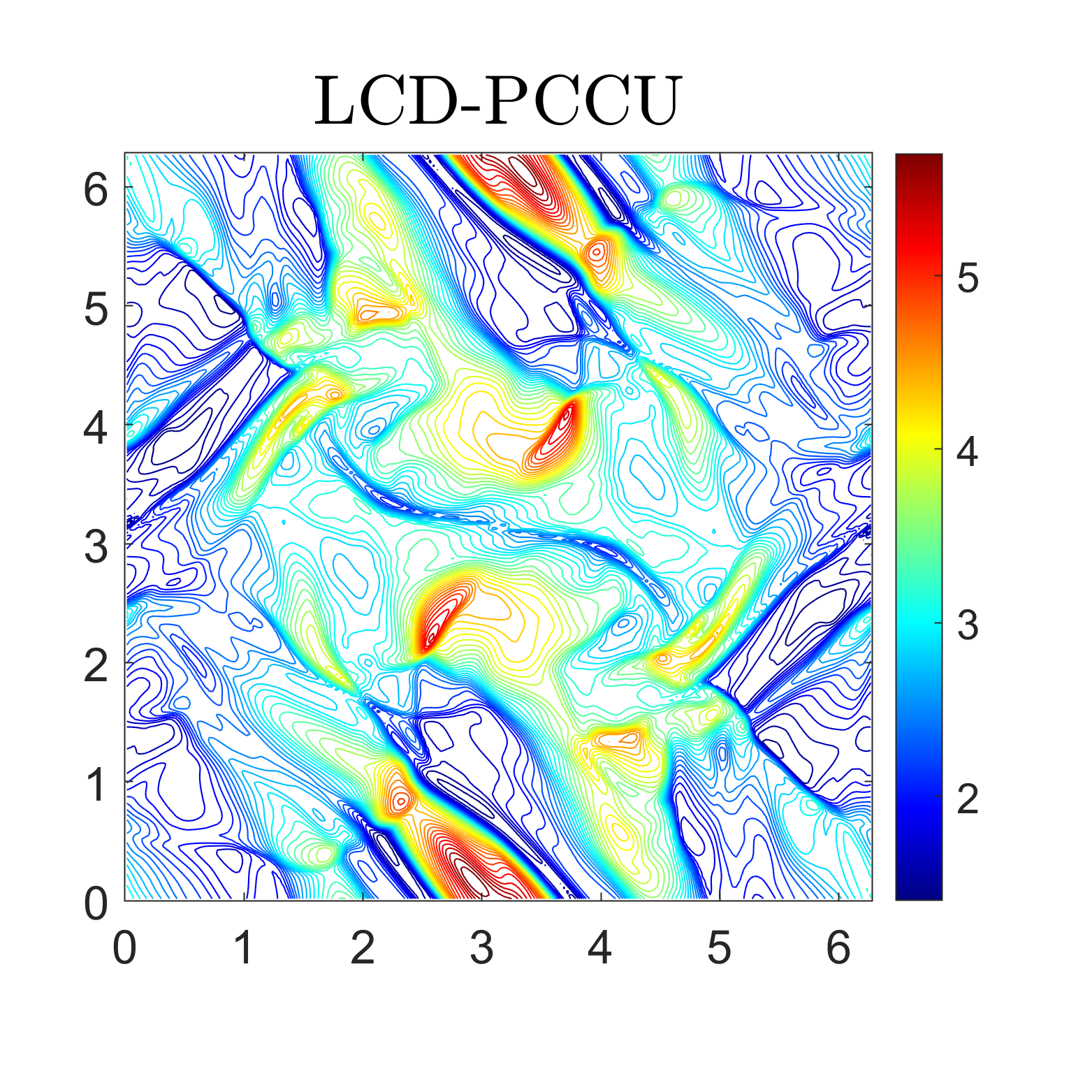}\hspace*{0.5cm}
            \includegraphics[trim=1.cm 1.5cm 0.9cm 0.9cm, clip, width=6.0cm]{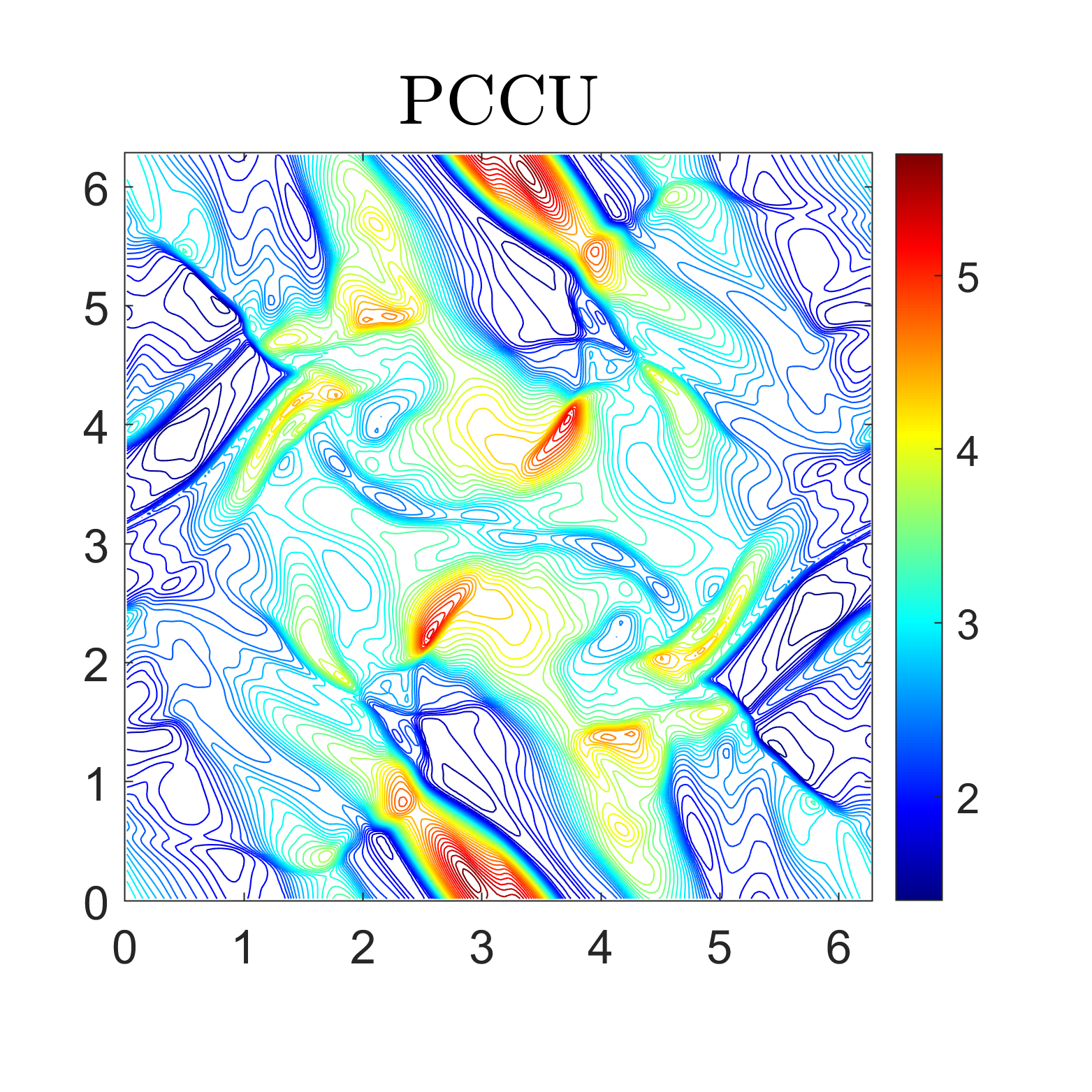}}
\caption{\sf Example 3: Density $\rho$ computed by the LCD-PCCU (left) and PCCU (right) schemes.\label{fig31}}
\end{figure}

\begin{figure}[ht!]
\centerline{\includegraphics[trim=0.9cm 0.4cm 1.4cm 0.3cm, clip, width=5.5cm]{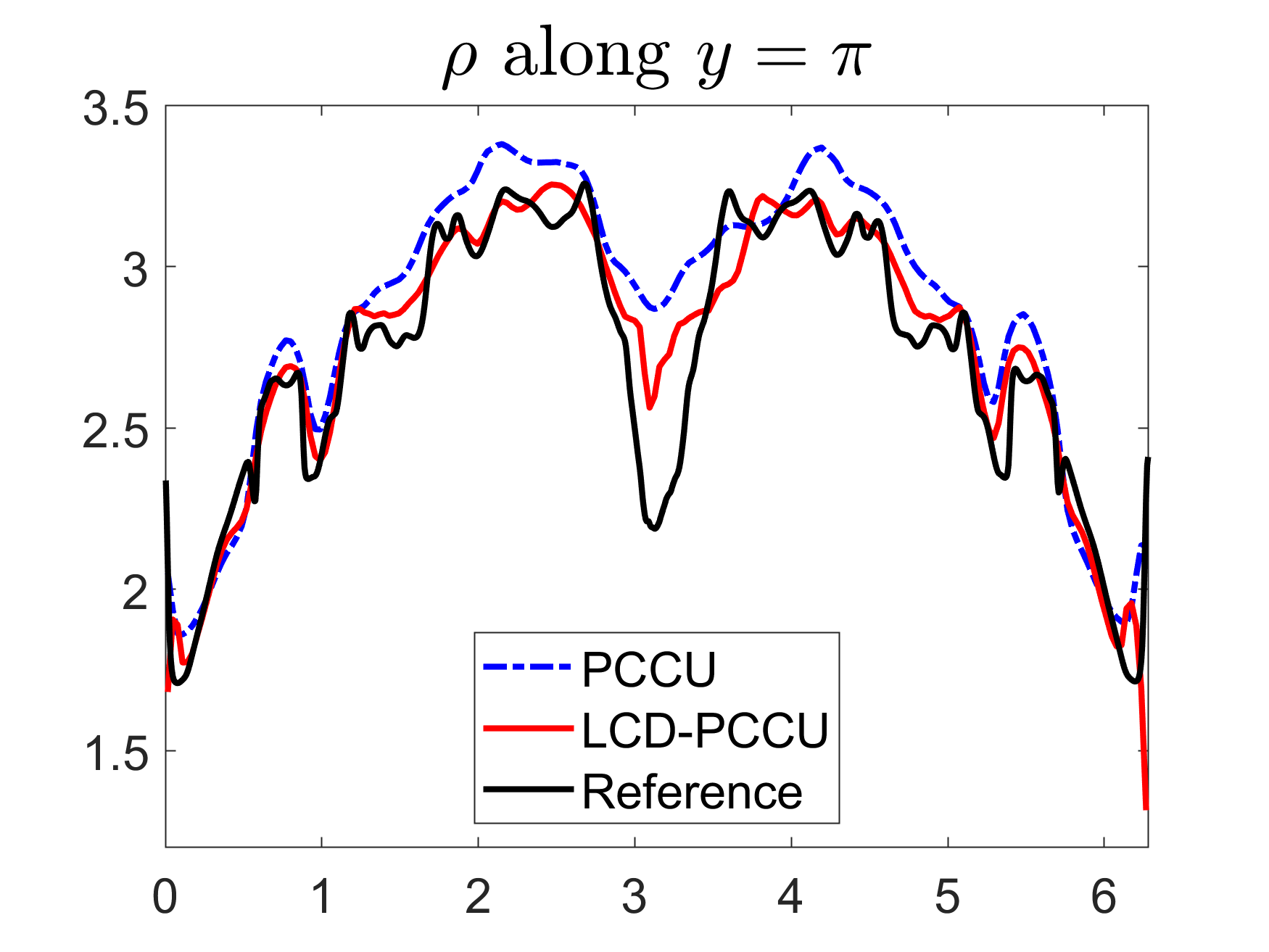}}
\caption{\sf Example 3: 1-D slices along the line $y=\pi$ of the solutions from Figure \ref{fig31} together with the reference solution.
\label{fig31a}}
\end{figure}

We compute the numerical solutions by both the LCD-PCCU and PCCU schemes until the final time $t=4$ using a uniform $200\times200$ mesh and
plot the obtained densities in Figure \ref{fig31}. As one can see, the LCD-PCCU solution is sharper, and this can be further seen in Figure
\ref{fig31a}, where we plot the 1-D slices of both densities along $y=\pi$ together with the reference solution computed by the PCCU scheme
on $1000\times1000$ uniform mesh.

\begin{figure}[ht!]
\centerline{\includegraphics[trim=0.6cm 0.4cm 1.cm 0.3cm, clip, width=6.0cm]{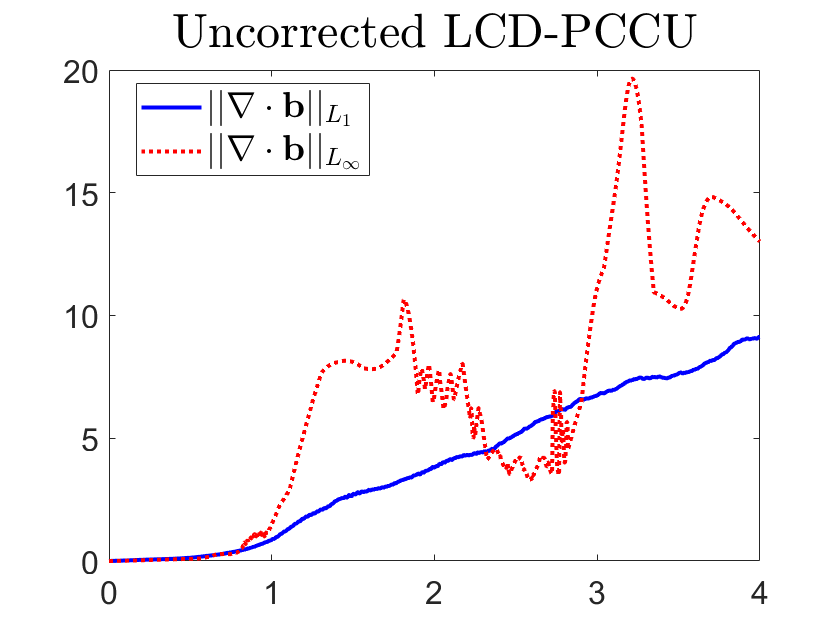}}
\caption{\sf Example 3: Time evolution of the $L^1$- and $L^\infty$-norms of $(\bnabla\cdot\bm b)_{j,k}$ computed by the Uncorrected
LCD-PCCU scheme.\label{fig32}}
\end{figure}

The time evolution of the $L^1$- and $L^\infty$-norms of $(\bnabla\cdot\bm b)_{j,k}$ computed by the Uncorrected LCD-PCCU scheme is
presented in Figure \ref{fig32}. As one can see, the magnitudes of both norms increase in time and in this example, they are several orders
of magnitude larger than the formal truncation error, which is about $10^{-3}$ on this grid. Consequently, applying the divergence-free
correction is essential for ensuring the physical consistency and long-term stability of the simulation.

\subsubsection*{Example 4---Rotor Problem} 
In this example, we study the ``second rotor problem'', originally introduced in \cite{BS1999,Toth2000} as a benchmark featuring a rapidly
rotating, dense fluid disk embedded in a stationary background. As time progresses, the disk undergoes both expansion and rotation. 

The initial conditions are
\begin{equation*} 
\begin{aligned} 
&(\rho,u,v)\Big|_{(x,y,0)}=\left\{\begin{aligned}&\left(10,\frac{0.5-y}{r_0},\frac{x-0.5}{r_0}\right),&&r<0.1,\\
&\left(1+9\mu,\frac{\mu(0.5-y)}{r},\frac{\mu(x-0.5)}{r}\right),&&0.1\le r\le0.115,\\&(1,0,0),&&r>0.115,\end{aligned}\right.\\
&w(x,y,0)=b_2(x,y,0)=b_3(x,y,0)\equiv0,\quad b_1(x,y,0)\equiv\frac{2.5}{\sqrt{4\pi}},\quad p(x,y,0)\equiv0.5, 
\end{aligned} 
\end{equation*}
where $r=\sqrt{(x-0.5)^2+(y-0.5)^2}$, $r_0=0.1$, and $\mu=(0.115-r)/0.015$. We use the periodic boundary conditions in the computational
domain $[0,1]\times[0,1]$.

\begin{figure}[ht!]
\centerline{\includegraphics[trim=0.5cm 1.6cm 0.9cm 1.2cm, clip, width=5.5cm]{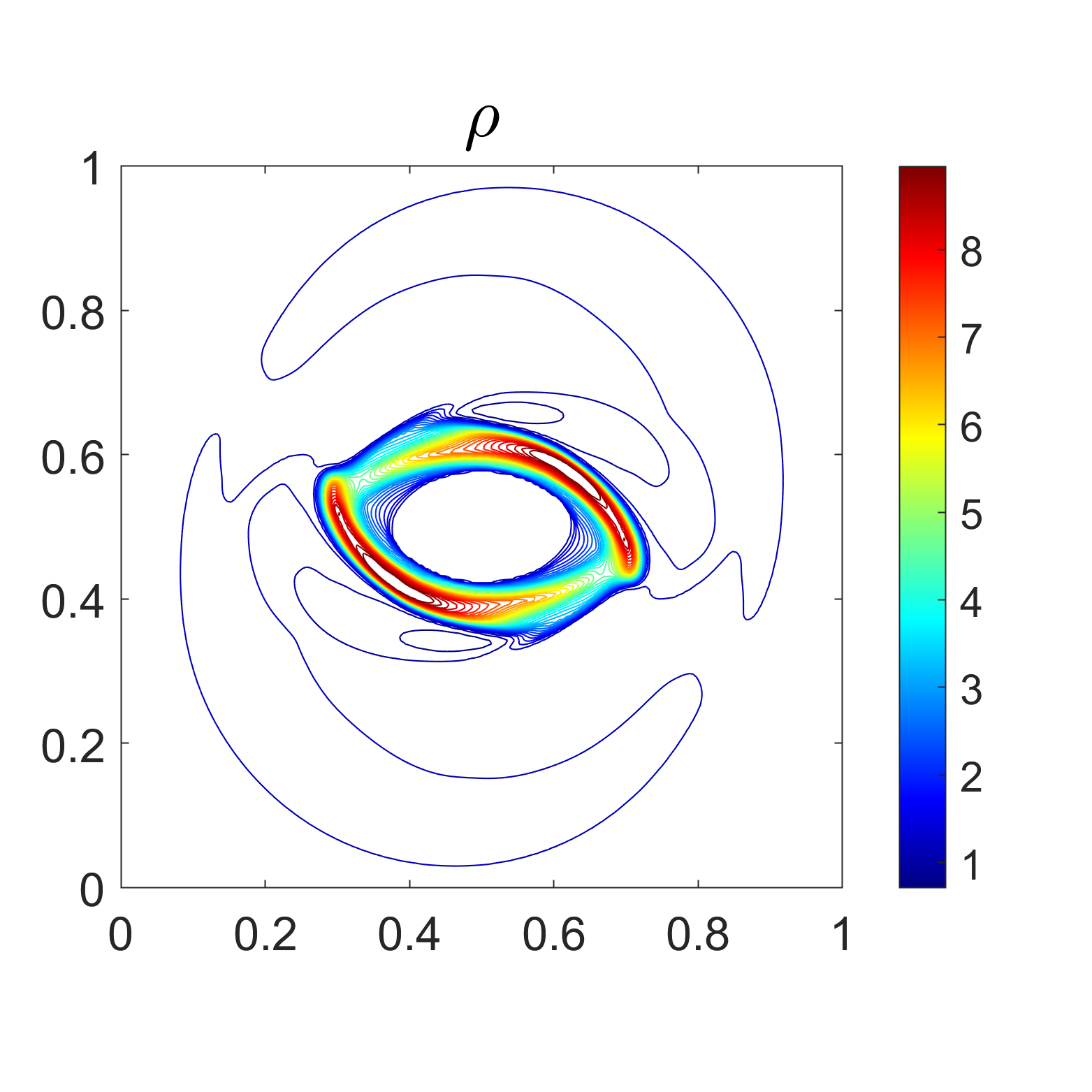}\hspace*{0.5cm}
            \includegraphics[trim=0.5cm 1.6cm 0.9cm 1.2cm, clip, width=5.5cm]{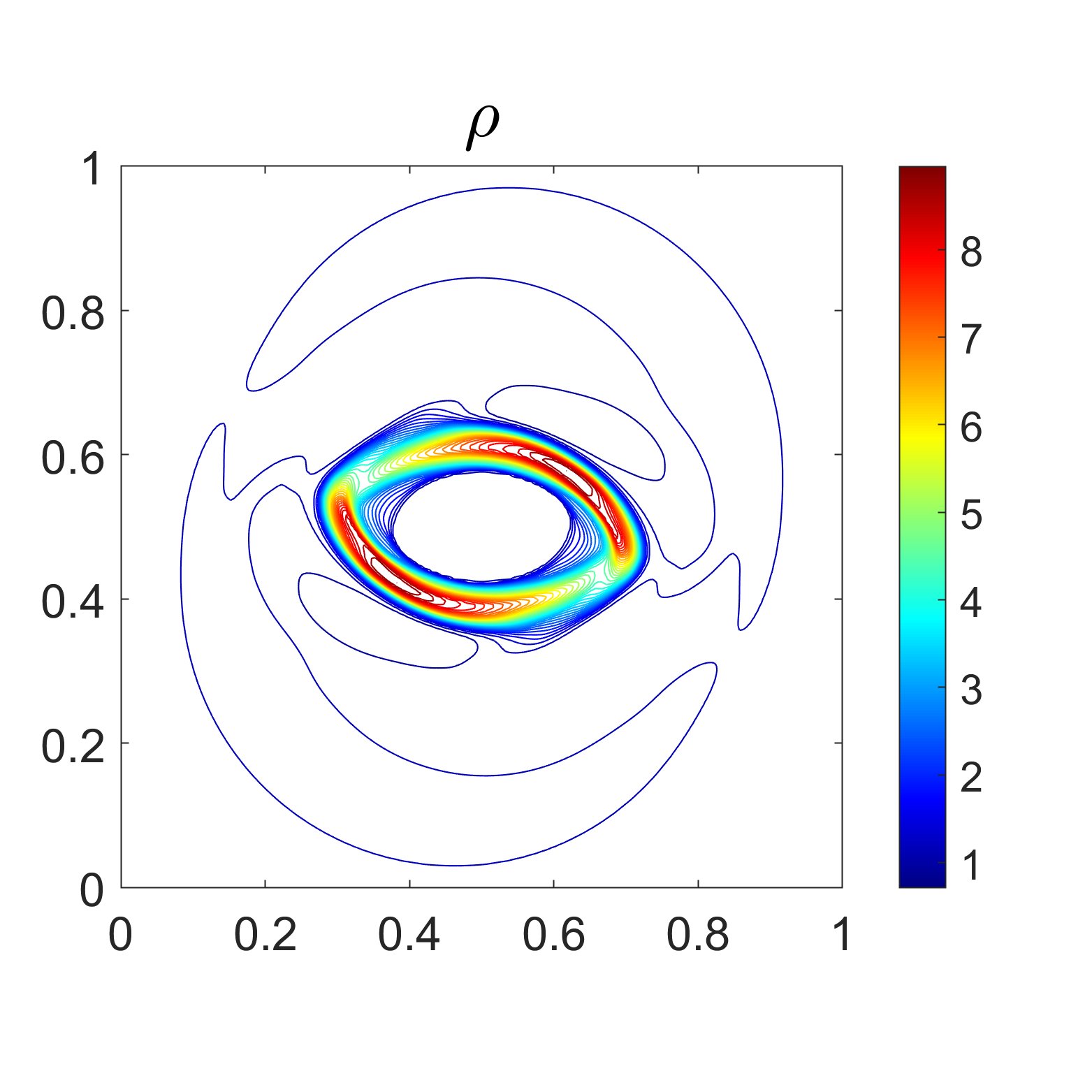}}
\vskip10pt
\centerline{\includegraphics[trim=0.5cm 1.6cm 0.9cm 1.2cm, clip, width=5.5cm]{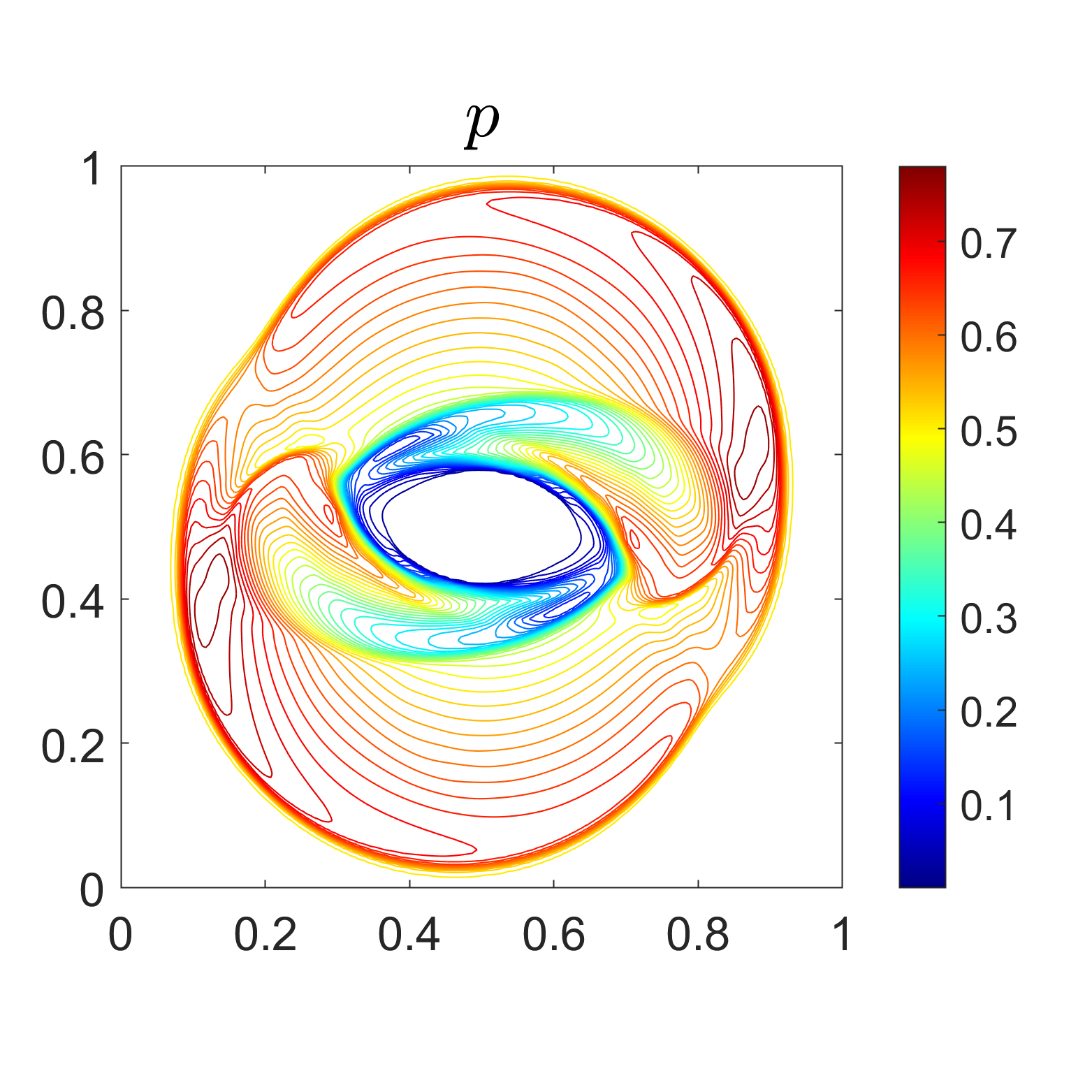}\hspace*{0.5cm}
            \includegraphics[trim=0.5cm 1.6cm 0.9cm 1.2cm, clip, width=5.5cm]{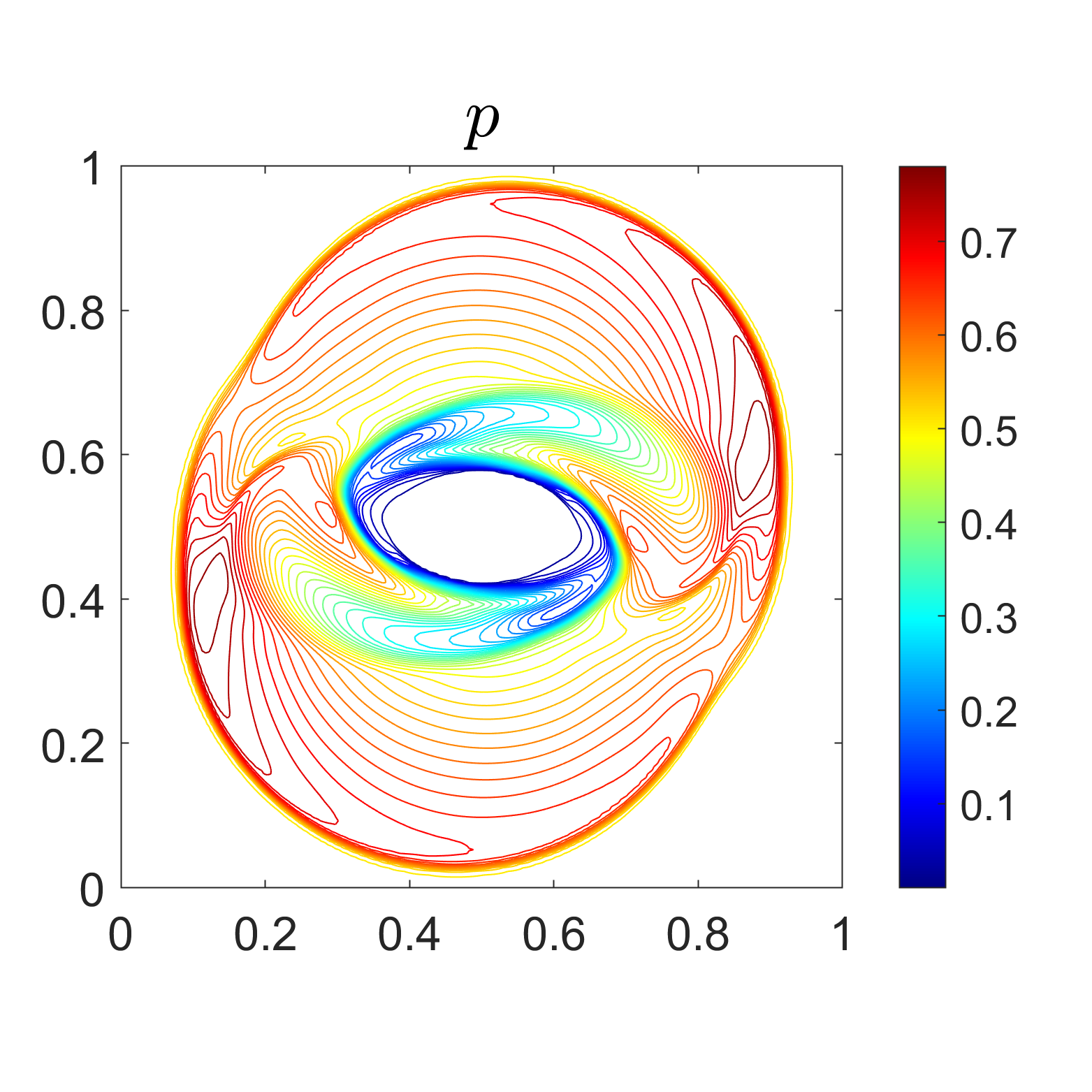}}
\caption{\sf Example 4: Density $\rho$ (top row) and pressure $p$ (bottom row) computed by the LCD-PCCU (left column) and PCCU (right
column) schemes.\label{fig41}}
\end{figure}

\begin{figure}[ht!]
\centerline{\includegraphics[trim=1.0cm 0.4cm 1.4cm 0.3cm, clip, width=5.5cm]{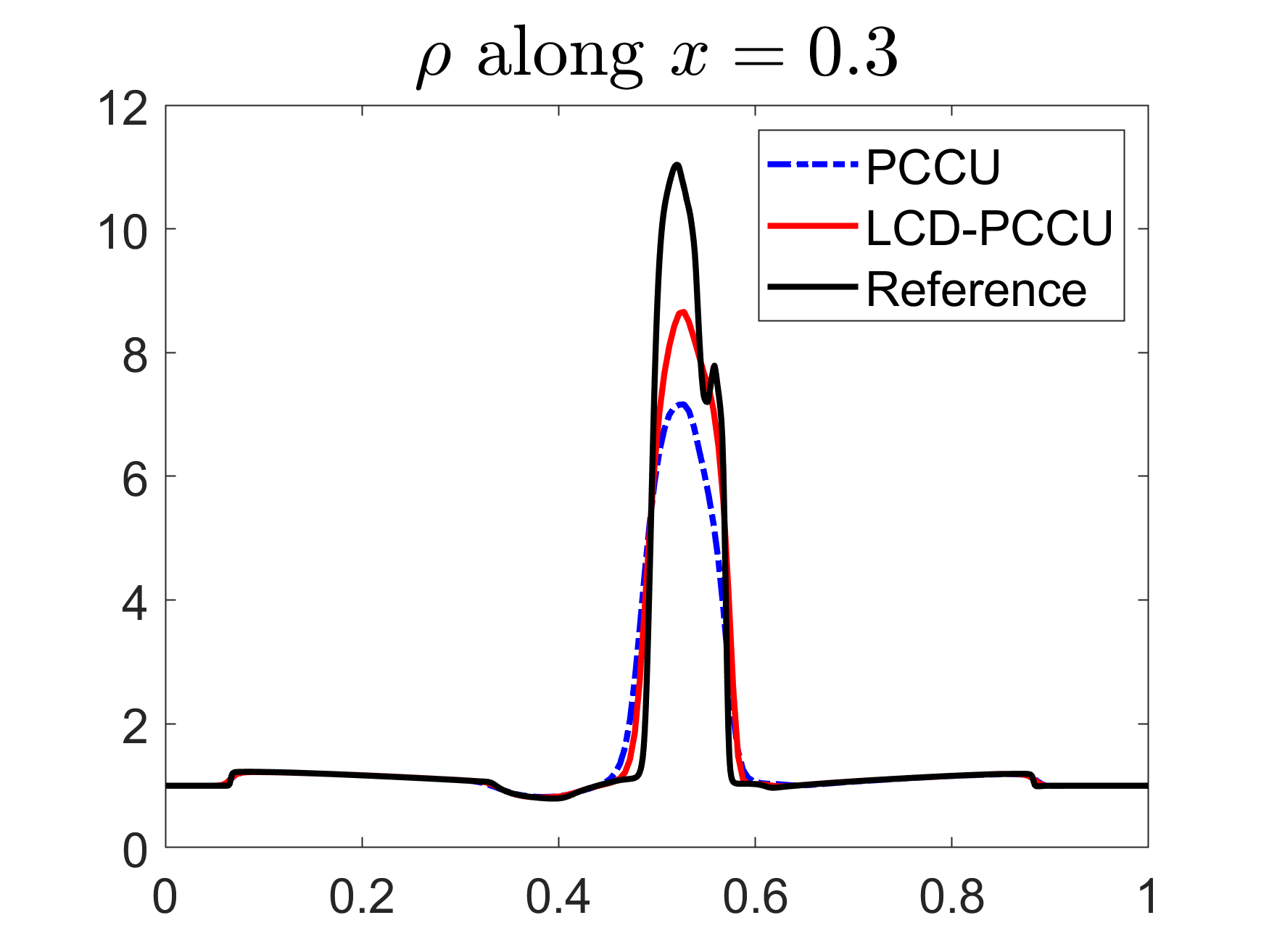}}
\caption{\sf Example 4: 1-D slices along the line $x=0.3$ of the densities from Figure \ref{fig41} together with the reference solution.
\label{fig41a}}
\end{figure}

We compute the numerical solutions by both the LCD-PCCU and PCCU schemes until the final time $t=0.295$ using a uniform $200\times200$ mesh
and plot the obtained $\rho$ and $p$ in Figure \ref{fig41}, where one can see that the LCD-PCCU scheme achieves higher resolution. To
further demonstrate this, we show (in Figure \ref{fig41a}) the 1-D slices of both densities along $x=0.3$ together with the reference
solution computed by the PCCU scheme on $1000\times1000$ uniform mesh.

\subsubsection*{Example 5---Blast Problem}
In the last example taken from \cite{BS1999}, we consider the blast problem, which poses a significant challenge due to the low gas pressure
and presence of strong magnetosonic shocks, which frequently lead to the occurrence of negative pressures near the shocks; see
\cite{LX2012,LXY2011} and references therein.

The initial conditions,
$$
(\rho,u,v,w,b_1,b_2,b_3)\Big|_{(x,y,0)}=\left(1,0,0,0,\frac{50}{\sqrt{\pi}},0,0\right),\quad
p(x,y,0)=\begin{cases}1000,&\sqrt{x^2+y^2}<0.1,\\0.1&\mbox{otherwise}, 
\end{cases} 
$$
are prescribed in the computational domain $[-0.5,0.5]\times[-0.5,0.5]$ subject to the zero-order extrapolation imposed at the boundary.

\begin{figure}[ht!]
\centerline{\includegraphics[trim=0.3cm 1.5cm 0.8cm 0.8cm, clip, width=5.5cm]{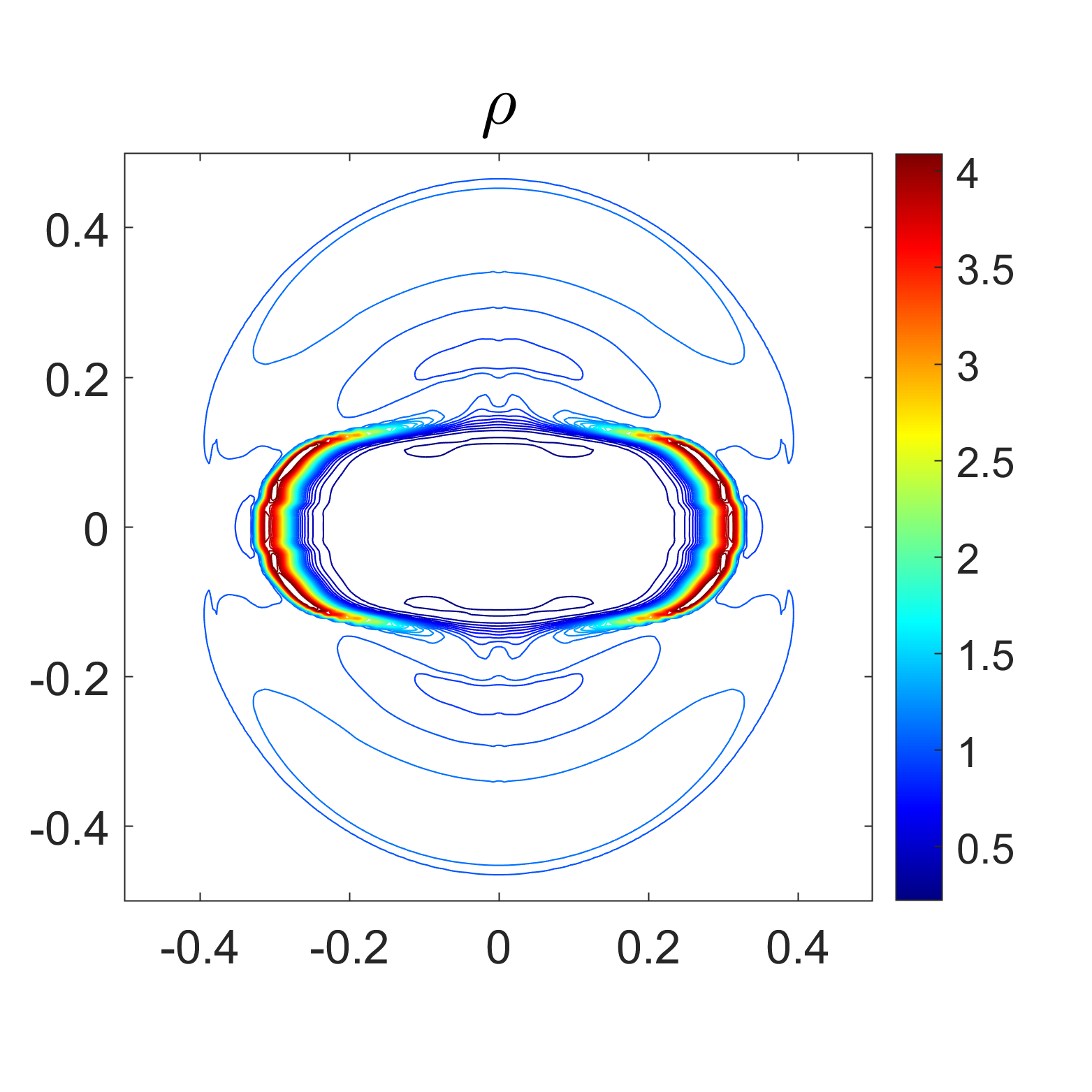}\hspace*{0.5cm}
            \includegraphics[trim=0.3cm 1.5cm 0.8cm 0.8cm, clip, width=5.5cm]{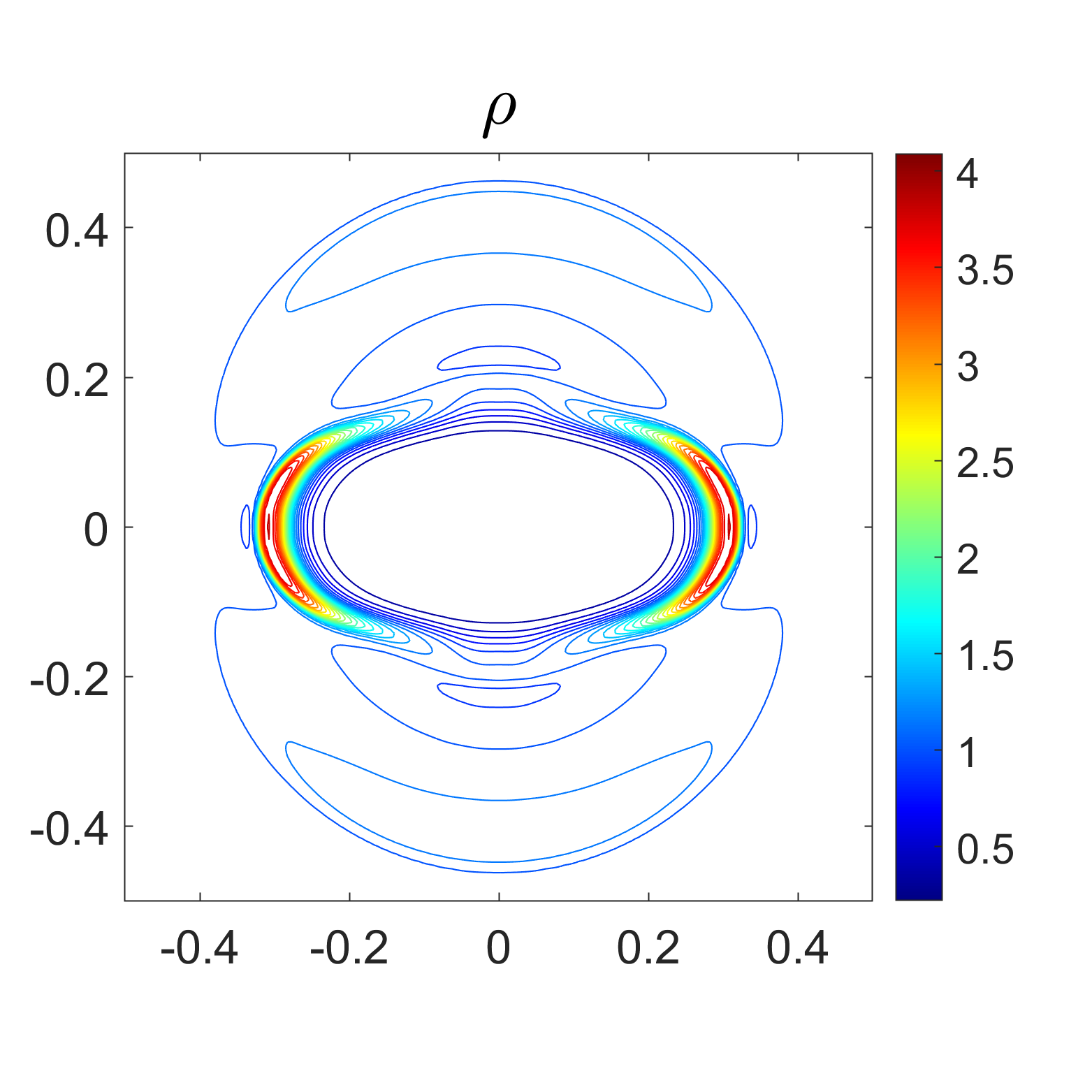}}
\vskip7pt
\centerline{\includegraphics[trim=0.3cm 1.5cm 0.8cm 0.8cm, clip, width=5.5cm]{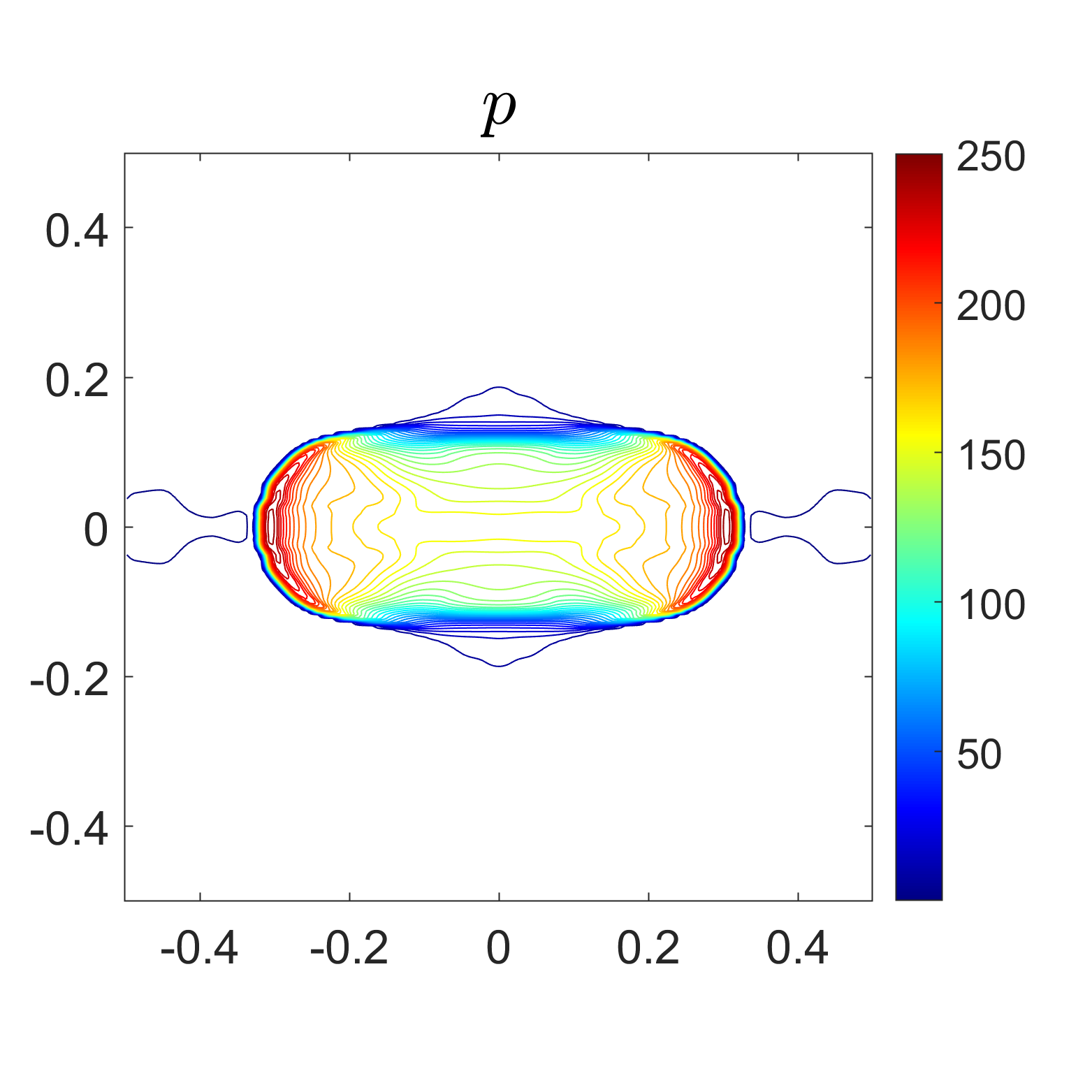}\hspace*{0.5cm}
            \includegraphics[trim=0.3cm 1.5cm 0.8cm 0.8cm, clip, width=5.5cm]{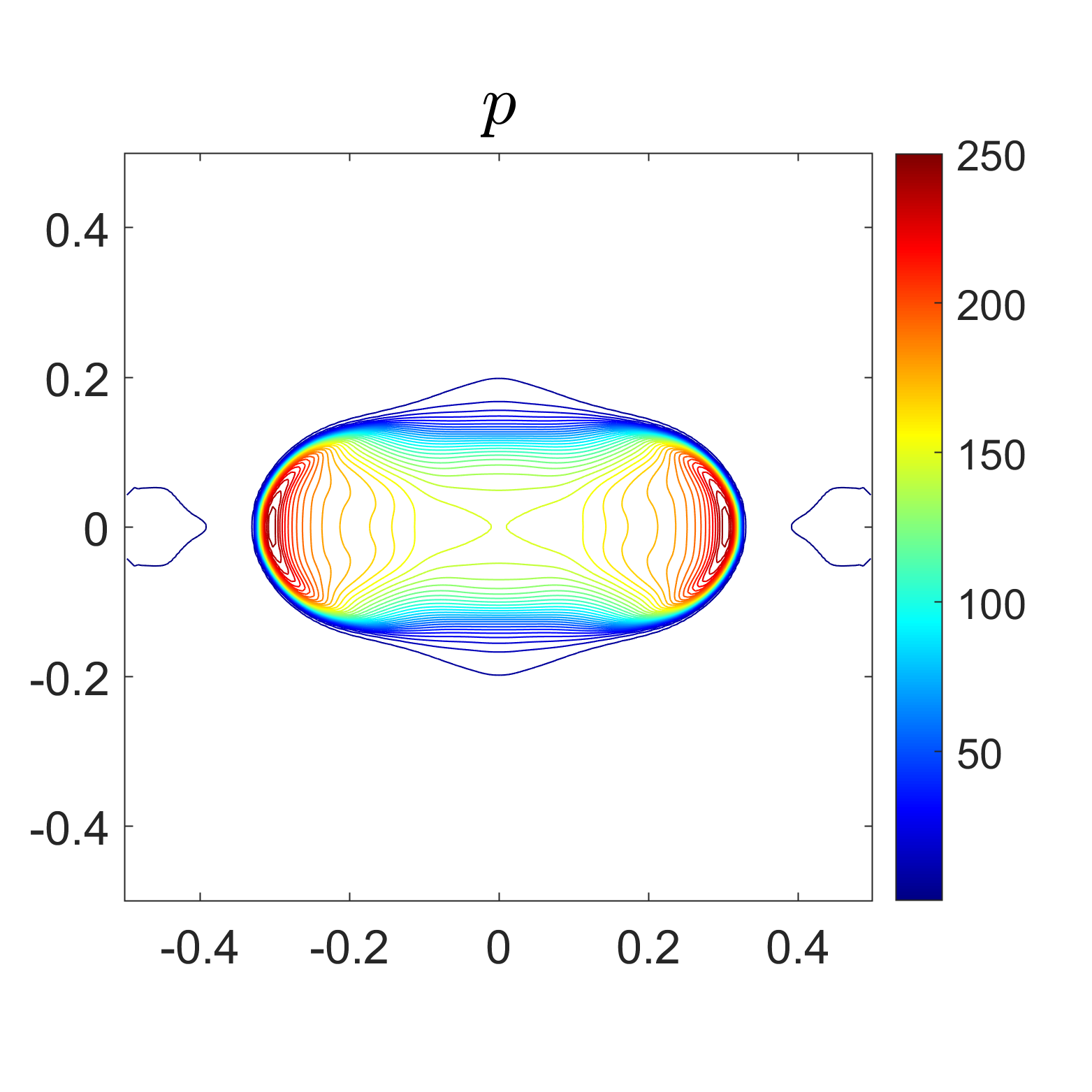}}
\vskip7pt
\centerline{\includegraphics[trim=0.3cm 1.5cm 0.8cm 0.8cm, clip, width=5.5cm]{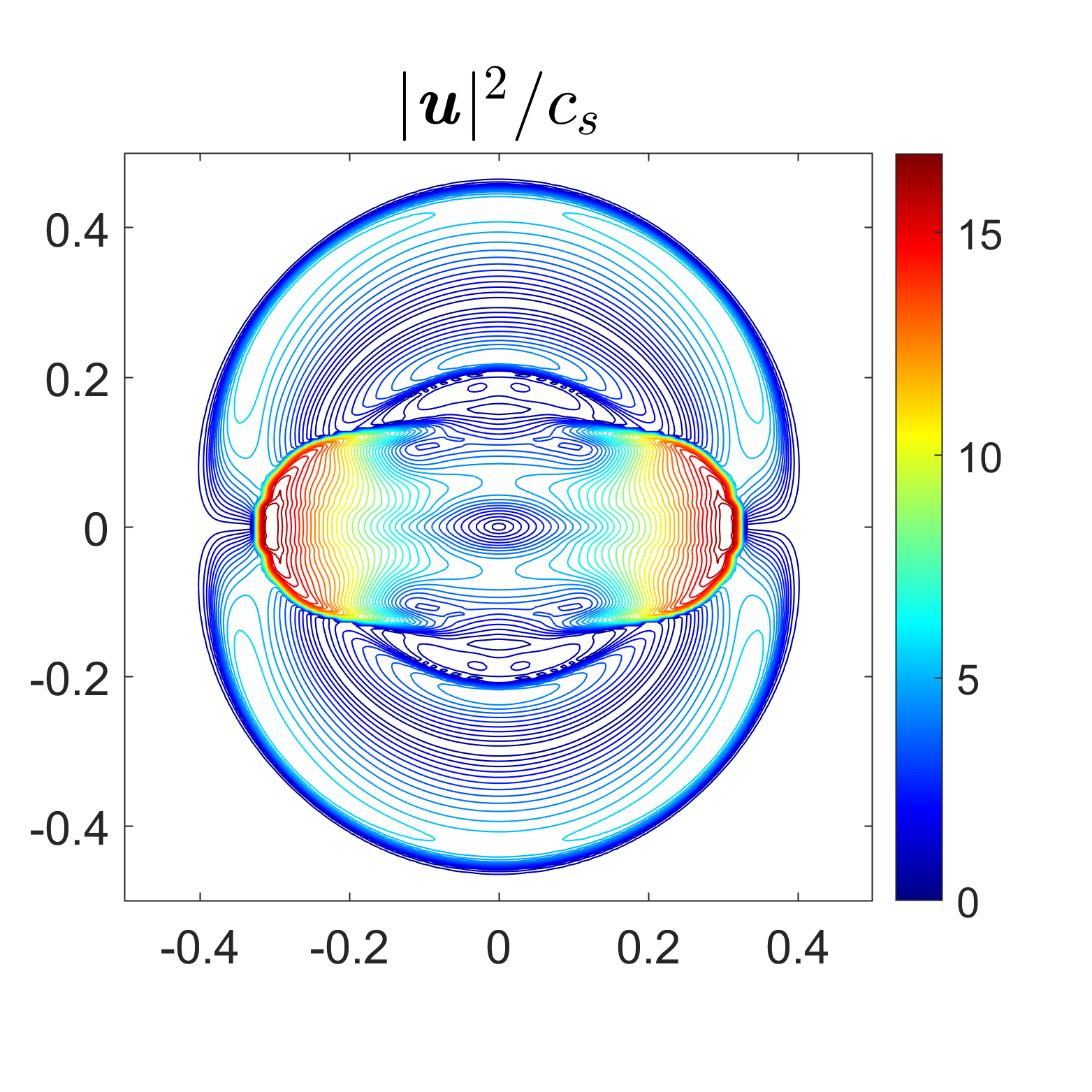}\hspace*{0.5cm}
            \includegraphics[trim=0.3cm 1.5cm 0.8cm 0.8cm, clip, width=5.5cm]{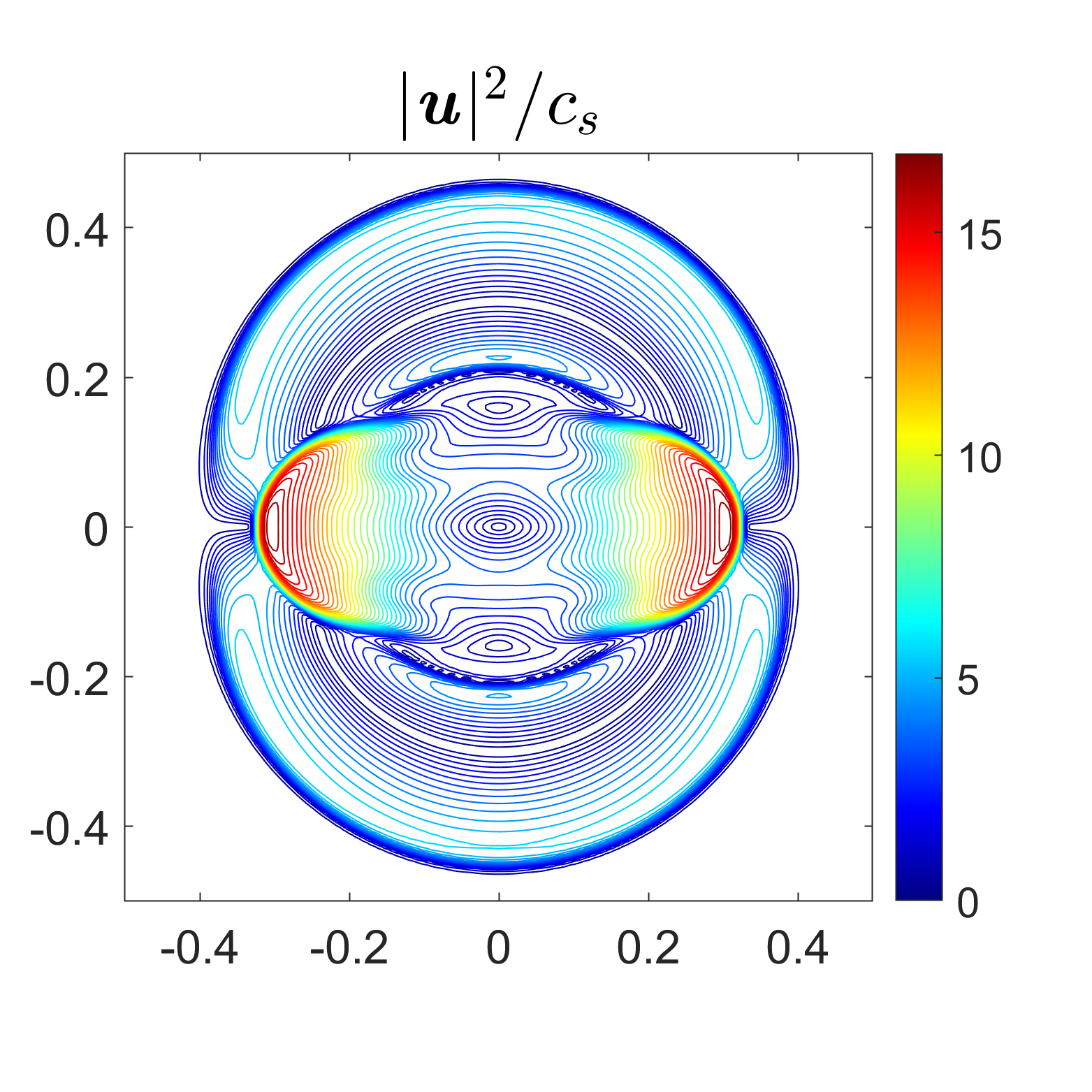}}
\vskip7pt
\centerline{\includegraphics[trim=0.3cm 1.5cm 0.8cm 0.8cm, clip, width=5.5cm]{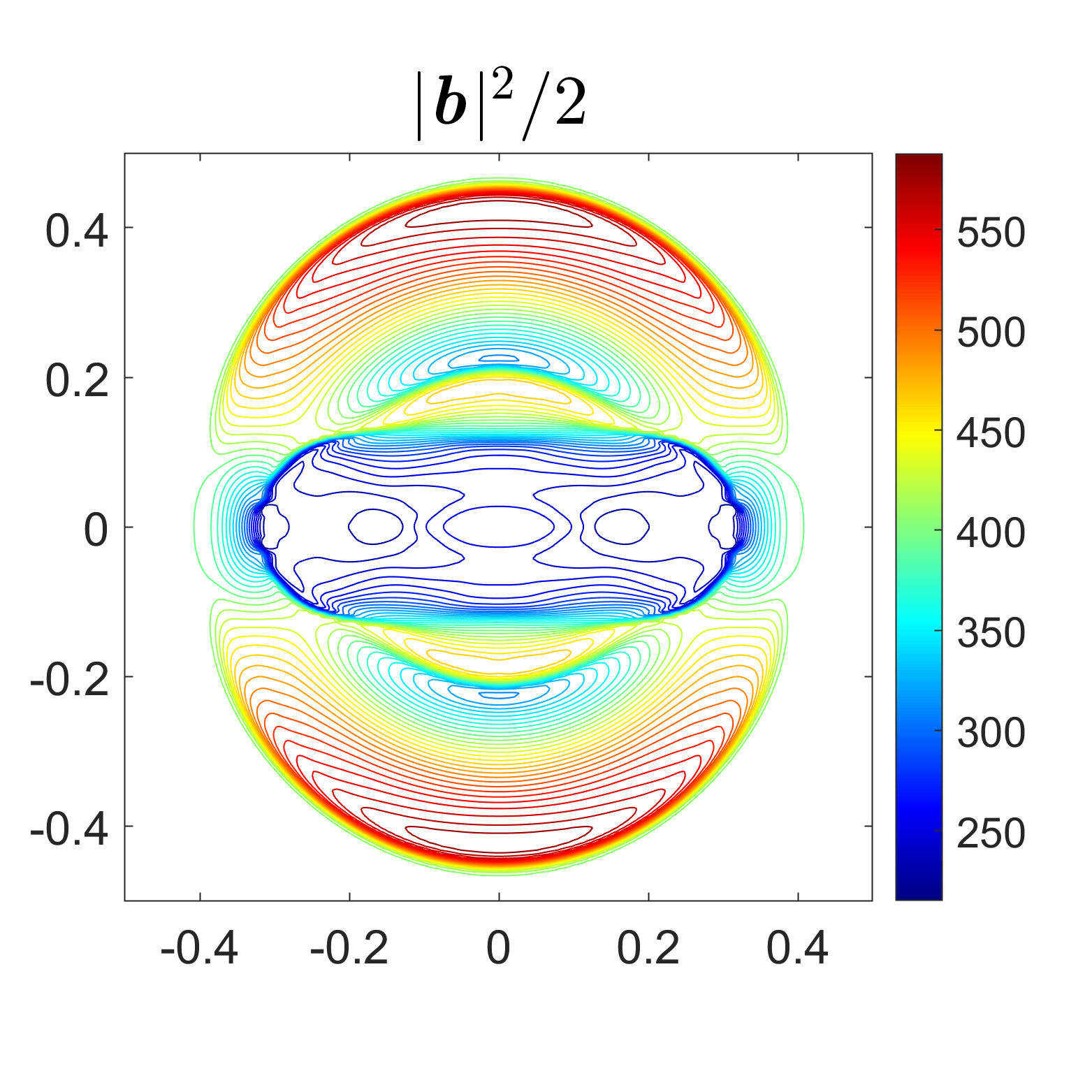}\hspace*{0.5cm}
            \includegraphics[trim=0.3cm 1.5cm 0.8cm 0.8cm, clip, width=5.5cm]{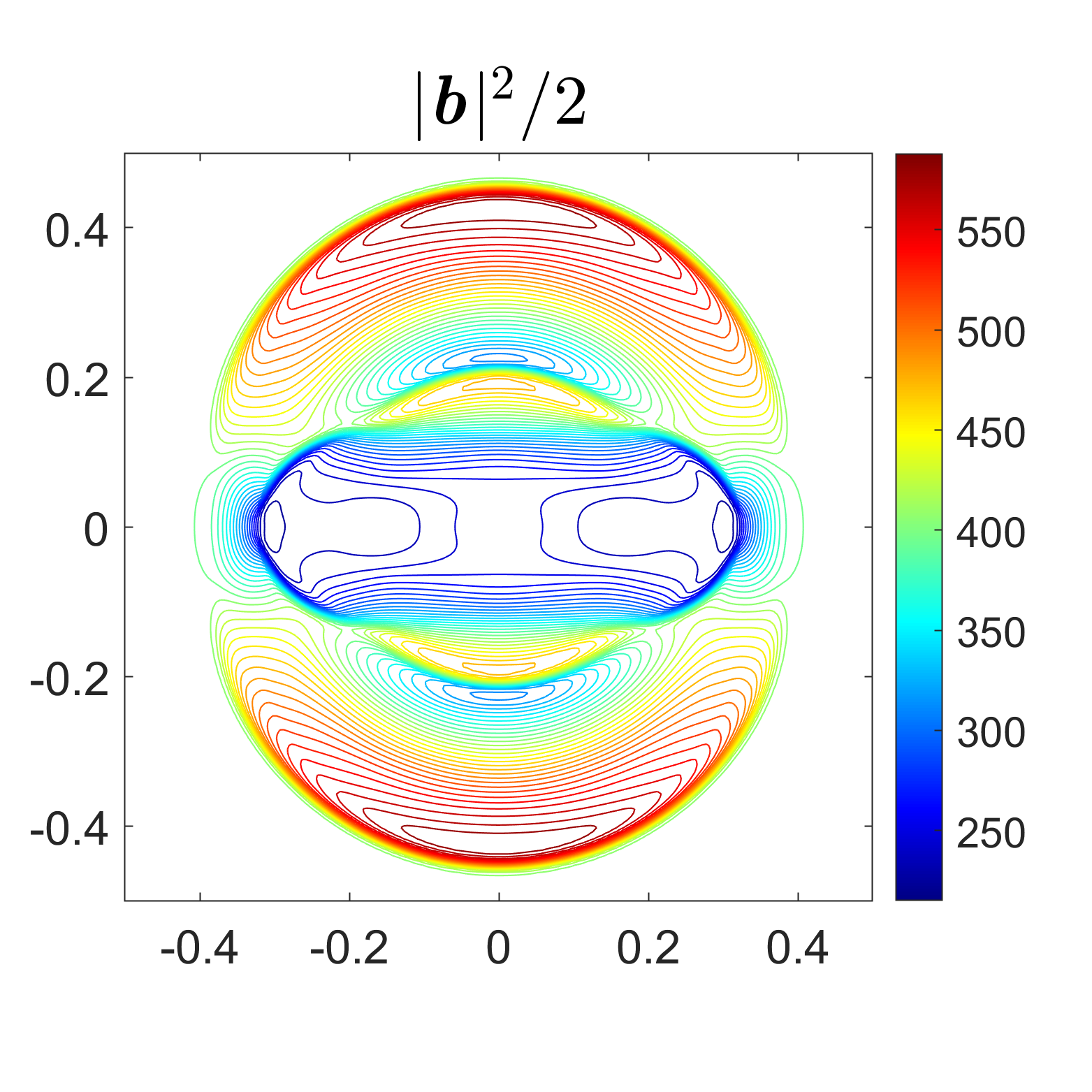}}
\caption{\sf Example 5: Density $\rho$ (top row), pressure $p$ (second row), velocity magnitude $|\bm u|$ (third row), and magnetic pressure
$|\bm b|^2/2$ (bottom row) computed by the LCD-PCCU (left column) and PCCU (right column) schemes on a uniform $200\times 200$ mesh.
\label{fig51}}
\end{figure}

\begin{figure}[ht!]
\centerline{\includegraphics[trim=0.3cm 1.5cm 0.8cm 0.8cm, clip, width=5.5cm]{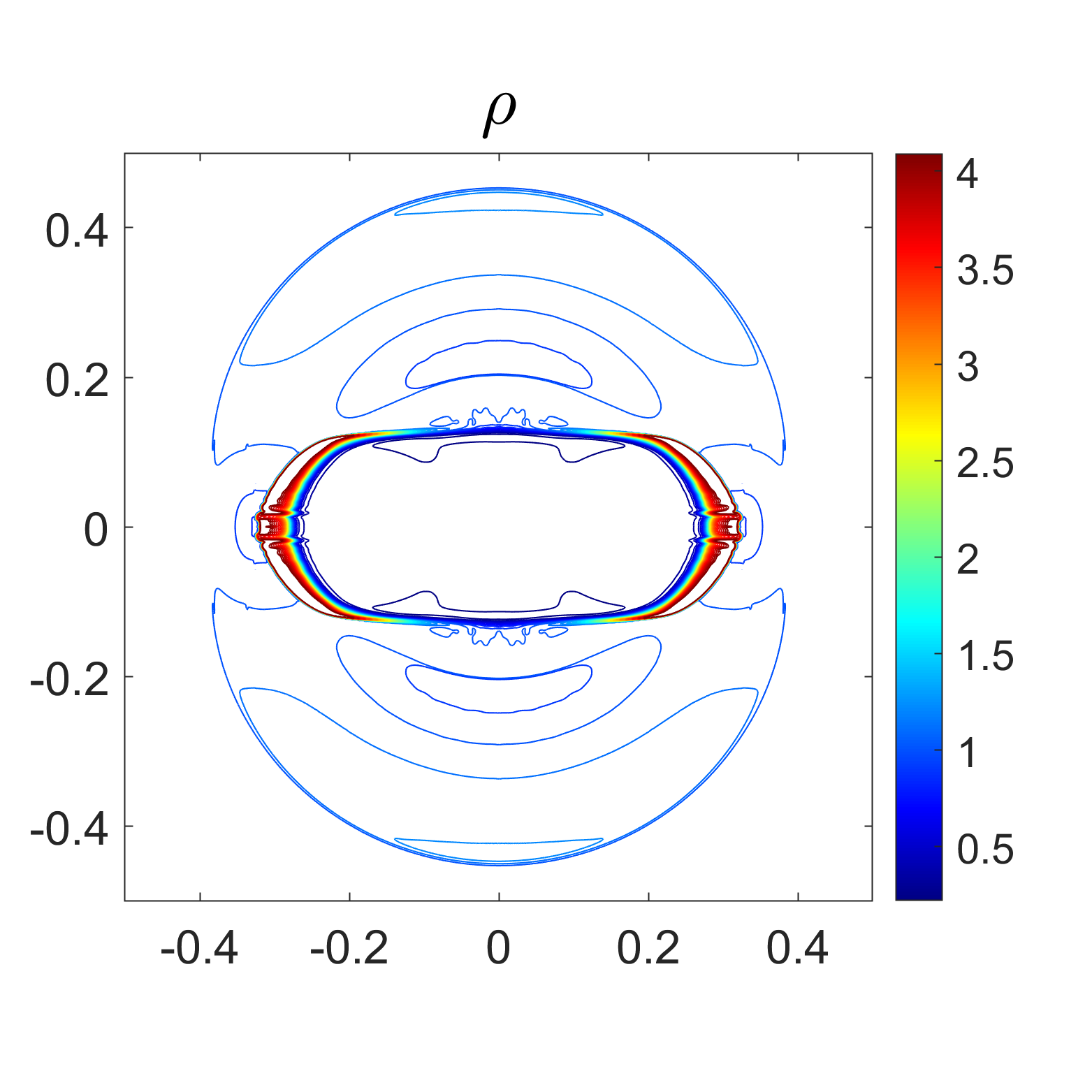}\hspace*{0.5cm}
            \includegraphics[trim=0.3cm 1.5cm 0.8cm 0.8cm, clip, width=5.5cm]{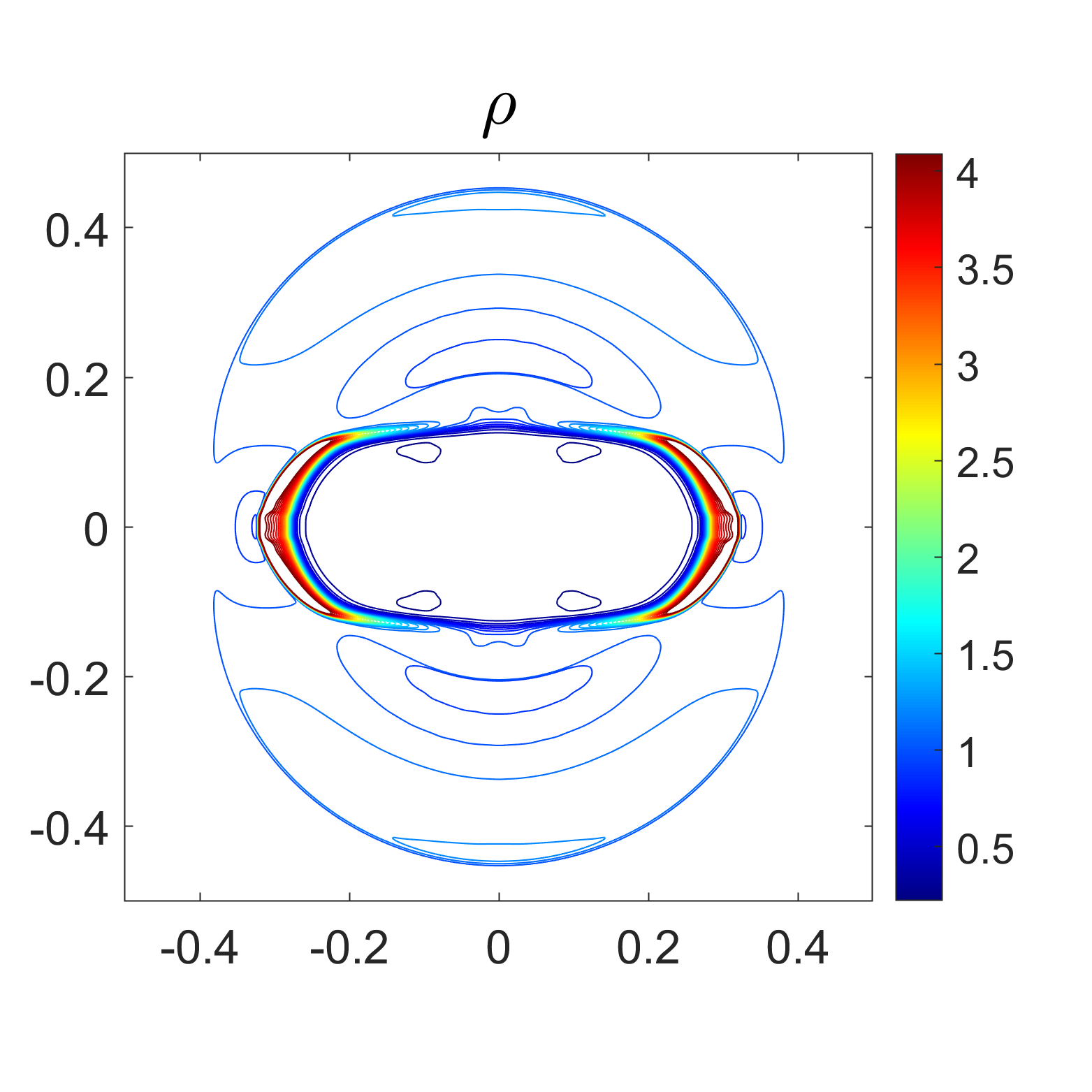}\hspace*{0.5cm}
            \includegraphics[trim=0.3cm 1.5cm 0.8cm 0.8cm, clip, width=5.5cm]{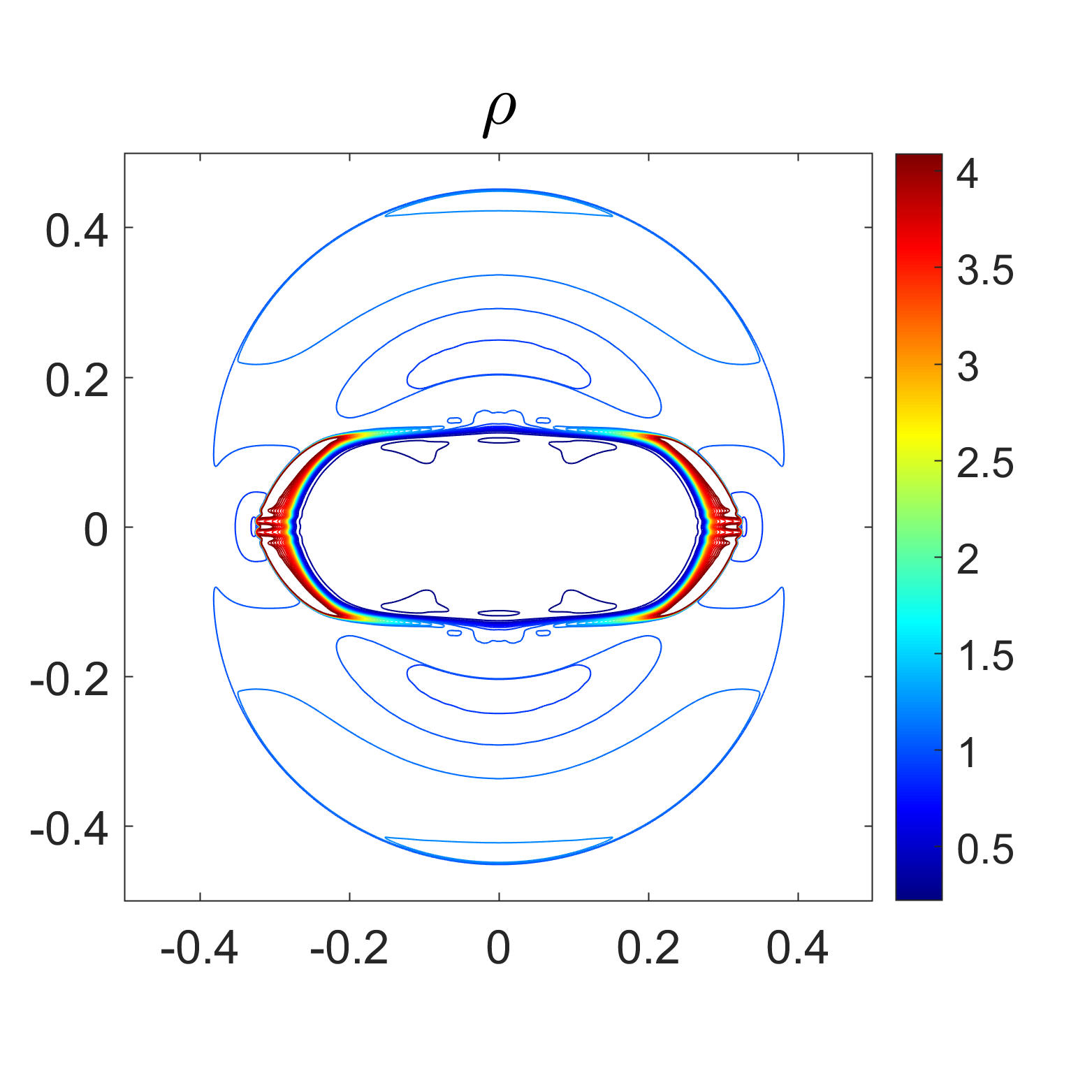}}
\vskip7pt
\centerline{\includegraphics[trim=0.3cm 1.5cm 0.8cm 0.8cm, clip, width=5.5cm]{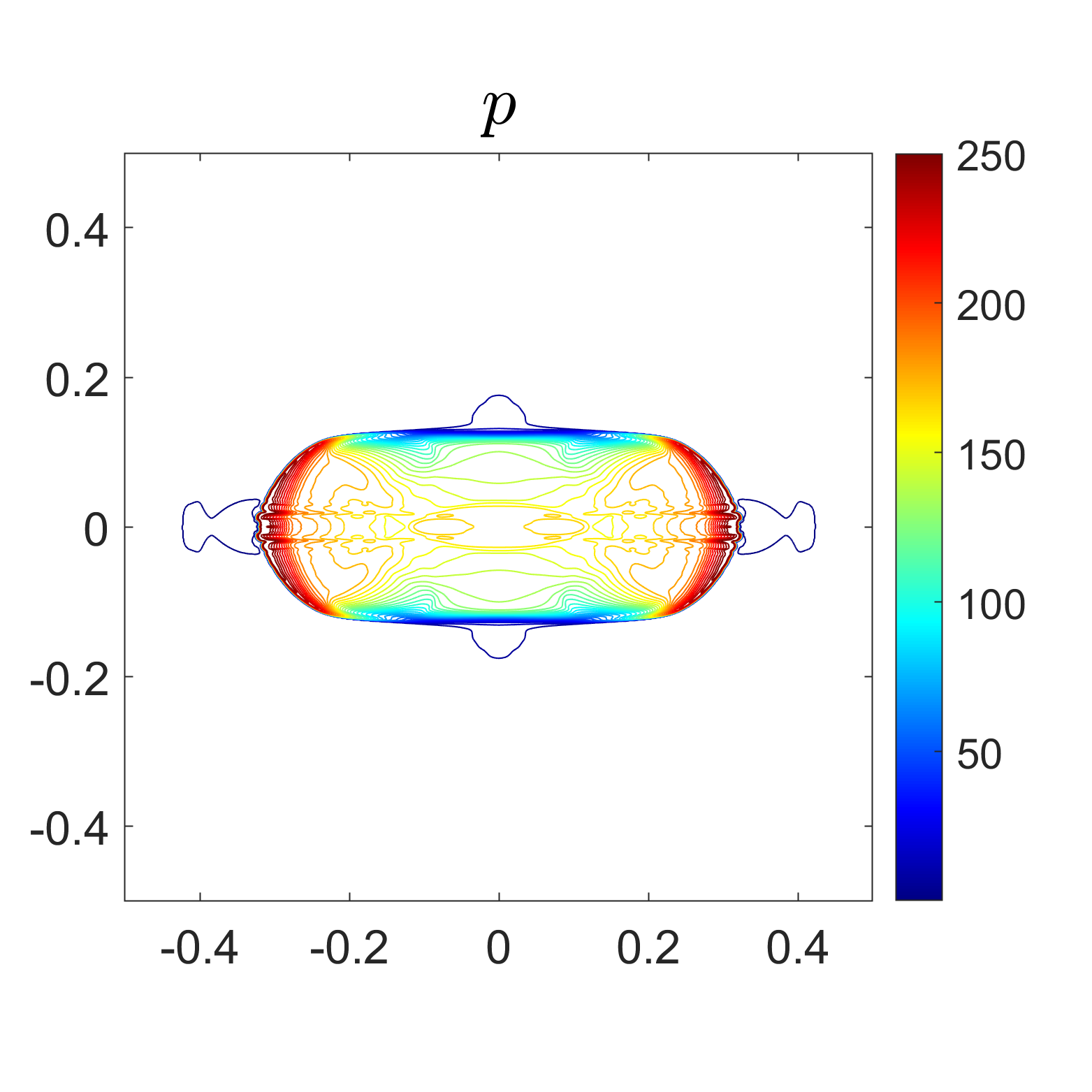}\hspace*{0.5cm}
            \includegraphics[trim=0.3cm 1.5cm 0.8cm 0.8cm, clip, width=5.5cm]{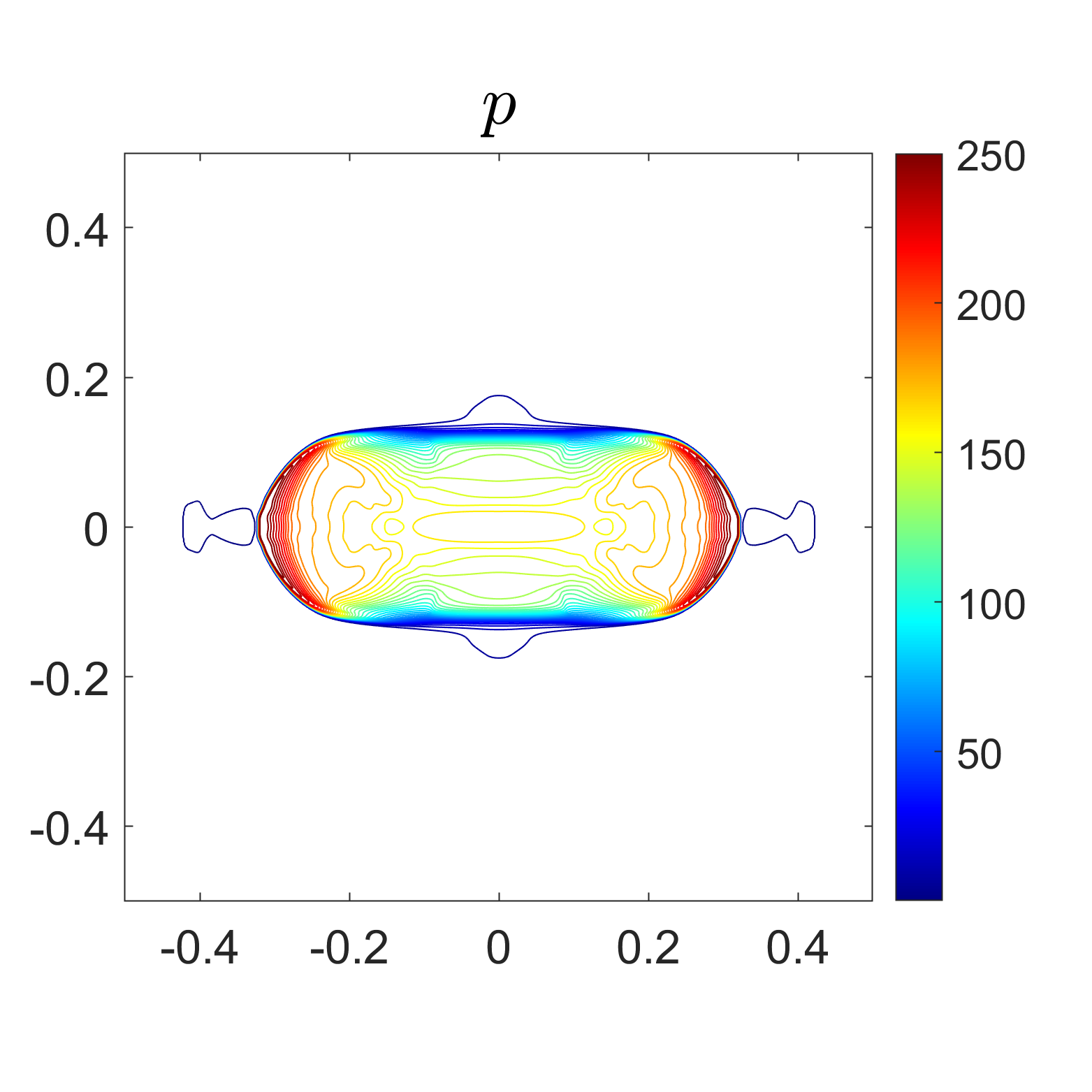}\hspace*{0.5cm}
            \includegraphics[trim=0.3cm 1.5cm 0.8cm 0.8cm, clip, width=5.5cm]{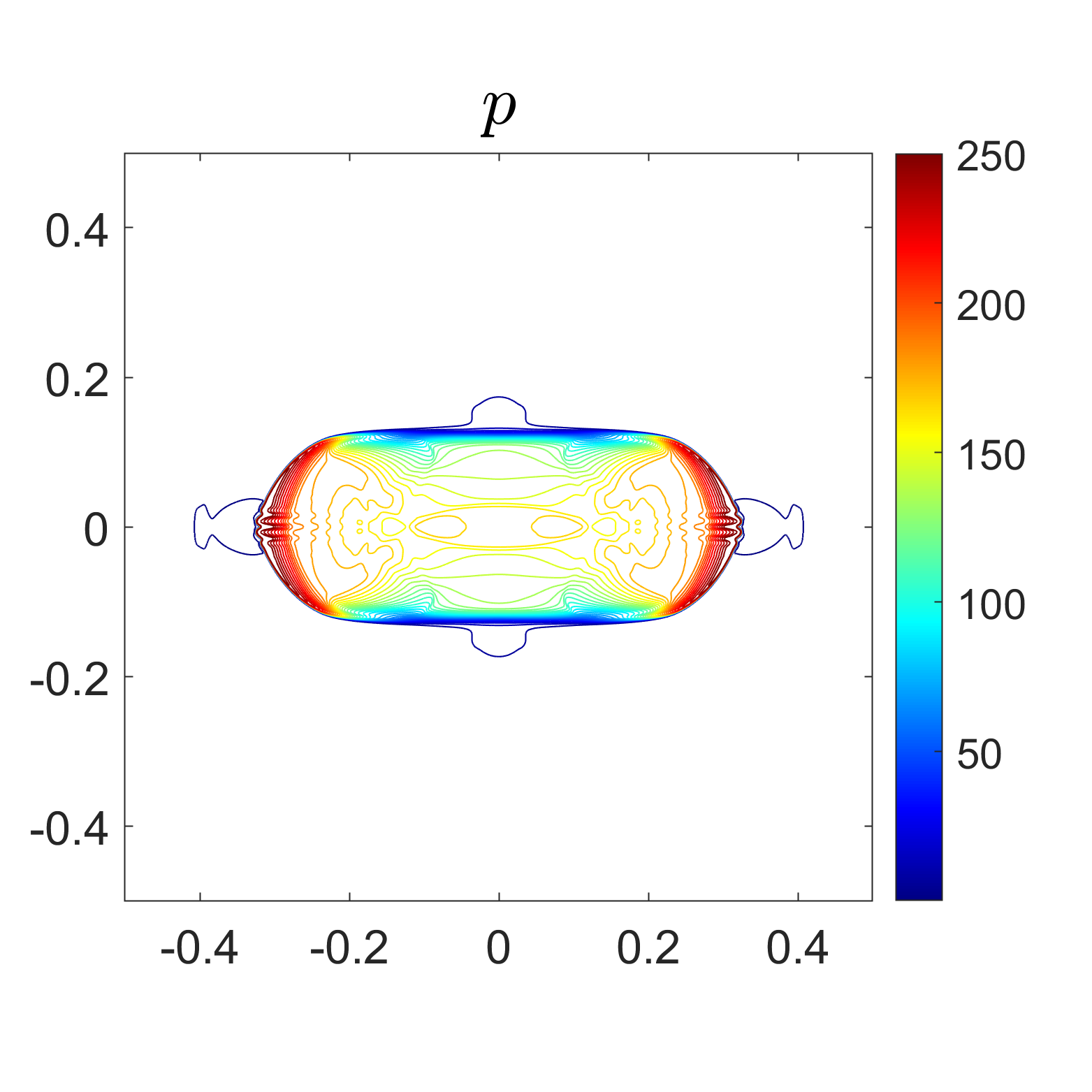}}
\vskip7pt
\centerline{\includegraphics[trim=0.3cm 1.5cm 0.8cm 0.8cm, clip, width=5.5cm]{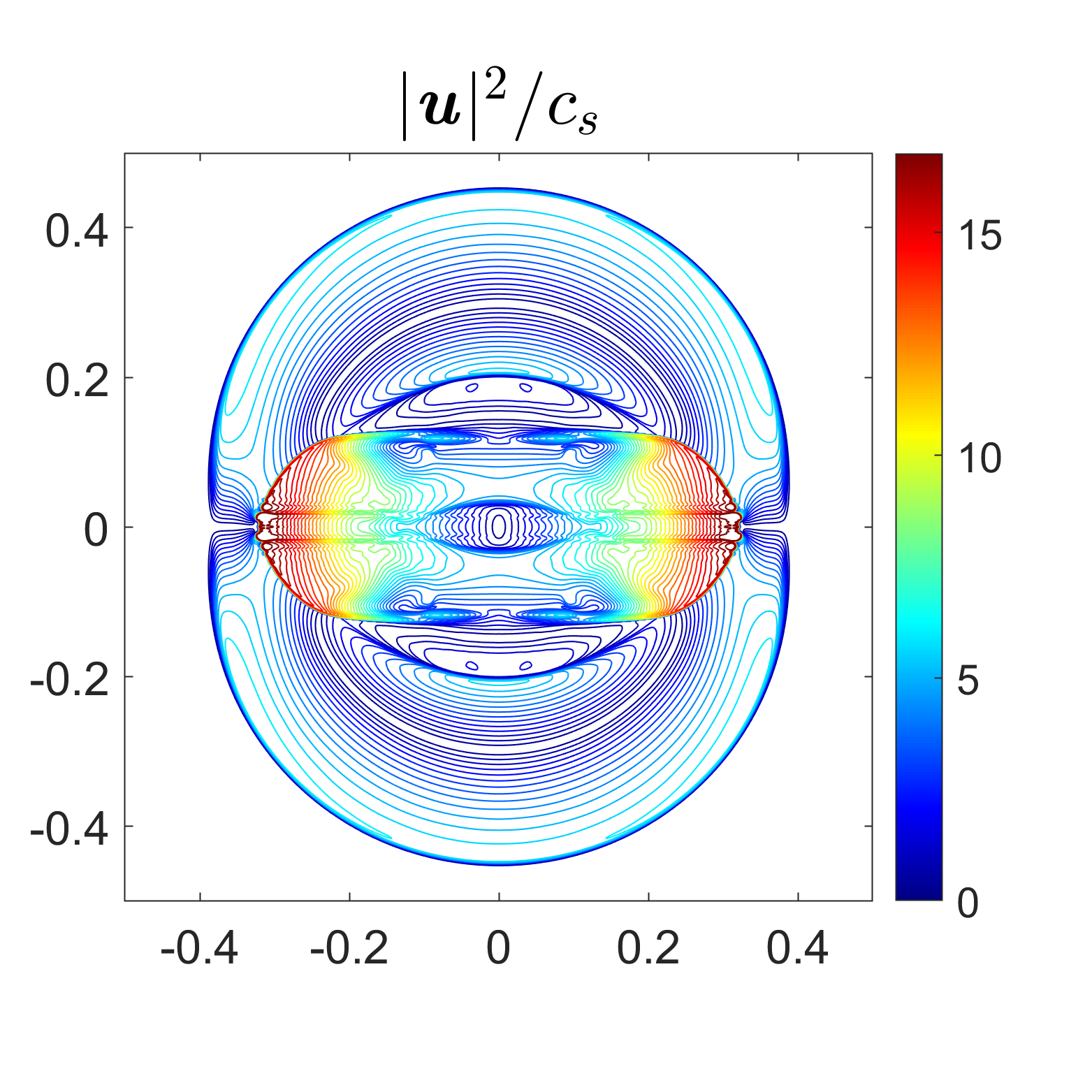}\hspace*{0.5cm}
            \includegraphics[trim=0.3cm 1.5cm 0.8cm 0.8cm, clip, width=5.5cm]{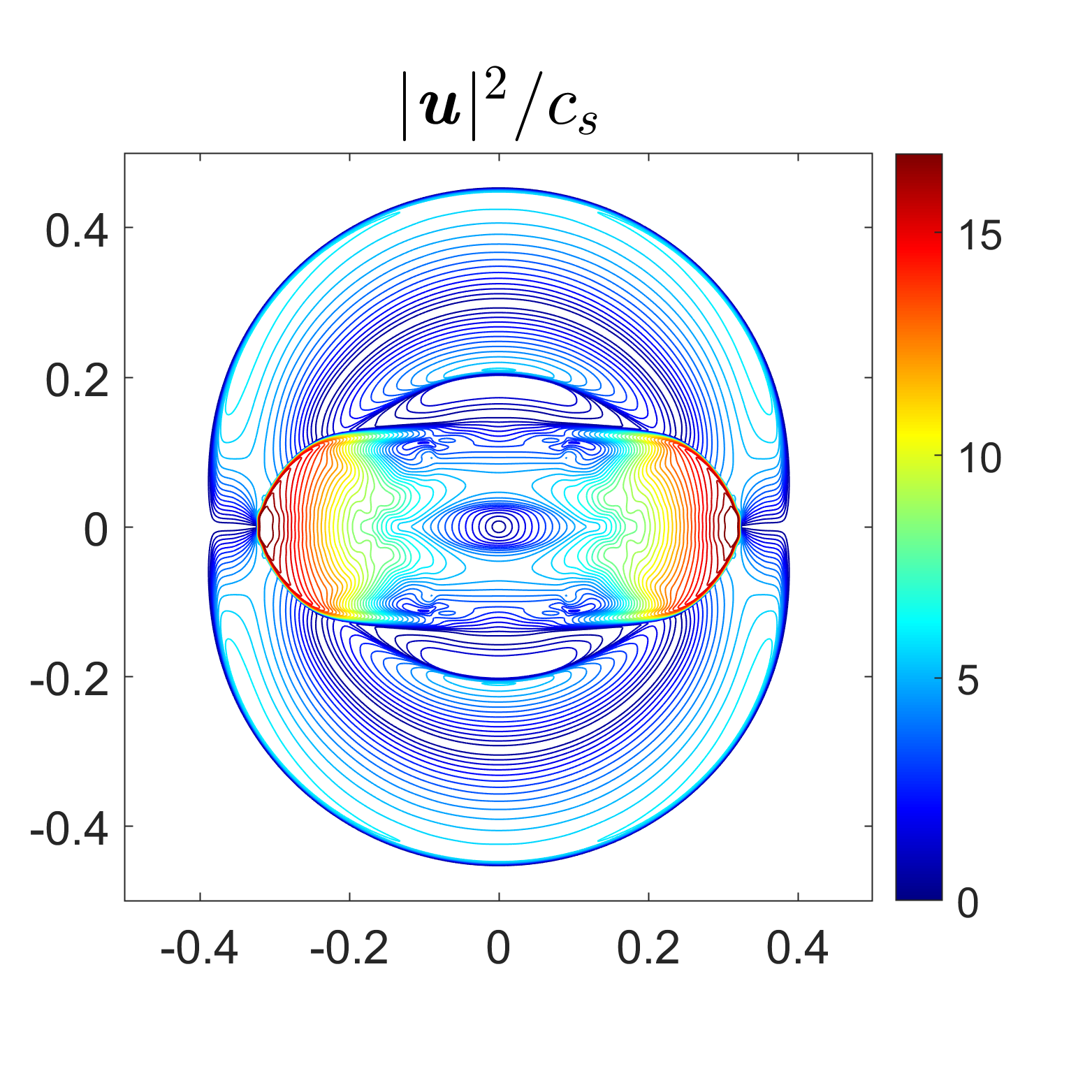}\hspace*{0.5cm}
            \includegraphics[trim=0.3cm 1.5cm 0.8cm 0.8cm, clip, width=5.5cm]{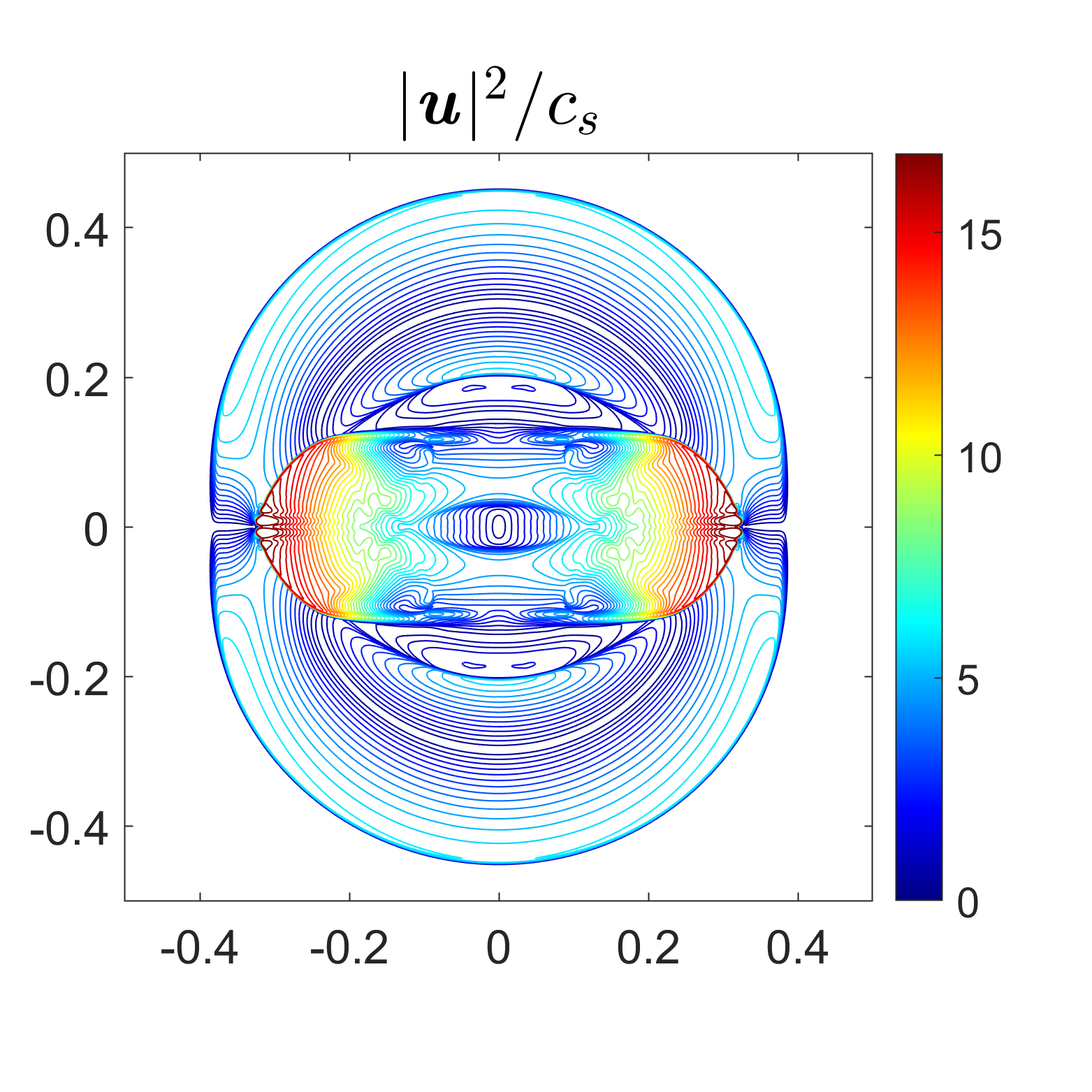}}
\vskip7pt
\centerline{\includegraphics[trim=0.3cm 1.5cm 0.8cm 0.8cm, clip, width=5.5cm]{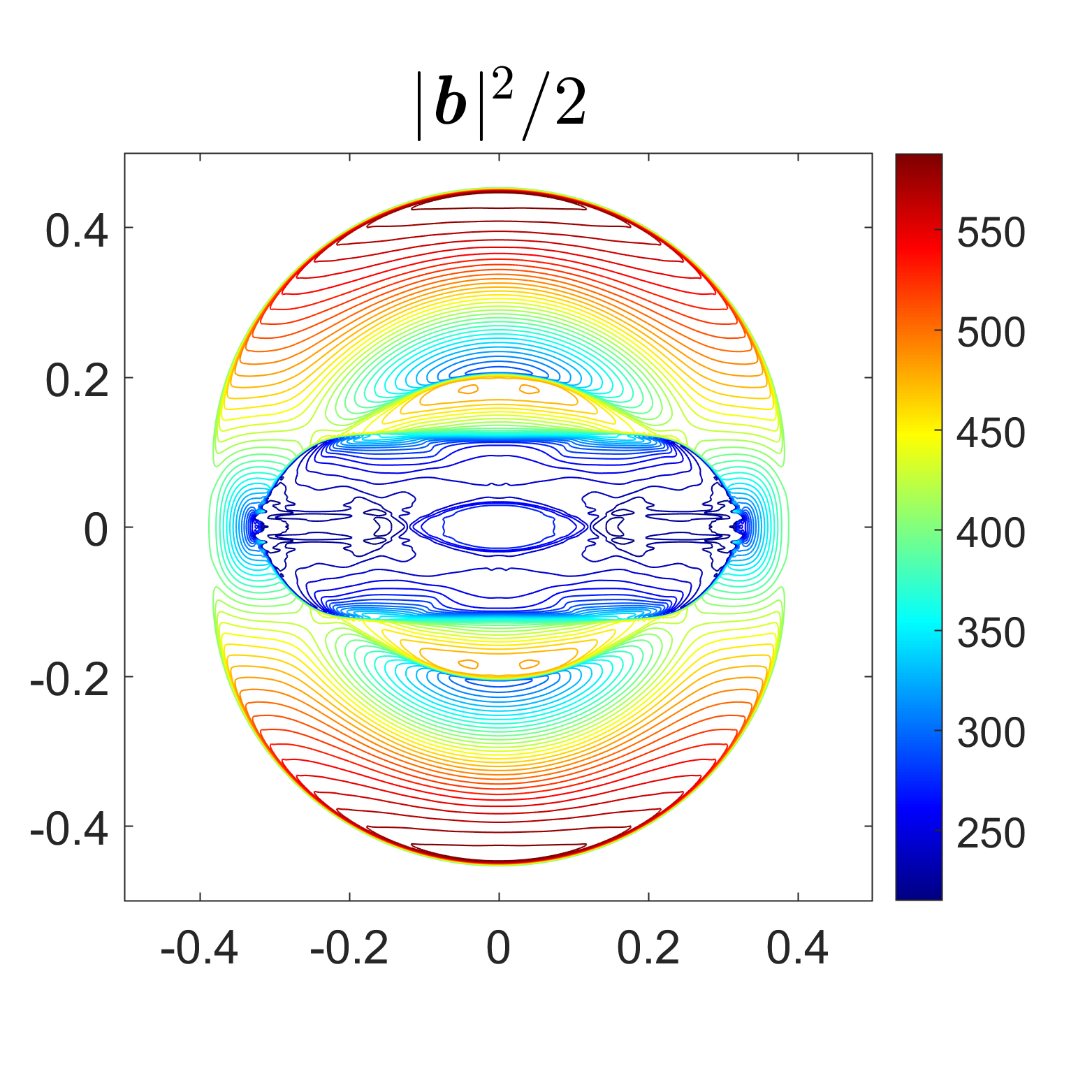}\hspace*{0.5cm}
            \includegraphics[trim=0.3cm 1.5cm 0.8cm 0.8cm, clip, width=5.5cm]{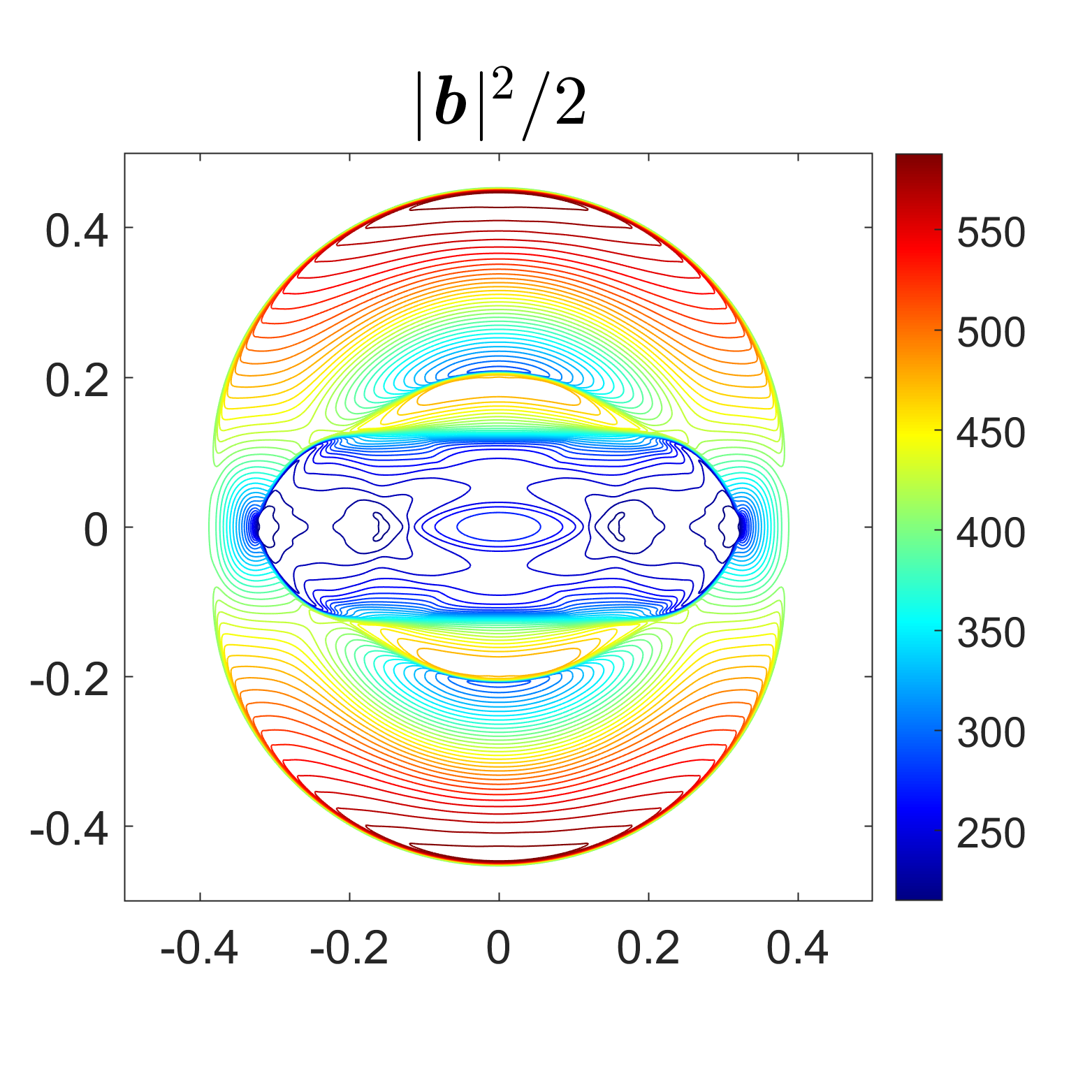}\hspace*{0.5cm}
            \includegraphics[trim=0.3cm 1.5cm 0.8cm 0.8cm, clip, width=5.5cm]{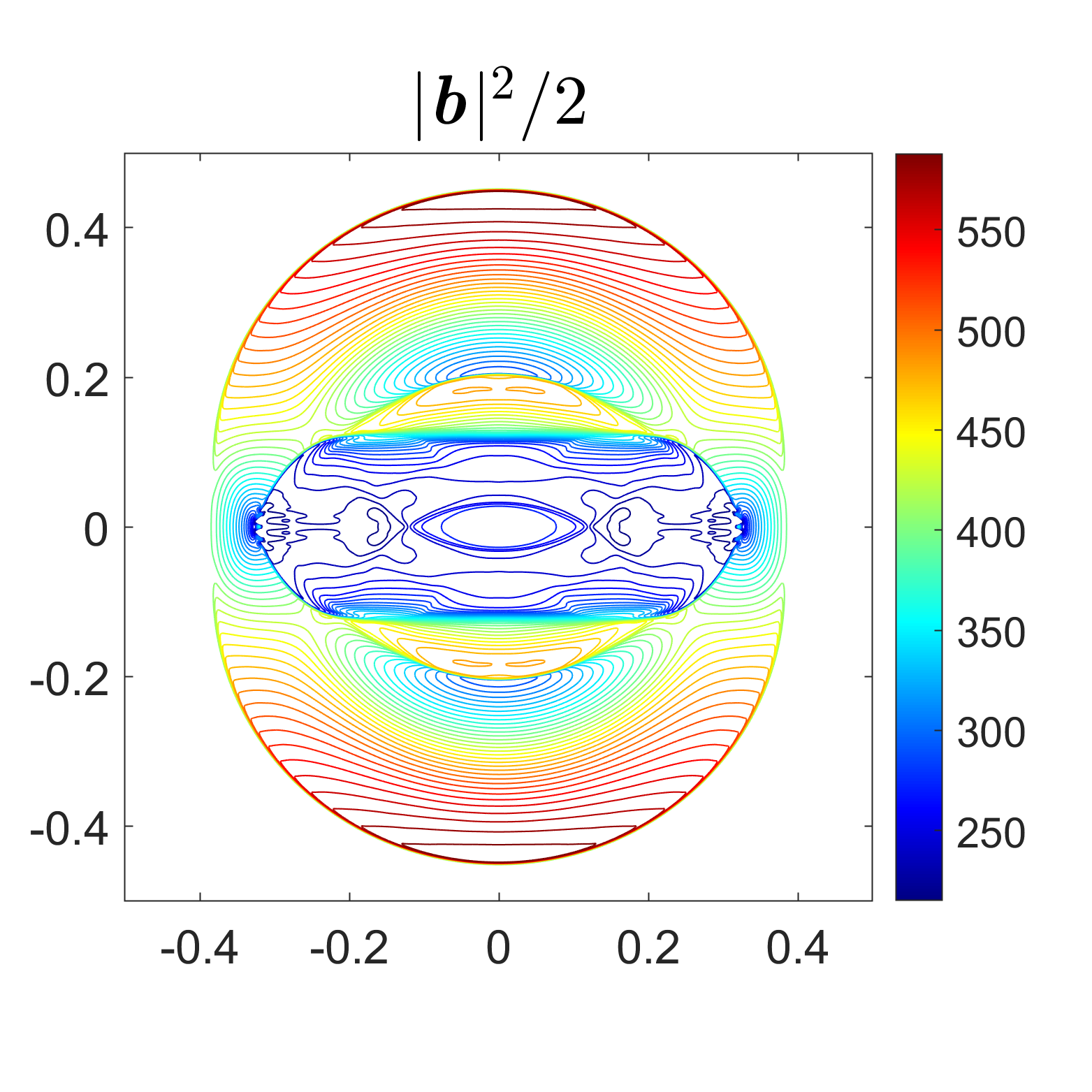}}
\caption{\sf Example 5: Density $\rho$ (top row), pressure $p$ (second row), velocity magnitude $|\bm u|$ (third row), and magnetic pressure
$|\bm b|^2/2$ (bottom row) computed by the LCD-PCCU scheme on a uniform $1000\times1000$ mesh (left column) and PCCU scheme on uniform
$1000\times1000$ (middle column) and $2000\times2000$ (right column) meshes.\label{fig52}}
\end{figure}

We compute the numerical solutions by both the LCD-PCCU and PCCU schemes until the final time $t=0.01$ on a uniform $200\times200$ mesh and
plot the obtained density $\rho$, pressure $p$, velocity magnitude $|\bm u|$, and magnetic pressure $|\bm b|^2/2$ in Figure \ref{fig51}. As
one can see, the numerical results computed by the LCD-PCCU scheme are visibly sharper. However, they contain wiggles in the areas of high
density and pressure. To experimentally verify that these structures are not numerical artifacts, we refine the mesh and perform the same
computations on a uniform $1000\times1000$ mesh. The obtained results reported in the left two columns of Figure \ref{fig52}, indicate that
similar structures start developing in the PCCU solution as well. We thus further refine the mesh and run the PCCU simulation on an even
finer $2000\times2000$ mesh; see Figure \ref{fig52} (right column). One can observe that those wiggly structures are now clearly present in
the PCCU results. One can notice that the resolution achieved by the LCD-PCCU scheme on the $1000\times1000$ mesh is practically the same as
that achieved by the PCCU scheme on the $2000\times2000$ mesh. This clearly indicates an advantage of the proposed LCD-PCCU scheme.

\section{Conclusions}
In this paper, we have developed a locally divergence-free local characteristic decomposition (LCD) based path-conservative central-upwind
(LCD-PCCU) scheme for the ideal magnetohydrodynamics (MHD) equations. The proposed scheme is applied to the Godunov-Powell nonconservative
modifications of the studied MHD systems, which have a complete eigenstructure required to derive LCD-based central-upwind numerical fluxes;
see \cite{CCHKL_22,CHK25}. In order to ensure the local divergence-free property, we have followed \cite{Chertock24} and augmented the
studied systems with the evolution equations for the corresponding derivatives of the magnetic field components and by using these evolved
quantities in the design of a special piecewise linear reconstruction of the magnetic field, which also guarantees a non-oscillatory nature
of the resulting scheme. The designed LCD-PCCU scheme has been tested on several benchmarks, and the obtained numerical results demonstrate
that the proposed scheme outperforms its PCCU counterpart from \cite{Chertock24}.

\appendix
\section{Eigendecomposition for Conservative Variables}\label{appA1}
In this appendix, we provide the reader with the matrices used in the LCD of the quasi-linear system \eref{2.7}; see \cite{Brio88} for
details.
 
First, the matrices $C^x$ and $C^y$ are 
\begin{equation}
C^x(\mU)=\begin{pmatrix}0&1&0&0&0&0&0&0\\a_1&\gamma_3u_N&\gamma_1u_T&\gamma_1w&-\gamma_1&b_1&\gamma_2b_2&\gamma_2b_3\\
-u_Nu_T&u_T&u_N&0&0&b_2&-b_N&0\\-u_Nw&w&0&u_N&0&b_3&0&-b_N\\a_2&a_3&a_4&a_5&\gamma u_N&\bm u\!\cdot\!\bm b&a_6&a_7\\0&0&0&0&0&u&0&0\\
\dfrac{u_Tb_N-u_Nb_T}{\rho}&\dfrac{b_T}{\rho}&-\dfrac{b_N}{\rho}&0&0&v&u_N&0\\[1.8ex]
\dfrac{wb_N-u_Nb_3}{\rho}&\dfrac{b_3}{\rho}&0&-\dfrac{b_N}{\rho}&0&w&0&u_N\end{pmatrix},
\label{A1}
\end{equation}
and 
\begin{equation}
C^y(\mU)=\begin{pmatrix}0&0&1&0&0&0&0&0\\-u_Nu_T&u_N&u_T&0&0&-b_N&b_1&0\\
a_1&\gamma_1u_T&\gamma_3u_N&\gamma_1w&-\gamma_1&\gamma_2b_2&b_2&\gamma_2b_3\\-u_Nw&0&w&u_N&0&0&b_3&-b_N\\
a_2&a_4&a_3&a_5&\gamma u_N&a_6&\bm u\!\cdot\!\bm b&a_7\\\dfrac{u_Tb_N-u_Nb_T}{\rho}&-\dfrac{b_N}{\rho}&\dfrac{b_T}{\rho}&0&0&u_N&u&0\\
0&0&0&0&0&0&v&0\\\dfrac{wb_N-u_Nb_3}{\rho}&0&\dfrac{b_3}{\rho}&-\dfrac{b_N}{\rho}&0&0&w&u_N\end{pmatrix},
\label{A2}
\end{equation}
where $u_N$, $b_N$ and $u_T$, $b_T$ are the normal and tangential components of $\bm u$ and $\bm b$ with respect to the $x$- and $y$-axis,
namely,
$$
(u_N,u_T,b_N,b_T):=\left\{\begin{aligned}
&(u,v,b_1,b_2)&&\mbox{in the $x$-direction (in \eref{A1})},\\
&(v,u,b_2,b_1)&&\mbox{in the $y$-direction (in \eref{A2})} ,
\end{aligned}\right.
$$
$\gamma_n:=n-\gamma$ for $n=1,2,3,$ and 
$$
\begin{aligned}
&a_1:=-\frac{\gamma_3}{2}u_N^2-\frac{\gamma_1}{2}(u_T^2+w^2),&&a_2:=
u_N\left(\frac{\gamma_1}{2\rho}\Big({\cal E}+p+\hf|\bm b|^2\Big)|\bm u|^2-\frac{b_N}{\rho}(\bm u\!\cdot\!\bm b)\right),\\
&a_3:=\frac{1}{\rho}\Big({\cal E}+p+\hf|\bm b|^2\Big)-\frac{b_N^2}{\rho}+\gamma_1u_N^2,&&a_4:=\gamma_1u_Nu_T-\frac{1}{\rho}b_Nb_T,\quad
a_5:=\gamma_1u_Nw-\frac{1}{\rho}b_Nb_3,\\
&a_6:=\gamma_2u_Nb_T-u_Tb_N,&&a_7:=\gamma_2u_Nb_3-wb_N.
\end{aligned}
$$

One can show that the eigenvalues of both $C^x$ and $C^y$ are given by
\begin{equation}
\lambda_{1,8}=u_N\mp c_f,\quad\lambda_{2,7}=u_N\mp c_a,\quad\lambda_{3,6}=u_N\mp c_s,\quad\lambda_4=\lambda_5=u_N,
\label{1.3a}
\end{equation}
where 
\begin{equation}
\begin{aligned}
c_a=\sqrt{\frac{b^2_N}{\rho}},\quad c_{f,s}=\left[\hf\left(c^2+\frac{|\bm b|^2}{\rho}\pm\sqrt{\Big(c^2+\frac{|\bm b|^2}{\rho}\Big)^2-
4c^2\frac{b_N^2}{\rho}}\right)\right]^\hf,
\end{aligned}
\label{A4}
\end{equation}
and $c:=\sqrt{\gamma p/\rho}$ is the speed of sound.

Finally, we provide the formula for the matrices $R^x$ and $R^y$, which diagonalize $C^x$ and $C^y$, that is, $(R^x)^{-1}C^xR^x$ and
$(R^y)^{-1}C^yR^y$ are diagonal. One can show that $R^x=\big[\bm r^x_1|\bm r^x_2|\dots|\bm r^x_7|\bm r^x_8\big]$ and
$(R^x)^{-1}=\big[\bm\ell^x_1|\bm\ell^x_2|\dots|\bm\ell^x_7|\bm\ell^x_8\big]^\top$, where $\bm r^x_i$ and $\bm\ell^x_i$, $i=1,\dots,8$ are
the right and left eigenvectors of $C^x$:
\allowdisplaybreaks
\begin{align*}
&\bm r^x_{1,8}=\begin{pmatrix}\alpha_f\\\alpha_f(u_N\mp c_f)\\\alpha_fu_T\pm\alpha_s\beta_1\beta_2c_a\\
\alpha_fw\pm\alpha_s\beta_1\beta_3c_a\\\hf\alpha_f|\bm u|^2+\mu_f^\mp\\0\\\frac{1}{\sqrt{\rho}}\alpha_s\beta_2c_f\\
\frac{1}{\sqrt{\rho}}\alpha_s\beta_3c_f\end{pmatrix},\quad
\bm r^x_{2,7}=\begin{pmatrix}0\\0\\\pm\beta_1\beta_3\\\mp\beta_1\beta_2\\\pm(\beta_3u_T-\beta_2w)\beta_1\\0\\\frac{1}{\sqrt{\rho}}\beta_3\\
-\frac{1}{\sqrt{\rho}}\beta_2\end{pmatrix},\quad
\bm r^x_{3,6}=\begin{pmatrix}\alpha_s\\\alpha_s(u_N\mp c_s)\\\alpha_su_T\mp\alpha_f\beta_1\beta_2c\\\alpha_sw\mp\alpha_f\beta_1\beta_3c\\
\hf\alpha_s|\bm u|^2+\mu_s^\mp\\0\\-\frac{1}{c_f\sqrt{\rho}}\alpha_f\beta_2c^2\\-\frac{1}{c_f\sqrt{\rho}}\alpha_f\beta_3c^2\end{pmatrix},\\
&\bm r^x_4=\Big(1,u_N,u_T,w,\hf|\bm u|^2,0,0,0\Big)^\top,\quad\bm r^x_5=(0,0,0,0,0,1,0,0)^\top,\quad\bm\ell^x_5=(0,0,0,0,0,1,0,0)^\top,\\
&\bm\ell^x_{1,8}=\begin{pmatrix}\hf\theta_1\alpha_fc^2|\bm u|^2\pm\theta_2\left[\alpha_fcu_N\beta_1-\alpha_sc_s\left(\beta_2u_T+\beta_3w
\right)\right]\\-\theta_1\alpha_fc^2u_N\mp\theta_2\alpha_fc\beta_1\\-\theta_1\alpha_fc^2u_T\pm\theta_2\alpha_sc_s\beta_2\\
-\theta_1\alpha_fc^2w\pm\theta_2\alpha_sc_s\beta_3\\\theta_1\alpha_fc^2\\0\\
\theta_1\sqrt{\rho}\alpha_s\beta_2c_f\left(c_s^2-\frac{\gamma_2}{\gamma_1}c^2\right)\\
\theta_1\sqrt{\rho}\alpha_s\beta_3c_f\left(c_s^2-\frac{\gamma_2}{\gamma_1}c^2\right)\end{pmatrix},\quad
\bm\ell^x_{2,7}=\begin{pmatrix}\mp\hf\beta_1\left(\beta_3u_T-\beta_2w\right)\\0\\\pm\hf\beta_1\beta_3\\\mp\hf\beta_1\beta_2\\0\\0\\
\hf\sqrt{\rho}\beta_3\\-\hf\sqrt{\rho}\beta_2\end{pmatrix},\\[1.0ex]
&\bm\ell^x_{3,6}=\begin{pmatrix}\hf\theta_1\alpha_sc_f^2|\bm u|^2\pm\theta_2\left[\alpha_sc_au_N\beta_1+\alpha_fc_f
\left(\beta_2u_T+\beta_3w\right)\right]\\-\theta_1\alpha_sc_f^2u_N\mp\theta_2\alpha_sc_a\beta_1\\
-\theta_1\alpha_sc_f^2u_T\mp\theta_2\alpha_fc_f\beta_2\\-\theta_1\alpha_sc_f^2w\mp\theta_2\alpha_fc_f\beta_3\\\theta_1\alpha_sc_f^2\\0\\
-\theta_1\sqrt{\rho}\alpha_f\beta_2c_f\left(c_f^2-\frac{\gamma_2}{\gamma_1}c^2\right)\\
-\theta_1\sqrt{\rho}\alpha_f\beta_3c_f\left(c_f^2-\frac{\gamma_2}{\gamma_1}c^2\right)\end{pmatrix},\quad
\bm\ell^x_4=\begin{pmatrix}1-\theta_1\left(\alpha_f^2c^2+\alpha_s^2c_f^2\right)|\bm u|^2\\
2\theta_1\left(\alpha_f^2c^2+\alpha_s^2c_f^2\right)u_N\\2\theta_1\left(\alpha_f^2c^2+\alpha_s^2c_f^2\right)u_T\\
2\theta_1\left(\alpha_f^2c^2+\alpha_s^2c_f^2\right)w\\-2\theta_1\left(\alpha_f^2c^2+\alpha_s^2c_f^2\right)\\0\\
2\theta_1\sqrt{\rho}\alpha_f\alpha_s\beta_2c_f\left(c_f^2-c_s^2\right)\\
2\theta_1\sqrt{\rho}\alpha_f\alpha_s\beta_3c_f\left(c_f^2-c_s^2\right)\end{pmatrix}.
\end{align*}
Here,
$$
\begin{aligned}
&\beta_1:={\rm sign}(b_N),\quad(\beta_2,\beta_3):=
\left\{\begin{aligned}&\bigg(\frac{1}{\sqrt{2}},\frac{1}{\sqrt{2}}\bigg),&&\mbox{if }b_T=b_3=0,\\
&\Bigg(\frac{b_T}{\sqrt{b_T^2+b_3^2}},\frac{b_3}{\sqrt{b_T^2+b_3^2}}\Bigg)&&\mbox{otherwise},\end{aligned}\right.\\[0.5ex]
&(\alpha_f,\alpha_s):=\left\{\begin{aligned}
&(1,1)&&\mbox{if }b_T=b_3=0,\\
&\left(\sqrt{\frac{(c_f^2-c_a^2)}{(c_f^2-c_s^2)}},\sqrt{\frac{(c_f^2-c^2)}{(c_f^2-c_s^2)}}\right)&&\mbox{otherwise},
\end{aligned}\right.\\[0.5ex] 
&\theta_1:=\hf\bigg[\alpha_f^2c^2\Big(c_f^2-\frac{\gamma_2}{\gamma_1}c^2\Big)+\alpha_s^2c_f^2\Big(c_s^2-\frac{\gamma_2}{\gamma_1}c^2\Big)
\bigg]^{-1},\quad\theta_2:=\hf\Big[\alpha_f^2c_fa\beta_1+\alpha_s^2c_sc_a\beta_1\Big]^{-1},\\
&\mu_f^\mp:=-\frac{\alpha_fc_f^2}{\gamma_1}\mp\alpha_fc_fu_N\pm\alpha_sc_a\beta_1(\beta_2u_T+\beta_3w)+
\frac{\gamma_2}{\gamma_1}\alpha_f(c_f^2-c^2),\\
&\mu_s^\mp:=-\frac{\alpha_sc_s^2}{\gamma_1}\mp\alpha_sc_su_N\mp\alpha_fc\beta_1(\beta_2u_T+\beta_3w)+
\frac{\gamma_2}{\gamma_1}\alpha_s(c_s^2-c^2).
\end{aligned}
$$

The structure of the matrices $R^y$ and $(R^y)^{-1}$ is similar, but in all of the right and left eigenvectors above, one needs to switch
the second and third components as well as the sixth and seventh components.

\section{Eigendecomposition for Primitive Variables}\label{appA2}
In this appendix, we provide the reader with the matrices used in the LCD for the quasi-linear system \eref{2.8}; see \cite{Balsara96} for
details.
 
First, the matrices $D^x$ and $D^y$ are 
\begin{equation*}
D^x(\mV)=\begin{pmatrix}u&\rho&0&0&0&0&0&0\\0&u&0&0&\dfrac{1}{\rho}&0&\dfrac{b_2}{\rho}&\dfrac{b_3}{\rho}\\[1.5ex]
0&0&u&0&0&0&-\dfrac{b_1}{\rho}&0\\0&0&0&u&0&0&0&-\dfrac{b_1}{\rho}\\0&\gamma p&0&0&u&0&0&0\\0&0&0&0&0&u&0&0\\0&b_2&-b_1&0&0&0&u&0\\
0&b_3&0&-b_1&0&0&0&u\end{pmatrix}
\end{equation*}
and 
\begin{equation*}
D^y(\mV)=\begin{pmatrix}v&0&\rho&0&0&0&0&0\\0&v&0&0&0&-\dfrac{b_2}{\rho}&0&0\\[1.5ex]
0&0&v&0&\dfrac{1}{\rho}&\dfrac{b_1}{\rho}&0&\dfrac{b_3}{\rho}\\0&0&0&v&0&0&0&-\dfrac{b_2}{\rho}\\0&0&\gamma p&0&v&0&0&0\\
0&-b_2&b_1&0&0&v&0&0\\0&0&0&0&0&0&v&0\\0&0&b_3&-b_2&0&0&0&v\end{pmatrix}.
\end{equation*}
We note that the matrices $D^x$ and $D^y$ have the same eigenvalues \eref{1.3a}--\eref{A4} as the matrices $C^x$ and $C^y$, but different
eigenvectors. The matrix $T^x$, which diagonalizes $D^x$ is $T^x=\big[\bm r^x_1|\bm r^x_2|\dots|\bm r^x_7|\bm r^x_8\big]$ and its inverse is
$(T^x)^{-1}=\big[\bm\ell^x_1|\bm\ell^x_2|\dots|\bm\ell^x_7|\bm\ell^x_8\big]^\top$, where $\bm r^x_i$ and $\bm \ell^x_i$, $i=1,\dots,8$ are
the right and left eigenvectors of $D^x$:
\allowdisplaybreaks
\begin{align*}
&\bm r^x_{1,8}=\begin{pmatrix}\widehat\alpha_f\rho\\\mp\widehat\alpha_fc_f\\\pm\widehat\alpha_sc_s\beta_1\beta_2\\
\pm\widehat\alpha_sc_s\beta_1\beta_3\\\widehat\alpha_f\rho c^2\\0\\\widehat\alpha_s\sqrt{\rho}c\beta_2\\\widehat\alpha_s\sqrt{\rho}c\beta_3
\end{pmatrix},\quad
\bm r^x_{2,7}=\begin{pmatrix}0\\0\\\mp\beta_3\\\pm\beta_2\\0\\0\\-\sqrt{\rho}\beta_1\beta_3\\\sqrt{\rho}\beta_1\beta_2\end{pmatrix},\quad
\bm r^x_{3,6}=\begin{pmatrix}\widehat\alpha_s\rho\\\mp\widehat\alpha_sc_s\\\mp\widehat\alpha_fc_f\beta_1\beta_2\\
\mp\widehat\alpha_fc_f\beta_1\beta_3\\\widehat\alpha_s\rho c^2\\0\\-\widehat\alpha_f\sqrt{\rho}a\beta_2\\
-\widehat\alpha_f\sqrt{\rho}a\beta_3 \end{pmatrix},\\
&\bm r^x_4=(1,0,0,0,0,0,0,0)^\top,\quad\bm r^x_5=(0,0,0,0,0,1,0,0)^\top,\\
&\bm \ell^x_{1,8}=\frac{1}{2c^2}\begin{pmatrix}0\\\mp\widehat\alpha_fc_f\\\pm\widehat\alpha_sc_s\beta_1\beta_2\\
\pm\widehat\alpha_sc_s\beta_1\beta_3\\\frac{1}{\rho}\widehat\alpha_f\\0\\\frac{1}{\sqrt{\rho}}\widehat\alpha_sc\beta_2\\
\frac{1}{\sqrt{\rho}}\widehat\alpha_sc\beta_3\end{pmatrix},\quad
\bm\ell^x_{2,7}=\hf\begin{pmatrix}0\\0\\\mp\beta_3\\\pm\beta_2\\0\\0\\-\frac{1}{\sqrt{\rho}}\beta_1\beta_3\\
\frac{1}{\sqrt{\rho}}\beta_1\beta_2\end{pmatrix},\quad
\bm \ell^x_{3,6}=\frac{1}{2c^2}\begin{pmatrix}0\\\mp\widehat\alpha_sc_s\\\mp\widehat\alpha_fc_f\beta_1\beta_2\\
\mp\widehat\alpha_fc_f\beta_1\beta_3\\\frac{1}{\rho}\widehat\alpha_s\\0\\-\frac{1}{\sqrt{\rho}}\widehat\alpha_fc\beta_2\\
-\frac{1}{\sqrt{\rho}}\widehat\alpha_fc\beta_3\end{pmatrix},\\
&\bm\ell^x_4=\Big(1,0,0,0,-\frac{1}{c^2},0,0,0\Big)^\top,\quad\bm\ell^x_5=(0,0,0,0,0,1,0,0)^\top.
\end{align*}
Here,
$$
(\widehat\alpha_f,\widehat\alpha_s)=\left\{\begin{aligned}
&\bigg(\frac{1}{\sqrt{2}},\frac{1}{\sqrt{2}}\bigg)&&\mbox{if }b_T=b_3=0,\\
&\left(\sqrt{\frac{(c^2-c_s^2)}{(c_f^2-c_s^2)}},\sqrt{\frac{(c_f^2-c^2)}{(c_f^2-c_s^2)}}\right)&&\mbox{otherwise},
\end{aligned}\right.
$$
and the other notations are the same as in Appendix \ref{appA1}. The structure of the matrices $T^y$ and $(T^y)^{-1}$, which diagonalize
$D^y$ is similar, but in all of the right and left eigenvectors above, one needs to switch the second and third components as well as the
sixth and seventh components.

\bibliographystyle{siamplain}
\bibliography{CKLN}

@article{BS1999,
     TITLE = {A staggered mesh algorithm using high order {G}odunov fluxes to ensure solenoidal magnetic fields in magnetohydrodynamic
simulations},
    AUTHOR = {Balsara, D. S. and Spicer, D. S.},
   JOURNAL = {J. Comput. Phys.},
    VOLUME = {149},
    NUMBER = {2},
     PAGES = {270--292},
      YEAR = {1999}}

@article{BB1980,
  TITLE={The effect of nonzero {$\nabla\cdot \bm B$} on the numerical solution of the magnetohydrodynamic equations},
  AUTHOR={Brackbill, J. U. and Barnes, D. C.},
  JOURNAL={J. Comput. Phys.},
  VOLUME={35},
  NUMBER={3},
  PAGES={426--430},
  YEAR={1980}}

@article{LS2005,
    AUTHOR = {Li, F. and Shu, C.-W.},
     TITLE = {Locally divergence-free discontinuous {G}alerkin methods for {MHD} equations},
   JOURNAL = {J. Sci. Comput.},
    VOLUME = {22/23},
      YEAR = {2005},
     PAGES = {413--442}}

@article{Toth2000,
  TITLE={The {$\nabla \cdot B = 0$} constraint in shock-capturing magnetohydrodynamics codes},
  AUTHOR={T{\'o}th, G.},
  JOURNAL={J. Comput. Phys.},
  VOLUME={161},
  NUMBER={2},
  PAGES={605--652},
  YEAR={2000}}

@article{EH1988,
    AUTHOR = {Evans, C. R. and Hawley, J. F.},
     TITLE = {Simulation of magnetohydrodynamic flows: {A} constrained transport method},
   JOURNAL = {Astrophys. J.},
    VOLUME = {332},
      YEAR = {1988},
    NUMBER = {2},
     PAGES = {659--677}}

@article{DeVore1991,
    AUTHOR = {DeVore, C. R.},
     TITLE = {Flux-corrected transport techniques for multidimensional compressible magnetohydrodynamics},
   JOURNAL = {J. Comput. Phys.},
    VOLUME = {92},
      YEAR = {1991},
    NUMBER = {2},
     PAGES = {142--160}}

@article{GS2005,
    AUTHOR = {Gardiner, T. A. and Stone, J. M.},
     TITLE = {An unsplit {G}odunov method for ideal {MHD} via constrained transport},
   JOURNAL = {J. Comput. Phys.},
    VOLUME = {205},
      YEAR = {2005},
    NUMBER = {2},
     PAGES = {509--539}}

@article{XBD2016,
    AUTHOR = {Xu, Z. and Balsara, D. S. and Du, H.},
     TITLE = {Divergence-free {WENO} reconstruction-based finite volume scheme for solving ideal {MHD} equations on triangular meshes},
   JOURNAL = {Commun. Comput. Phys.},
    VOLUME = {19},
      YEAR = {2016},
    NUMBER = {4},
     PAGES = {841--880}}

@article{YXL2013,
    AUTHOR = {Yakovlev, S. and Xu, L. and  Li, F.},
     TITLE = {Locally divergence-free central discontinuous {G}alerkin methods for ideal {MHD} equations},
   JOURNAL = {J. Comput. Sci.},
    VOLUME = {4},
      YEAR = {2013},
    NUMBER = {4},
     PAGES = {80--91}}

@article{PRLGD1999,
    AUTHOR = {Powell, K. G. and Roe, P. L. and Linde, T. J. and Gombosi, T. I. and De Zeeuw, D. L.},
     TITLE = {A solution-adaptive upwind scheme for ideal magnetohydrodynamics},
   JOURNAL = {J. Comput. Phys.},
    VOLUME = {154},
      YEAR = {1999},
    NUMBER = {2},
     PAGES = {284--309}}

@article{Balsara2009,
    AUTHOR = {Balsara, D. S.},
     TITLE = {Divergence-free reconstruction of magnetic fields and {WENO} schemes for magnetohydrodynamics},
   JOURNAL = {J. Comput. Phys.},
    VOLUME = {228},
      YEAR = {2009},
    NUMBER = {14},
     PAGES = {5040--5056}}

@article{BKC2021,
    AUTHOR = {Balsara, D. S. and Kumar, R. and Chandrashekar, P.},
     TITLE = {Globally divergence-free {DG} scheme for ideal compressible {MHD}},
   JOURNAL = {Commun. Appl. Math. Comput. Sci.},
    VOLUME = {16},
      YEAR = {2021},
    NUMBER = {1},
     PAGES = {59--98}}

@article{DBTF2019,
    AUTHOR = {Dumbser, M. and Balsara, D. S. and Tavelli, M. and Fambri, F.},
     TITLE = {A divergence-free semi-implicit finite volume scheme for ideal, viscous, and resistive magnetohydrodynamics},
   JOURNAL = {Internat. J. Numer. Methods Fluids},
    VOLUME = {89},
      YEAR = {2019},
    NUMBER = {1-2},
     PAGES = {16--42}}

@article{FLX2018,
    AUTHOR = {Fu, P. and Li, F. and Xu, Y.},
     TITLE = {Globally divergence-free discontinuous {G}alerkin methods for ideal magnetohydrodynamic equations},
   JOURNAL = {J. Sci. Comput.},
    VOLUME = {77},
      YEAR = {2018},
    NUMBER = {3},
     PAGES = {1621--1659}}

@article{LX2012,
    AUTHOR = {Li, F. and Xu, L.},
     TITLE = {Arbitrary order exactly divergence-free central discontinuous {G}alerkin methods for ideal {MHD} equations},
   JOURNAL = {J. Comput. Phys.},
    VOLUME = {231},
      YEAR = {2012},
    NUMBER = {6},
     PAGES = {2655--2675}}

@inproceedings{PRMGD1995,
  TITLE={An upwind scheme for magnetohydrodynamics},
  AUTHOR={Powell, K. G. and Roe, P. L. and Myong, R. S. and Gombosi, T. and De Zeeuw, D.},
  BOOKTITLE={12th Computational Fluid Dynamics Conference: AIAA Paper 95-1704-CP},
  PAGES={661--674},
  YEAR={1995}}

@article{CRT2014,
    AUTHOR = {Christlieb, A. J. and Rossmanith, J. A. and Tang, Q.},
     TITLE = {Finite difference weighted essentially non-oscillatory schemes with constrained transport for ideal magnetohydrodynamics},
   JOURNAL = {J. Comput. Phys.},
    VOLUME = {268},
      YEAR = {2014},
     PAGES = {302--325}}

@article{HRT2013,
    AUTHOR = {Helzel, C. and Rossmanith, J. A. and Taetz, B.},
     TITLE = {A high-order unstaggered constrained-transport method for the three-dimensional ideal magnetohydrodynamic equations based on
              the method of lines},
   JOURNAL = {SIAM J. Sci. Comput.},
    VOLUME = {35},
      YEAR = {2013},
    NUMBER = {2},
     PAGES = {A623--A651}}

@article{MT2012,
    AUTHOR = {Mishra, S. and Tadmor, E.},
     TITLE = {Constraint preserving schemes using potential-based fluxes. {III}. {G}enuinely multi-dimensional schemes for {MHD} equations},
   JOURNAL = {ESAIM Math. Model. Numer. Anal.},
    VOLUME = {46},
      YEAR = {2012},
    NUMBER = {3},
     PAGES = {661--680}}

@article{Rossmanith2006,
    AUTHOR = {Rossmanith, J. A.},
     TITLE = {An unstaggered, high-resolution constrained transport method for magnetohydrodynamic flows},
   JOURNAL = {SIAM J. Sci. Comput.},
    VOLUME = {28},
      YEAR = {2006},
    NUMBER = {5},
     PAGES = {1766--1797}}

@article{CK2016,
    AUTHOR = {Chandrashekar, P. and Klingenberg, C.},
     TITLE = {Entropy stable finite volume scheme for ideal compressible {MHD} on 2-{D} {C}artesian meshes},
   JOURNAL = {SIAM J. Numer. Anal.},
    VOLUME = {54},
      YEAR = {2016},
    NUMBER = {2},
     PAGES = {1313--1340}}

@article{DGWW2018,
    AUTHOR = {Derigs, D. and  Gassner, G. J. and  Walch, S. and  Winters, A. R.},
     TITLE = {Entropy stable finite volume approximations for ideal magnetohydrodynamics},
   JOURNAL = {Jahresber. Dtsch. Math.-Ver.},
    VOLUME = {120},
      YEAR = {2018},
    NUMBER = {1},
     PAGES = {153--219}}

@article{Godunov2025,
  TITLE={Symmetric form of the equations of magnetohydrodynamics (in {R}ussian)},
  AUTHOR={Godunov, S. K.},
  JOURNAL={Numerical Methods for Mechanics of Continuum Medium},
  VOLUME={1},
  PAGES={26--34},
  YEAR={1972}}

@article{CCHKL_22,
   AUTHOR = {Chertock, A. and Chu, S. and Herty, M. and Kurganov, A. and Luk\'{a}\v{c}ov\'{a}-Medvi{\softd}ov\'{a}, M.},
  journal = {J. Comput. Phys.},
  title   = {Local characteristic decomposition based central-upwind scheme},
  year    = {2023},
  note    = {Paper No. 111718},
  volume  = {473}}

@article{CKRZ2024,
    AUTHOR = {Chertock, A. and Kurganov, A. and Redle, M. and Zeitlin, V.},
     TITLE = {Locally divergence-free well-balanced path-conservative central-upwind schemes for rotating shallow water {MHD}},
   JOURNAL = {J. Comput. Phys.},
    VOLUME = {518},
      YEAR = {2024},
      NOTE = {Paper No. 113300}}

@article{CKM,
    AUTHOR = {Castro D{\'{\i}}az, M. J. and Kurganov, A. and Morales de Luna, T.},
     TITLE = {Path-conservative central-upwind schemes for nonconservative hyperbolic systems},
   JOURNAL = {ESAIM Math. Model. Numer. Anal.},
    VOLUME = {53},
      YEAR = {2019},
    NUMBER = {3},
     PAGES = {959--985}}

@book{Gottlieb11,
     AUTHOR = {Gottlieb, S. and Ketcheson, D. and Shu, C.-W.},
     TITLE = {Strong stability preserving {R}unge-{K}utta and multistep time discretizations},
     PUBLISHER = {World Scientific Publishing Co. Pte. Ltd., Hackensack, NJ},
     YEAR = {2011},
     PAGES = {xii+176}}

@article {Gottlieb12,
	Author = {Gottlieb, S. and Shu, C.- W. and Tadmor, E.},
	Journal = {SIAM Rev.},
	Number = {1},
	Pages = {89--112},
	Title = {Strong stability-preserving high-order time discretization
	methods},
	Volume = {43},
	Year = {2001}}

@incollection{CK2023,
    AUTHOR = {Chu, S. and Kurganov, A.},
     TITLE = {Local characteristic decomposition based central-upwind scheme for compressible multifluids},
 BOOKTITLE = {Finite volumes for complex applications {X}. {V}ol. 2. {H}yperbolic and related problems},
    SERIES = {Springer Proc. Math. Stat.},
    VOLUME = {433},
     PAGES = {73--81},
 PUBLISHER = {Springer, Cham},
      YEAR = {2023}}

@article{CHK25,
  author  = {Chu, S. and Herty, M. and Kurganov, A.},
  journal = {J. Comput. Phys.},
  title   = {Novel local characteristic decomposition based path-conservative central-upwind schemes},
  year    = {2025},
  note    = {Paper No. 113692},
  volume  = {524}}

@article{LSZ2018,
    AUTHOR = {Liu, Y. and Shu, C.-W. and Zhang, M.},
     TITLE = {Entropy stable high order discontinuous {G}alerkin methods for ideal compressible {MHD} on structured meshes},
   JOURNAL = {J. Comput. Phys.},
    VOLUME = {354},
      YEAR = {2018},
     PAGES = {163--178}}

@article{LXY2011,
    AUTHOR = {Li, F. and Xu, L. and Yakovlev, S.},
     TITLE = {Central discontinuous {G}alerkin methods for ideal {MHD} equations with the exactly divergence-free magnetic field},
   JOURNAL = {J. Comput. Phys.},
    VOLUME = {230},
      YEAR = {2011},
    NUMBER = {12},
     PAGES = {4828--4847}}

@incollection{Powell1997,
author={Powell, K. G.},
editor={Hussaini, M. Y. and van Leer, B. and Van Rosendale, J.},
title={An approximate {R}iemann solver for magnetohydrodynamics},
bookTitle={Upwind and High-Resolution Schemes},
year={1997},
publisher={Springer Berlin Heidelberg},
address={Berlin, Heidelberg},
pages={570--583},
isbn={978-3-642-60543-7}}

@article{Chertock24,
    AUTHOR = {Chertock, A. and Kurganov, A. and Redle, M. and Wu, K.},
     TITLE = {A new locally divergence-free path-conservative central-upwind
              scheme for ideal and shallow water magnetohydrodynamics},
   JOURNAL = {SIAM J. Sci. Comput.},
    VOLUME = {46},
      YEAR = {2024},
    NUMBER = {3},
     PAGES = {A1998--A2024}}

@article {Balsara96,
    AUTHOR = {Roe, P. L. and Balsara, D. S.},
     TITLE = {Notes on the eigensystem of magnetohydrodynamics},
   JOURNAL = {SIAM J. Appl. Math.},
  FJOURNAL = {SIAM Journal on Applied Mathematics},
    VOLUME = {56},
      YEAR = {1996},
    NUMBER = {1},
     PAGES = {57--67},}

@article {Brio88,
    AUTHOR = {Brio, M. and Wu, C. C.},
     TITLE = {An upwind differencing scheme for the equations of ideal
              magnetohydrodynamics},
   JOURNAL = {J. Comput. Phys.},
  FJOURNAL = {Journal of Computational Physics},
    VOLUME = {75},
      YEAR = {1988},
    NUMBER = {2},
     PAGES = {400--422},}

@article{Orszag_Tang_1979,
author={Orszag, S. A. and Tang, C.-M.},
title={Small-scale structure of two-dimensional magnetohydrodynamic turbulence},
volume={90},
DOI={10.1017/S002211207900210X},
number={1},
journal={J. Fluid Mech.},
year={1979},
pages={129--143}}

@article{KNP,
    AUTHOR = {Kurganov, A. and Noelle, S. and Petrova, G.},
     TITLE = {Semidiscrete central-upwind schemes for hyperbolic conservation laws and {H}amilton-{J}acobi equations},
   JOURNAL = {SIAM J. Sci. Comput.},
    VOLUME = {23},
      YEAR = {2001},
    NUMBER = {3},
     PAGES = {707--740}}
\end{document}